\documentclass{article}

\usepackage[margin=1in]{geometry}
\usepackage{amsmath, mathtools, amsfonts, amsthm, amssymb, cite, suffix, enumitem, bm, graphbox}
\mathtoolsset{showonlyrefs}
\usepackage[colorlinks]{hyperref}
\hypersetup{linkcolor=blue}
\usepackage{xcolor}

\usepackage{multirow}
\usepackage{pgfplots}
\pgfplotsset{compat=1.18}

\newtheorem{theorem}{Theorem}[section]
\newtheorem{lemma}[theorem]{Lemma}
\newtheorem{corollary}[theorem]{Corollary}
\newtheorem{assumption}[theorem]{Assumption}
\newtheorem{proposition}[theorem]{Proposition}

\theoremstyle{definition}

\newtheorem{remark}[theorem]{Remark}

\usepackage{mleftright}
\renewcommand{\left}{\mleft}
\renewcommand{\right}{\mright}

\renewcommand{\emptyset}{\varnothing}

\renewcommand{\P}{\mathbb{P}}
\newcommand{\E}{\mathbb{E}}
\newcommand{\R}{\mathbb{R}}
\newcommand{\N}{\mathbb{N}}
\renewcommand{\L}{\mathcal{L}}
\newcommand{\cL}{\mathcal{L}}

\newcommand{\cG}{\mathcal{G}}
\newcommand{\cB}{\mathcal{B}}
\newcommand{\cA}{\mathcal{A}}
\newcommand{\cC}{\mathcal{C}}
\newcommand{\cM}{\mathcal{M}}
\newcommand{\cN}{\mathcal{N}}
\newcommand{\cV}{\mathcal{V}}

\renewcommand{\c}{\mathsf{c}}

\newcommand{\cE}{\mathcal{E}}

\newcommand{\Var}{\operatorname{Var}}

\newcommand{\Cov}{\operatorname{Cov}}
\newcommand{\TV}{\mathrm{TV}}
\renewcommand{\d}{\mathbf{d}}
\newcommand{\ind}[1]{\mathbf{1}_{\{#1\}}}

\newcommand{\fr}{\frac}
\newcommand{\eps}{\varepsilon}
\newcommand{\Exp}[1]{\exp\left(#1\right)}
\newcommand{\Econd}[2]{\E\left[ #1 \,\middle|\, #2 \right]}
\newcommand{\Pcond}[2]{\P\left[ #1 \,\middle|\, #2 \right]}
\WithSuffix\newcommand\ind*[1]{\mathbf{1}_{#1}}
\newcommand{\sse}{\subseteq}

\newcommand{\polylog}{\operatorname{polylog}}
\newcommand{\gampe}{\Gamma_p^\eps}
\newcommand{\condmeas}{\mu_\sH^\eta(n,p)}

\newcommand{\pball}{\cB_p^\eta}
\newcommand{\phalfball}{\cB_p^{\eta/2}}

\newcommand{\dtv}{\d_\TV}
\renewcommand{\dh}{\d_\mathrm{Ham}}
\newcommand{\dK}{\d_\mathrm{Kol}}
\newcommand{\dW}{\d_\mathrm{Was}}

\newcommand{\dls}{\d_{\Lambda^*}}

\newcommand{\db}{\d_\square}
\newcommand{\dloc}{\dh^{\mathrm{loc}(v)}}

\newcommand{\pto}{\overset{\P}{\longrightarrow}}

\newcommand{\sg}{\mathsf{g}}
\newcommand{\sr}{\mathsf{r}}
\newcommand{\st}{\mathsf{t}}
\newcommand{\sv}{\mathsf{v}}
\newcommand{\se}{\mathsf{e}}
\newcommand{\sd}{\mathsf{d}}
\newcommand{\sE}{\mathsf{E}}
\newcommand{\sN}{\mathsf{N}}
\newcommand{\hsN}{\widehat{\sN}}
\newcommand{\sH}{\mathsf{H}}
\newcommand{\sR}{\mathsf{R}}
\renewcommand{\deg}{\mathsf{deg}}

\renewcommand{\tilde}{\widetilde}
\newcommand{\edgeset}{\binom{[n]}{2}}

\newcommand{\X}{\widetilde{X}}

\let\temp\phi
\let\phi\varphi
\let\varphi\temp

\pgfmathsetmacro{\betazero}{-1}
\pgfmathsetmacro{\betaone}{0.53}
\pgfmathsetmacro{\betatwo}{0.5}
\pgfmathsetmacro{\pstar}{0.17865}

\pgfmathdeclarefunction{H}{1}{% one argument
  \pgfmathparse{\betazero * #1 + \betaone * #1^2 + \betatwo * #1^3}%
}

\pgfmathdeclarefunction{I}{1}{% one argument
  \pgfmathparse{0.5 * #1 * ln(#1) + 0.5 * (1 - #1) * ln(1 - #1)}%
}

\pgfmathdeclarefunction{L}{1}{% one argument
  \pgfmathparse{H(#1) - I(#1)}%
}

\pgfmathdeclarefunction{Hprime}{1}{% one argument
  \pgfmathparse{\betazero + 2 * \betaone * #1 + 3 * \betatwo * #1^2}%
}

\pgfmathdeclarefunction{phi}{1}{% one argument
  \pgfmathparse{exp(2 * Hprime(#1)) / (1 + exp(2 * Hprime(#1)))}%
}

\title{Quantitative central limit theorems for exponential random graphs}
\author{Vilas Winstein \\ University of California, Berkeley \\ \texttt{vilas@berkeley.edu}}
\begin{document}

\maketitle

\begin{abstract}
Ferromagnetic exponential random graph models (ERGMs) are nonlinear exponential tilts of Erd\H{o}s--R\'enyi
models, under which the presence of certain subgraphs such as triangles may be emphasized.
These models are mixtures of \emph{metastable wells} which each behave \emph{macroscopically} like new
Erd\H{o}s--R\'enyi models themselves, exhibiting the same laws of large numbers for the overall \emph{edge count}
as well as all \emph{subgraph counts}.
However, the microscopic fluctuations of these quantities had remained elusive for some time.
A recent breakthrough of Fang, Liu, Shao, and Zhao \cite{fang2024normal} introduces a new application of
\emph{Stein's method} to nonlinear exponential tilts of product measures, and uses it to prove \emph{quantitative}
central limit theorems (CLTs) for edge and subgraph counts under ferromagnetic ERGMs in the \emph{Dobrushin regime}
of parameters, which is a \emph{perturbative} or \emph{very high-temperature} regime.
A key input in their proof is the fact that the fluctuations of subgraph counts are driven by those
of the overall edge count, which was first proved by Sambale and Sinulis \cite{sambale2020logarithmic}
for triangle counts in the Dobrushin regime using functional-analytic tools such as modified log-Sobolev inequalities.
Additionally, during the preparation of the present article, \cite{fang2024normal} was updated to include a separate
argument which covers the full \emph{subcritical}, or \emph{high-temperature} regime of parameters, albeit
with worse error bounds than in the Dobrushin regime.

In the present article, we develop a novel probabilistic technique, based on the careful
analysis of the evolution of relevant quantities under the \emph{Glauber dynamics}, which we feel
clarifies the underlying mechanisms at play.
This technique allows us to improve upon the result of \cite{sambale2020logarithmic}
in multiple ways, including extending its validity beyond the Dobrushin
regime to all parameter regimes including the \emph{supercritical} or \emph{low-temperature} regime.
As a consequence we can extend the quantitative CLT of \cite{fang2024normal} to the supercritical regime,
and we also obtain improved error bounds in the subcritical regime (beyond the Dobrushin regime).
Furthermore, our technique is flexible enough that we can also prove quantitative CLTs for \emph{vertex degrees}
and \emph{local subgraph counts}, e.g.\ the number of triangles at a particular vertex.
To the best of our knowledge, these local results have not appeared before in \emph{any} parameter regime.
\end{abstract}

\setcounter{tocdepth}{2}
\tableofcontents

\section{Introduction}
\label{sec:intro}

\emph{Exponential random graph models} (ERGMs) are exponential families of random graph models.
In other words, the $n$-vertex random graph $X$ has an ERGM distribution if there is a function
$\sH$, called the \emph{Hamiltonian}, such that
\begin{equation}
\label{eq:gibbsmeasure}
    \P[X = x] \propto e^{n^2 \sH(x)}
\end{equation}
for any fixed simple $n$-vertex graph $x$.

Typically, the Hamiltonian $\sH$ is required to be a continuous function on the space of \emph{graphons},
which will be introduced properly in Section \ref{sec:review_leadingorder_graphons} below.
This means that $\sH(x)$ is of order $1$ for typical dense graphs, and in \eqref{eq:gibbsmeasure}
we scale $\sH$ by $n^2$ so that we do not get a trivial measure.
Indeed, note that the number of graphs on $n$ vertices is itself exponential in $n^2$, so scaling $\sH$
by something $\ll n^2$ would lead to an approximately uniform distribution, whereas scaling $\sH$ by
something $\gg n^2$ would lead to a distribution concentrated on the maximizer(s) of $\sH$.
When properly scaled, Gibbs measures (i.e.\ exponential families)
such as \eqref{eq:gibbsmeasure} interpolate naturally between these two extremes.

ERGMs have been examined by various practitioners, initially inspired by social networks
\cite{frank1986markov,holland1981exponential,wasserman1994social}, and also studied early on
from a statistical physics perspective
\cite{burda2004network,park2004solution,park2005solution}.
Since then ERGMs have attracted the attention of probabilists, statisticians,
and computer scientists who have investigated them from a variety of angles including
large deviations and phase transitions \cite{chatterjee2013estimating,radin2013phaseComplex,radin2013phaseERGM,yin2013critical,eldan2018exponential},
distributional limit theorems \cite{mukherjee2023statistics,fang2024normal,fang2025conditionalcentrallimittheorems},
sampling \cite{bhamidi2008mixing,reinert2019approximating,bresler2024metastable}, and
concentration of measure \cite{sambale2020logarithmic,ganguly2024sub,winstein2025concentration}.
Certain ERGMs are also useful, via ideas such as importance sampling, for studying random graphs constrained to
have an atypical number of triangles or other subgraphs \cite{radin2014asymptotics,lubetzky2015replica}.
Since we consider Hamiltonians $\sH$ which are naturally thought of as functions on the space of graphons,
i.e.\ \emph{dense} graph limits, typical samples are dense graphs; some recent work has also studied
modifications which yield sparse graphs instead \cite{yin2017asymptotics,cook2024typical}.
Needless to say, this subject is quite vast and there are many important works not
mentioned in this list.

In the present article, we consider Hamiltonians $\sH$ which can be written as finite linear combinations
of \emph{subgraph densities}, and moreover we restrict ourselves to the \emph{ferromagnetic regime},
where the coefficients of these densities are nonnegative, apart from the coefficient of the overall
\emph{edge} density.
This is a positive-correlation condition, and it turns out to result in homogeneous large-scale behavior,
meaning that edge and subgraph counts satisfy the same \emph{laws of large numbers} as the corresponding counts
in an Erd\H{o}s--R\'enyi graph $\cG(n,p)$, for some value of $p$ determined by the Hamiltonian.

However, such large-scale results leave open the question of the microscopic fluctuations of these counts;
these fluctuations are the focus of the present article.
As it turns out, many observables of ERGMs such as edge counts, subgraph counts, degrees, and more, have
approximately Gaussian fluctuations, meaning they obey \emph{central limit theorems} (CLTs).
Moreover, the distance between the distributions of these observables and a Gaussian may be quantified using
various natural metrics on the space of probability distributions.

Such distributional limit results have already appeared in a few contexts, as will be discussed in more detail
in Section \ref{sec:intro_setup_clt} below.
Our approach follows that of \cite{fang2024normal}, who proved the first \emph{quantitative} CLT for the edge
count $\sE(X)$ of an ERGM sample $X$ in the \emph{subcritical} or \emph{high-temperature} regime of parameters.
This implies a similar CLT for subgraph counts, as observed first by \cite{sambale2020logarithmic}
in the \emph{Dobrushin} or \emph{very-high-temperature} regime; moreover, \cite{fang2024normal}
makes use of the result of \cite{sambale2020logarithmic} in proving the CLT for $\sE(X)$ itself.
Almost all previous work on this subject has been restricted to high temperature, or in many cases \emph{very} high temperature,
and the \emph{supercritical} or \emph{low-temperature} regime has been relatively untouched.
\begin{center}
    \textbf{There are three main contributions of the present article:}
\end{center}
\begin{enumerate}
    \item We extend the quantitative CLT for the edge count $\sE(X)$ of \cite{fang2024normal}, as well as the
    resulting CLTs for subgraph counts due originally to \cite{sambale2020logarithmic} in the triangle case,
    to the supercritical regime.
    \item We improve the quantitative error bounds on the distance between the distribution of $\sE(X)$ and
    a Gaussian given by \cite{fang2024normal} in the subcritical regime (outside of the Dobrushin regime).
    \item We prove a quantitative CLT for the degree $\deg_v(X)$ of a vertex under the ERGM measure, as well
    as local versions of subgraph counts such as the number of triangles at a vertex, in all parameter regimes.
    To the best of our knowledge, even a \emph{non-quantitative} CLT for these observables has not appeared before in
    \emph{any} parameter regime, except for in the special case of the $2$-star ERGM \cite{mukherjee2023statistics}.
\end{enumerate}
A precise definition of the parameter regimes mentioned above will be provided below in
Section \ref{sec:review_fluctuations_dynamical}; additionally, a numerical plot of these regimes is given
in Figure \ref{fig:regimes} for a particular example of an ERGM.

\subsection{Problem setup and related works}
\label{sec:intro_setup}

Let us now precisely define the ERGMs we consider, and give a brief overview of previous work
relevant to the CLTs we will discuss.

\subsubsection{Definition of the Hamiltonian}
\label{sec:intro_setup_ergm}

Given a fixed graph $G = (\cV, \cE)$ and an $n$-vertex graph $x$,
we denote by $\sN_G(x)$ the number of \emph{labeled copies of $G$ in $x$}, also called the
\emph{homomorphism count}; this is the number of maps $\cV \to [n]$ (the vertex set of $x$) which map
edges to edges.
Importantly, these homomorphisms are \emph{not induced}, meaning we do not require that non-edges map
to non-edges.
Moreover, we \emph{do not} require that the homomorphisms be injective; this is simply to be consistent with
relevant previous works such as \cite{chatterjee2013estimating,winstein2025concentration}.
We remark that \cite{fang2024normal} did assume injectivity, but this assumption is mild in both directions
and our non-assumption of injectivity will only lead to a small modification in the proof.

Note that $\sN_G(x)$ is also what we mean by \emph{subgraph count}, i.e.\ the term ``subgraph''
and ``homomorphism'' are synonymous for the purposes of the present article.
This may have some counterintuitive implications, as in $\sN_G(x)$, symmetries are included; for instance,
if $G$ consists of a single edge then $\sN_G(x) = 2 \sE(x)$, and the homomorphism count of a triangle
would be multiplied by $6$ as compared to what one might typically think of as the triangle count.
In general, the multiplier is the number of \emph{automorphisms} of the graph $G$.

In any case, we define $\st(G,x)$ to be the \emph{homomorphism density}
\begin{equation}
\label{eq:homdensity}
    \st(G,x) \coloneqq \fr{\sN_G(x)}{n^{|\cV|}},
\end{equation}
which is the probability that a uniformly selected map $\cV \to [n]$ is a homomorphism.
This notation comes from the theory of graphons, which will be discussed in Section
\ref{sec:review_leadingorder_graphons} below.

With the definition \eqref{eq:homdensity} in hand, we may define the Hamiltonian of a ferromagnetic ERGM.
Let us fix a sequence of graphs $G_0, G_1, \dotsc, G_K$, where we always take $G_0$ to be a single edge.
Then, for any parameters $\beta_0 \in \R$ and $\beta_1, \dotsc, \beta_K \geq 0$, let us define the Hamiltonian
\begin{equation}
\label{eq:hdef}
    \sH(x) = \sH_\beta(x) \coloneqq \sum_{j=0}^K \beta_j \, \st(G_j, x).
\end{equation}
Here $\beta$ is shorthand for the vector of parameters $\beta = (\beta_0, \beta_1, \dotsc,\beta_K)$.
The condition that $\beta_j \geq 0$ for all $j > 0$ defines the so-called \emph{ferromagnetic} regime of
parameters; note that we still allow $\beta_0 < 0$ in this regime.
For the moment, let us denote by $\X$ an $n$-vertex sample from the ERGM measure \eqref{eq:gibbsmeasure} with this
choice \eqref{eq:hdef} of Hamiltonian $\sH$; we will shortly replace this with $X$ from the true distribution
we are after, which is \emph{conditioned} on a metastable well (to be defined shortly).

Note that if $K = 0$ we recover the Erd\H{o}s-R\'enyi model $\cG(n,q)$ with
$q = \fr{e^{2\beta_0}}{1 + e^{2\beta_0}}$.
So every ERGM we consider is a \emph{nonlinear exponential tilt} of such an Erd\H{o}s-R\'enyi model.
For instance, if $K=1$ and $G_1$ is a triangle, we recover the \emph{edge-triangle}
model, where the number of triangles is increased from $\cG(n,q)$ via tilting.

\subsubsection{Macroscopic behavior}
\label{sec:intro_setup_macroscopic}

It was shown in \cite{chatterjee2013estimating} via a large deviations principle that,
apart from a Lebesgue-measure-zero set of parameters
$\beta \in \R \times \R_{\geq0}^K$, the ferromagnetic condition implies the existence of some $p$
for which
\begin{equation}
\label{eq:macroscopicp}
    \st(G,\X) \pto p^{|\cE|}
\end{equation}
for every fixed graph $G = (\cV, \cE)$, as $n \to \infty$.
Note that in general this $p$ is \emph{not the same} as $\fr{e^{2\beta_0}}{1 + e^{2\beta_0}}$, the density of
the base Erd\H{o}s-R\'enyi model away from which the ERGM was tilted.

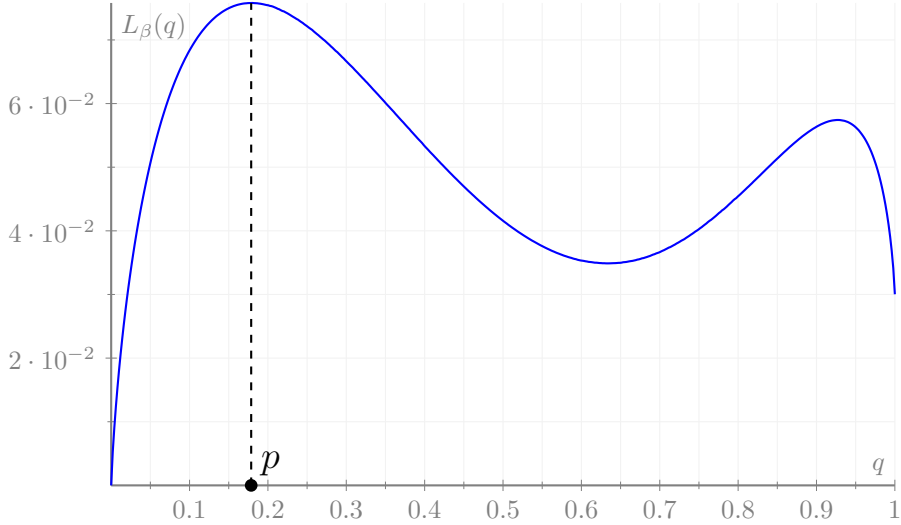
\begin{figure}
\centering
\begin{tikzpicture}
  \begin{axis}[
    width=12cm,
    height=8cm,
    axis lines=middle,
    xlabel={$q$},
    ylabel={$L_\beta(q)$},
    domain=0:1,
    samples=500,
    thick,
    grid=both,
    minor tick num=1,
    grid style={gray!10},
    enlargelimits=false,
    clip=false,
    scaled ticks=false,
    axis line style={-,gray},
    label style={gray},
    tick label style={gray},
  ]

  \addplot[blue, thick] {L(x)};

  \addplot[dashed, black] coordinates {(\pstar,0) (\pstar,{L(\pstar)})};
  \addplot[mark=*, black] coordinates {(\pstar,0)};

  \node[above right, black] at (axis cs:\pstar,0) {\Large $p$};
  \end{axis}
\end{tikzpicture}
\caption{Plot of $L_\beta(q)$, as defined in \eqref{eq:Ldef}, for the ERGM with $G_0, G_1, G_2$
an edge, a wedge ($2$-star), and a triangle respectively, with parameters
$\beta = (\beta_0,\beta_1,\beta_2) = (-1, 0.53, 0.5)$.
The unique optimal density $p$, which maximizes $L_\beta$ globally, is highlighted.
There are multiple local maximizers, so this choice of $\beta$ lies in the \emph{supercritical}
or \emph{low-temperature} regime for this ERGM, which will be introduced
in Section \ref{sec:review_fluctuations_dynamical}.}
\label{fig:lbeta}
\end{figure}

Now, if $Y$ is a sample from the Erd\H{o}s--R\'enyi model $\cG(n,p)$, then we also have
$\st(G,Y) \pto p^{|\cE|}$; therefore, \eqref{eq:macroscopicp}
is the statement that the ERGM sample $\X$ behaves like the Erd\H{o}s--R\'enyi sample $Y$, in a certain
\emph{macroscopic} sense: both models have the same laws of large numbers for all subgraph counts.
This will be further elucidated in Section \ref{sec:review_leadingorder_ldp}.
As will also be explained in that section,
the density $p$ in \eqref{eq:macroscopicp} is exhibited as the global maximizer in $[0,1]$ of the function
\begin{equation}
\label{eq:Ldef}
    L_\beta(q) \coloneqq \sH(q) - I(q),
\end{equation}
where $I(q) \coloneqq \fr{1}{2} q \log q + \fr{1}{2} (1-q) \log(1-q)$, and we evaluate the Hamiltonian
$\sH$ at a number $q \in [0,1]$ by interpreting $\st(G,q)$ as $q^{|\cE|}$.
See Figure \ref{fig:lbeta} for a sample plot of $L_\beta(q)$.

The Lebesgue-measure-zero set of exceptional parameters $\beta$ are those for which $L_\beta$ has
\emph{multiple} global maxima; let us denote by $M_\beta$ the (finite) set of global maximizers of $L_\beta$. 
If $|M_\beta| > 1$, we may decompose the ERGM measure into a mixture of measures, one for each
$p \in M_\beta$, under which \eqref{eq:macroscopicp} holds.
In this case, we say that there is \emph{phase coexistence}, while if $|M_\beta| = 1$ then there is
\emph{phase uniqueness}.
From now on, we let $X$ denote a sample from \emph{one of the measures} in this mixture, so that
there is a unique $p$ for which $X$ looks, macroscopically, like a sample $Y$ from the
Erd\H{o}s--R\'enyi model $\cG(n,p)$.

\subsubsection{Approximate independence}
\label{sec:intro_setup_indep}

\begin{figure}
\centering
$\underset{\text{Dobrushin}}{\includegraphics[scale=0.18,align=c]{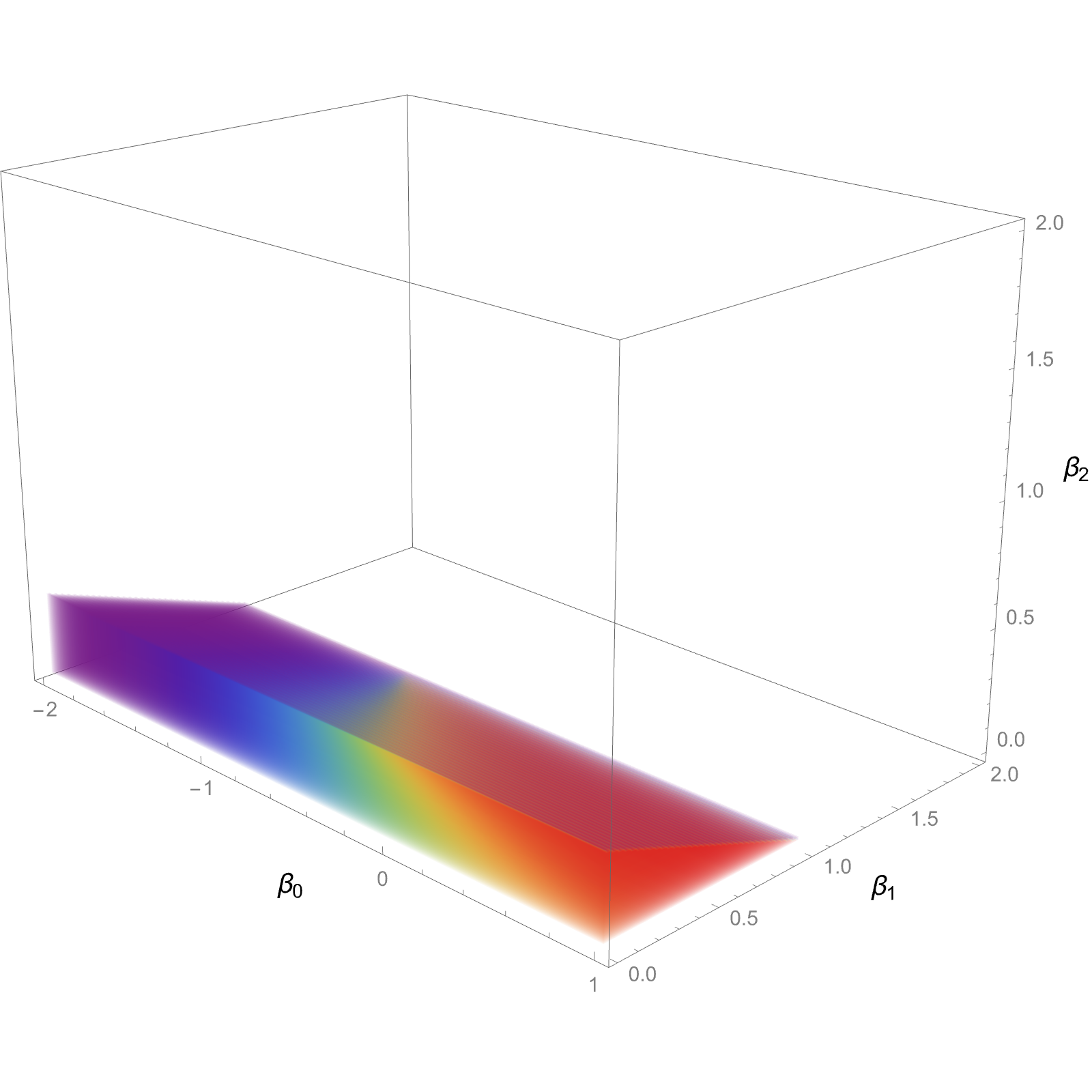}}$
\hspace{0.1cm}
$\underset{\text{Subcritical}}{\includegraphics[scale=0.18,align=c]{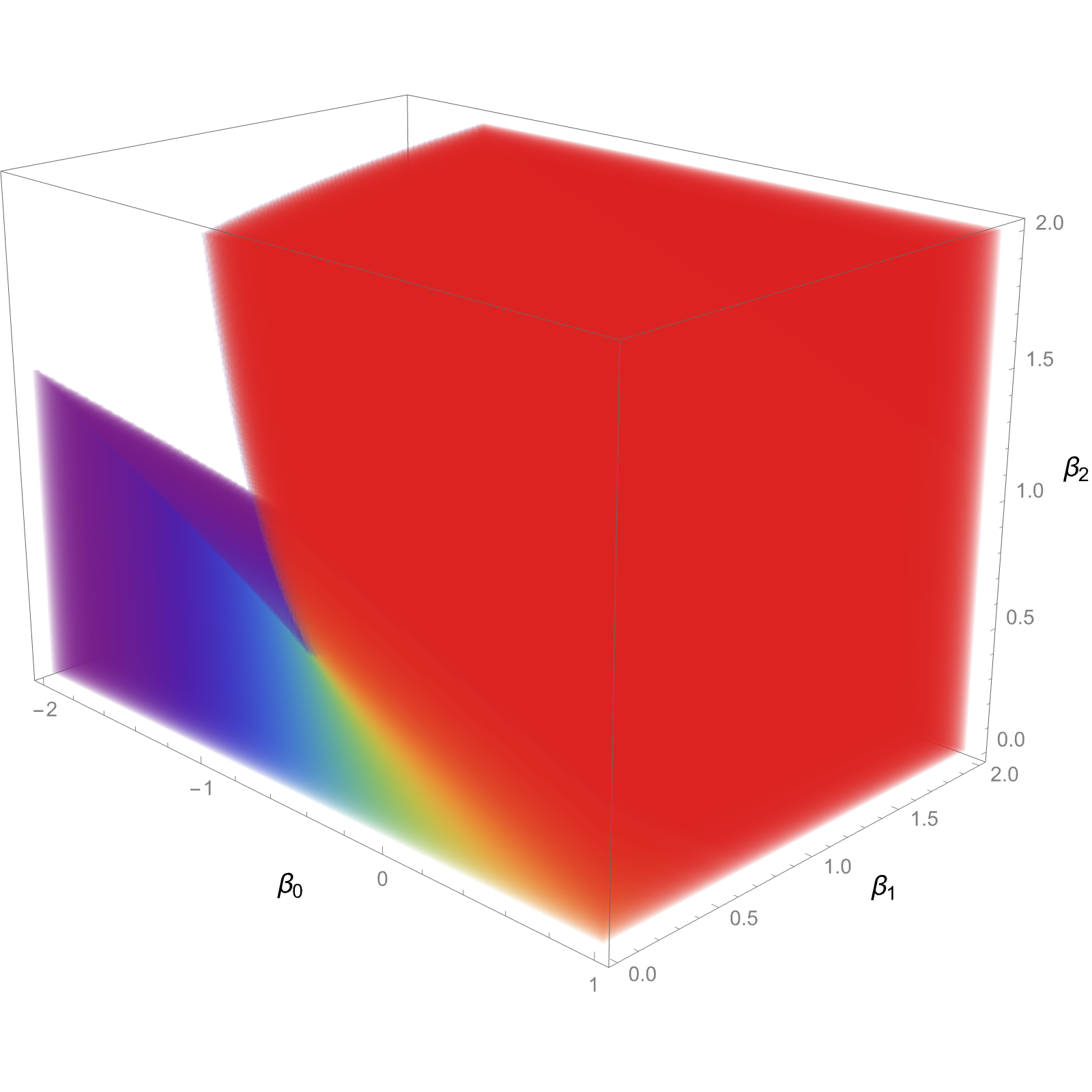}}$
\hspace{0.1cm}
$\underset{\text{Supercritical}}{\includegraphics[scale=0.18,align=c]{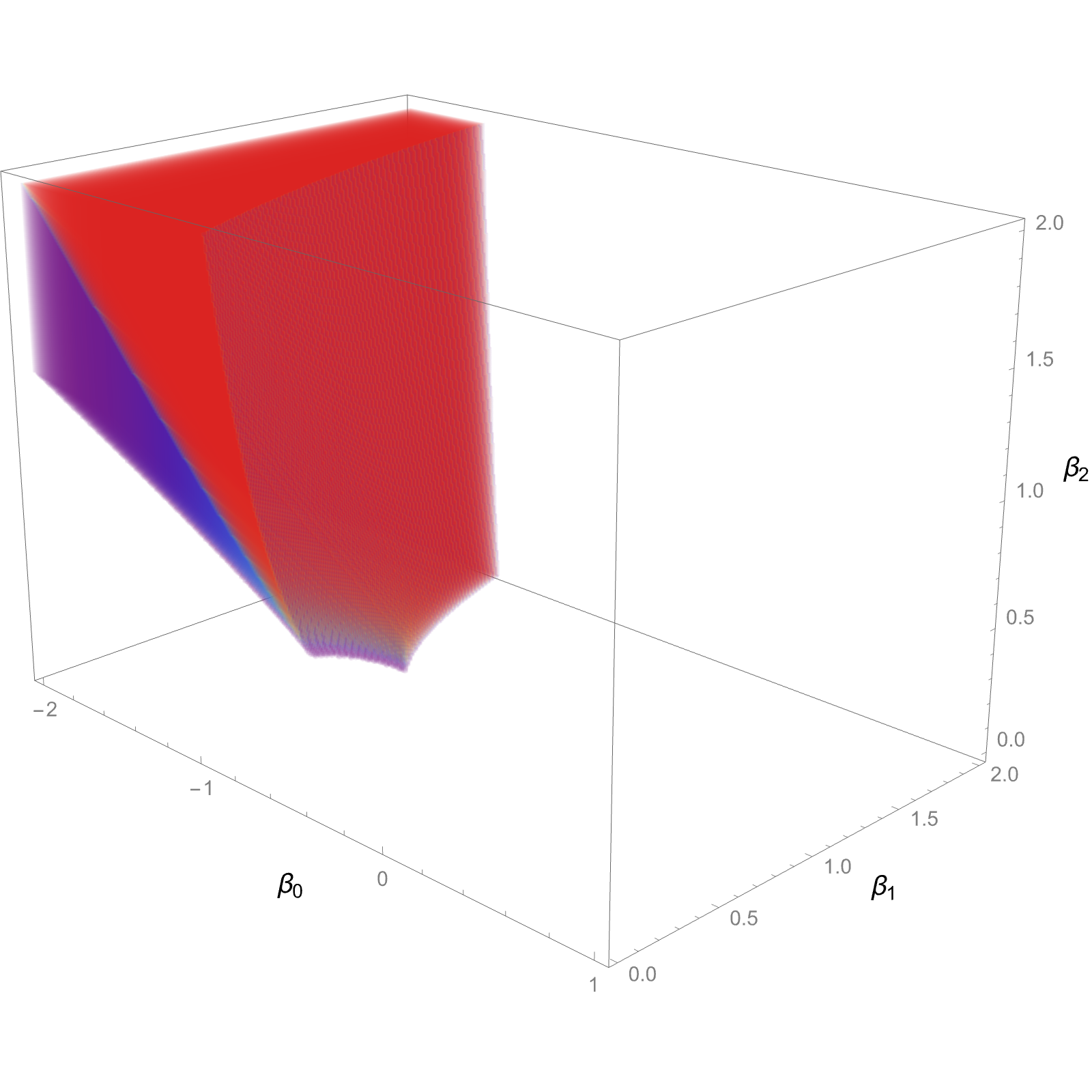}}$
\caption{Plots of the three relevant parameter regimes, for the same ERGM as in Figure \ref{fig:lbeta}
(with $G_0,G_1,G_2$ being an edge, wedge, and triangle) as $\beta = (\beta_0,\beta_1,\beta_2)$ varies
in $[-2,1] \times [0,2] \times [0,2]$.
The colors indicate the value of $p$, the optimizer of $L_\beta$, with purple denoting $p=0$ and red denoting $p=1$,
other colors interpolating between these two extremes.
Left: Dobrushin or ``very high temperature'' regime, where the natural dynamics exhibit uniform contraction.
Middle: subcritical or ``high temperature'' regime, where the dynamics exhibit rapid mixing.
Right: supercritical or ``low temperature'' regime, where the dynamics mix slowly from a worst-case initialization
but exhibit \emph{metastable mixing} from a warm start within a metastable well.
The supercritical regime also includes the phase transition surface, on which there is
\emph{phase coexistence}.}
\label{fig:regimes}
\end{figure}

Unlike in $Y$, the sample from the Erd\H{o}s--R\'enyi model $\cG(n,p)$, the edges of $X$ are \emph{not} independent.
Nevertheless, various approximate independence results do indeed hold \cite{bhamidi2008mixing,reinert2019approximating,ganguly2024sub,winstein2025concentration}.
Of particular relevance to the present work is a \emph{Gaussian concentration inequality} for Lipschitz
observables of $X$.
This was first shown in \cite{ganguly2024sub} in the \emph{subcritical} regime of parameters, which will
be described properly in Section \ref{sec:review_fluctuations_dynamical} below, but briefly it is the regime
of $\beta$ for which $L_\beta$ has a unique \emph{local} maximum at which the function is strictly concave,
and it should be thought of as a ``high temperature'' regime.
In particular, throughout the entire subcritical regime, we have phase uniqueness.
See Figure \ref{fig:regimes} (middle) for a sample plot of the subcritical regime.

A key input to this argument was the rapid mixing of \emph{Glauber dynamics} (a natural Markov chain for which
the ERGM measure is stationary) in the subcritical regime, as shown by \cite{bhamidi2008mixing}.
The Gaussian concentration inequality was later extended to the \emph{supercritical} or ``low-temperature''
regime by \cite{winstein2025concentration}, using the recent \emph{metastable mixing} result of
\cite{bresler2024metastable}, which proves rapid mixing from a \emph{warm start} within a
\emph{large metastable well} of the Glauber dynamics;
for more on this, see Section \ref{sec:review_fluctuations_dynamical}.
The supercritical regime is where $L_\beta$ has multiple local maxima, and where $L_\beta$ is strictly
concave at its \emph{global} maximum.
Note that throughout most of the supercritical regime we still have phase uniqueness,
i.e.\ a unique \emph{global} maximum for $L_\beta$.
However, the supercritical regime also contains the surface of phase transition, on which we have
phase coexistence.
See Figure \ref{fig:regimes} (right) for a sample plot of the supercritical regime,
where one edge of the phase transition surface may also be seen;
this surface was characterized in some cases by \cite{yin2013critical}.
Finally, we remark that in the phase coexistence case, the mixture measures and the large
metastable wells just mentioned are essentially the same thing.

\subsubsection{Central limit theorems}
\label{sec:intro_setup_clt}

Despite the aforementioned advances in understanding the independence properties of the edges of $X$ in a
quantitative way, proving a central limit theorem for the edge count in ERGMs remained an elusive open
problem for some time.
Progress towards this goal was made by \cite{mukherjee2023statistics}, which covers a special case,
the \emph{$2$-star} ERGM where, in our notation, $K = 1$ and $G_1$ is a $2$-star or wedge graph.
This ERGM may also be viewed as an \emph{Ising model} as it only has pairwise interactions, and this structure
allows for specialized analysis leading to central limit theorems for the edge count and vertex degree,
as well as other results.
Another line of reasoning in \cite{bianchi2024limit} used analyticity properties of the \emph{free energy}
to study the \emph{edge-triangle} model where, as mentioned above, $K = 1$ and $G_1$ is a triangle.
Among other results, they achieved a CLT for the edge count in this model, and they also claimed
that their methods extend to more general ERGMs.
However, it should be noted that all of these previous results are \emph{non-quantitative},
simply providing convergence in distribution and giving
no bound on any distance between the distribution of the edge count and a Gaussian.

A recent breakthrough \cite{fang2024normal} developed a novel application of \emph{Stein's method},
which yields quantitative central limit theorems for certain observables of general \emph{nonlinear exponential
families} including ERGMs as well as other spin systems (which are sometimes also called graphical models).
In that work, the authors applied their general result to the edge count in ferromagnetic ERGMs in the
subcritical regime of parameters, and derived a central limit theorem which is quantitative in terms of
both the \emph{Kolmogorov distance} and the \emph{Wasserstein distance}.
They also obtain a similar result for subgraph counts.

A key input to their argument is the fact that the fluctuations of $\sN_G(X)$ are driven by those of the
total edge count $\sE(X)$.
This was initially proved by \cite{sambale2020logarithmic} in the case where $G$ is a triangle,
by using a \emph{multilevel concentration inequality} derived from functional-analytic tools such as a
\emph{modified log-Sobolev inequality} in the \emph{Dobrushin regime}, which should be thought of
as a ``very high temperature'' regime and is a subset of the subcritical regime.
As such, the work of \cite{fang2024normal} was initially restricted to the Dobrushin regime.
However, during the preparation of the present article, \cite{fang2024normal} was updated to include a
new argument which covers the full subcritical regime, albeit with worse quantitative bounds than
in the Dobrushin regime.
See Figure \ref{fig:regimes} (left) for a sample plot of the Dobrushin regime.

In the present article, we introduce a new probabilistic technique which allows us to extend these results to the
supercritical regime of parameters, with a caveat in the phase coexistence case which will be clarified
in Section \ref{sec:intro_results} below.
Our technique also improves the quantitative Kolmogorov and Wasserstein distance bounds
obtained by \cite{fang2024normal} in the subcritical regime (beyond the Dobrushin regime),
and moreover we believe it provides more
intuition for the mechanisms at hand, as compared to the more algebraically technical proofs of
\cite{sambale2020logarithmic,fang2024normal}.
Finally, we also use our methods to obtain a quantitative CLT for the \emph{degree} of a particular
vertex in $X$, as well as \emph{local subgraph counts}, i.e.\ the number of subgraphs adjacent to a particular
vertex.
To the best of our knowledge, these local results are novel in all parameter regimes.
\subsection{Statement of results}
\label{sec:intro_results}

Our results are quantitative in both the Wasserstein distance and the Kolmogorov distance between
$\R$-valued random variables $S$ and $Z$, so we first recall these standard definitions:
\begin{align}
    \label{eq:wasdist}
    \dW(S,Z) &\coloneqq \sup_{\substack{\psi : \R \to \R \\ 1\text{-Lipschitz}}} \;
        \left| \E[\psi(S)] - \E[\psi(Z)] \right|, \\
    \label{eq:koldist}
    \dK(S,Z) &\coloneqq \sup_{s \in \R} \; \left| \P[S \leq s] - \P[Z \leq s] \right|.
\end{align}
We remark that the Wasserstein distance may also be expressed as
$\E[|S'-Z'|]$, where $S'$ and $Z'$ have the same marginal 
distributions as $S$ and $Z$ but are coupled in a way which minimizes this expectation.

Let us fix a ferromagnetic ERGM, meaning that we have fixed graphs $G_0, G_1, \dotsc, G_K$ with $G_0$ a single
edge, and fixed parameters $\beta = (\beta_0, \beta_1, \dotsc, \beta_K) \in \R \times \R_{\geq0}^K$.
We assume that $L_\beta$, as defined in \eqref{eq:Ldef}, has at least one \emph{strictly concave global
maximum}, meaning a density $p \in M_\beta$ (the set of global maxima of $L_\beta$) which also satisfies
$L_\beta''(p) < 0$.
Let us denote by $U_\beta \sse M_\beta$ the set of such strictly concave global maxima of $L_\beta$;
the strict concavity has implications for the mixing and concentration properties of the ERGM, which
will be explained further in Section \ref{sec:review_fluctuations_dynamical}.

To make the statements easier to digest, we first consider the \emph{non-critical phase uniqueness}
case, where $U_\beta = M_\beta = \{p\}$; note that this is how we define the value of $p$ used
in the statements.
Additionally, note that this includes most of the supercritical regime, only avoiding the surface of phase transition
(see Figure \ref{fig:regimes} (right)) which has Lebesgue measure zero.
In this case, we simply let $X$ denote a sample from the full ERGM measure \eqref{eq:gibbsmeasure}
with Hamiltonian \eqref{eq:hdef}.
Moreover, here and throughout the work, for the graphs $G_j$ in the definition \eqref{eq:hdef},
we write $G_j = (\cV_j, \cE_j)$ and set $\sv_j = |\cV_j|$ and $\se_j = |\cE_j|$.
Additionally, we let $Z \sim \cN(0,1)$.

\subsubsection{Global results}
\label{sec:intro_results_global}

Our first main result is the following quantitative central limit theorem for the \emph{edge count}
$\sE(X)$.

\begin{theorem}
\label{thm:global_clt_intro}
Define
\begin{equation}
\label{eq:sigmadef_global_intro}
    \sigma_n^2 = 
    \left( 1 - 2 p (1-p) \sum_{j=1}^K \beta_j \se_j (\se_j-1) p^{\se_j-2} \right)^{-1}
    \times p(1-p) \binom{n}{2}.
\end{equation}
Then for any $\eps > 0$, we have
\begin{equation}
    \d \left( \fr{\sE(X) - \E[\sE(X)]}{\sigma_n}, Z \right) \lesssim n^{-\fr{1}{2} + \eps},
\end{equation}
for both $\d = \dW$ and $\d = \dK$.
\end{theorem}

Note that in the definition of $\sigma_n^2$ in \eqref{eq:sigmadef_global_intro},
the prefactor is strictly greater than $1$, meaning that the fluctuations of
the edge count are strictly larger under the ERGM than they are under $\cG(n,p)$, 
where the variance is exactly $p (1-p) \binom{n}{2}$.
This should be expected due to the ferromagnetic condition which is a sort of positive-correlation
condition between the edges.

Theorem \ref{thm:global_clt_intro} improves upon \cite[Theorem 3.1]{fang2024normal} in multiple ways.
First of all, it covers the supercritical regime whereas \cite[Theorem 3.1]{fang2024normal} was restricted
to the subcritical regime.
Secondly, the Wasserstein and Kolmogorov distance bounds obtained by \cite[Theorem 3.1]{fang2024normal} are
of order $n^{-\fr{1}{2}}$ in the \emph{Dobrushin} regime, but only order $n^{-\fr{1}{4}}$ in the remainder
of the subcritical regime, whereas our Theorem \ref{thm:global_clt_intro} yields a bound of
order $n^{-\fr{1}{2}+\eps}$ throughout all parameter regimes.
Finally, we remark that we feel our proof strategy provides a more hands-on approach than that of
\cite{fang2024normal} and gives more intuition about the mechanism involved in the key input,
namely the fact that the fluctuations of $\sN_G(X)$ are controlled by those of $\sE(X)$,
as will be discussed in Section \ref{sec:intro_iop_hajek} and presented explicitly as
Proposition \ref{prop:global_hajek_intro}.

We remark that in the subcritical regime, as brought to our attention by Xiao Fang,
we may replace $\E[\sE(X)]$ by an explicit closed-form expression in this CLT.
Indeed, in \cite{fang2025conditionalcentrallimittheorems} a constant $c_*$
(with an explicit formula depending on the ERGM specification) was found which satisfies
\begin{equation}
    \left| \E[\sE(X)] - \binom{n}{2} \left( p - \fr{c_*}{n} \right) \right| \lesssim \sqrt{n}
\end{equation}
in the subcritical regime.
This difference is small enough to be washed out upon rescaling by $\sigma_n$, which is of order $n$.
The formula for $c_*$ is rather complicated, however, so we refer the reader to
\cite[Lemma B.1]{fang2025conditionalcentrallimittheorems} for the exact expression.

As a corollary of Theorem \ref{thm:global_clt_intro}, also using Proposition \ref{prop:global_hajek_intro} below,
we obtain the following central limit theorem for the subgraph count $\sN_G(X)$.
For the following result, we assume $G = (\cV, \cE)$ is a fixed graph with $\sv = |\cV|$
and $\se = |\cE|$, and as before we set $Z \sim \cN(0,1)$.

\begin{corollary}
\label{cor:global_subgraph_clt_intro}
Recall the definition of $\sigma_n^2$ from \eqref{eq:sigmadef_global_intro}.
For any $\eps > 0$, 
\begin{equation}
    \dW \left( \fr{\sN_G(X) - \E[\sN_G(X)]}{2 \se p^{\se-1} n^{\sv-2} \sigma_n}, Z \right) \lesssim n^{-\fr{1}{2}+\eps}.
\end{equation}
\end{corollary}

In this corollary, we obtain the \emph{same} error bound in terms of
the Wasserstein distance as in Theorem \ref{thm:global_clt_intro}.
Due to our proof strategy, we do not provide the same error bound for the Kolmogorov distance,
but one may apply the standard inequality
\begin{equation}
\label{eq:kolwas}
    \dK(f(X), Z) \lesssim \sqrt{\dW(f(X),Z)},
\end{equation}
which holds since $Z \sim \cN(0,1)$,
to obtain quantitative CLTs for $\sN_G(X)$ in terms of the Kolmogorov distance, albeit with
a worse error bound than what appears for $\sE(X)$ in Theorem \ref{thm:global_clt_intro}.

Note that even though it does improve upon the result of \cite{fang2024normal},
it is not clear if the Kolmogorov distance bound appearing in Theorem \ref{thm:global_clt_intro} is tight.
Indeed, simulations presented in \cite{winstein2025concentration} indicate that the true
Kolmogorov distance between the shifted and scaled $\sE(X)$ and a Gaussian is of order $n^{-1}$,
which is what one would obtain for truly independent edges.
Note however that such a result may require a lower-order correction to the variance proxy $\sigma_n^2$.

Finally, note that a bound of $n^{-\fr{1}{2}} \polylog n$ on both distances should be achievable using
the same methods presented in this work, but we do not focus on optimizing the bound to this degree
since, as mentioned, such a result would still likely be far from optimal.

\subsubsection{Local results}
\label{sec:intro_results_local}

The strategy we develop for Theorem \ref{thm:global_clt_intro} also applies to the case of the \emph{degree}
of a particular (deterministic) vertex $v \in [n]$.
We denote by $\deg_v(X)$ this degree of $v$ in $X$, which is the number of edges in $X$ adjacent to $v$.
As previously mentioned, to the best of our knowledge, the following result, or even a non-quantitative version,
has not appeared in \emph{any} parameter regime, aside from in the special case of the
$2$-star ERGM \cite{mukherjee2023statistics}.
To state the result, in addition to the notation introduced above Theorem \ref{thm:global_clt_intro},
for any vertex $\rho \in \cV_j$, we set $\sd_\rho$ to be the degree of $\rho$ in the graph $G_j$.

\begin{theorem}
\label{thm:local_clt_intro}
Define
\begin{equation}
\label{eq:sigmadef_local_intro}
    \varsigma_n^2 = \left( 1 - p(1-p) \sum_{j=1}^K \beta_j p^{\se_j-2} \sum_{\rho \in \cV_j} \sd_\rho (\sd_\rho-1) \right)^{-1} \times p (1-p) (n-1).
\end{equation}
Then for any deterministic vertex $v \in [n]$ and any $\eps > 0$, we have
\begin{equation}
    \d \left( \fr{\deg_v(X) - \E[\deg_v(X)]}{\varsigma_n}, Z \right) \lesssim n^{-\fr{1}{4}+\eps},
\end{equation}
for both $\d = \dW$ and $\d = \dK$.
\end{theorem}

Again, $\varsigma_n^2$ is strictly greater than the variance of the degree of a vertex in $\cG(n,p)$,
meaning the fluctuations of the degree of a vertex in an ERGM are larger than those in the corresponding
Erd\H{o}s--R\'enyi model.
As with the CLT for $\sE(X)$, we may replace $\E[\deg_v(X)]$ in the above theorem by an explicit
(but rather complicated) expression in the subcritical regime using
\cite[Lemma B.1]{fang2025conditionalcentrallimittheorems}.

Like in the global setting, we also obtain the following CLT for \emph{local} subgraph counts as a corollary of
Theorem \ref{thm:local_clt_intro} by using Proposition \ref{prop:local_hajek_intro} below,
which is the local analog of Proposition \ref{prop:global_hajek_intro}.
In the following statement, the quantity $\sN_G^v(x)$ is defined to be the number of homomorphisms of $G$
in $x$ which have $v \in [n]$ in their image.

\begin{corollary}
\label{cor:local_subgraph_clt_intro}
Recall the definition of $\varsigma_n^2$ from \eqref{eq:sigmadef_local_intro}.
For any $\eps > 0$,
\begin{equation}
    \dW \left( \fr{\sN_G^v(X) - \E[\sN_G^v(X)]}{2 \se p^{\se-1} n^{\sv-2} \varsigma_n}, Z \right) \lesssim n^{-\fr{1}{4} + \eps}.
\end{equation}
\end{corollary}

Recalling how homomorphism counts lead to multipliers based on automorphisms of $G$, the homomorphism count
$\sN_\triangle^v(x)$ (for instance) is $6$ times what is usually thought of as the number of triangles adjacent to $v$.
Be warned, though, that when the graph $G$ is \emph{not} vertex-transitive (i.e.\ there are two vertices which
cannot be mapped to one another by an automorphism), $\sN_G^v(x)$ may not be simply a multiple of
the typical count of copies of $G$ adjacent to $v$ in $x$.
For such cases, a better understanding is provided by considering the number of homomorphisms
which map a \emph{particular vertex} $\rho \in \cV$ to $v$.
A CLT similar to Corollary \ref{cor:local_subgraph_clt_intro} but for these restricted counts
is also provided, as Corollary \ref{cor:local_subgraph_clt} below.

As with Corollary \ref{cor:global_subgraph_clt_intro}, the Wasserstein distance bound in Corollary \ref{cor:local_subgraph_clt_intro}
is the same as for $\deg_v(X)$, and a quantitative Kolmogorov
distance bound may be obtained by applying the inequality \eqref{eq:kolwas}.
As was the case for $\sE(X)$, it is not clear if the Kolmogorov distance bound obtained in
Theorem \ref{thm:local_clt_intro} is optimal, as the simulations of \cite{winstein2025concentration}
indicate that the true Kolmogorov distance between a Gaussian and the distribution of the shifted and
scaled $\deg_v(X)$ is of order $n^{-\fr{1}{2}}$, again possibly requiring
a lower-order correction to $\varsigma_n^2$.

\subsubsection{Phase coexistence and criticality}
\label{sec:intro_results_pcc}

Let us now address what happens in the \emph{phase coexistence} case where $|M_\beta| > 1$.
In this case, as mentioned in Section \ref{sec:intro_setup_macroscopic}, the ERGM decomposes as a mixture of measures,
one for each $p \in M_\beta$, and we must work with one measure in the mixture at a time.
Note that our results only apply for mixture measures corresponding to $p \in U_\beta \sse M_\beta$,
namely we require that $L_\beta''(p) < 0$.
See the last paragraph in this section for a discussion of the alternative case, which we call \emph{critical}.

In practice, working with one of the mixture measures amounts to taking $X$ to be a sample from the ERGM,
conditioned so that the \emph{cut distance} between $X$ and a constant graphon with density $p \in U_\beta$
is at most some small number $\eta > 0$.
Graphons and the cut distance will be properly introduced below, in Section \ref{sec:review_leadingorder_graphons}.
Note that the set of graphs which we condition on here is often referred to as a \emph{metastable well around $p$}
due to its behavior under Glauber dynamics, which will be explained in Section \ref{sec:review_fluctuations_dynamical}.

Due to the lack of the \emph{FKG inequality} on these conditioned measures, we must restrict our scope somewhat,
and only consider ERGMs where the defining graphs $G_0, G_1, \dotsc, G_K$ are all \emph{forests}.
This assumption arises as a requirement for a natural method for getting around the
use of the FKG inequality, as first presented in \cite[Proposition 6.2]{winstein2025concentration},
which is restated as Proposition \ref{prop:Eproduct} in the present work.
Under this assumption, if $X$ denotes a sample from the ERGM \emph{conditioned on some metastable well around $p \in U_\beta$}
then Theorems \ref{thm:global_clt_intro} and \ref{thm:local_clt_intro} hold as stated.
We note that follow-up work by Mukherjee and the author \cite{mukherjee2026approximate} removes
this assumption by proving an approximate version of the FKG inequality.
However, for completeness we retain this assumption in the present work, and state it explicitly as
follows for future reference.

\begin{assumption}
\label{assumption}
We always assume the data specifying the Hamiltonian \eqref{eq:hdef}, i.e.\ the graphs $G_0, G_1, \dotsc, G_K$
and parameters $\beta = (\beta_0, \beta_1, \dotsc, \beta_K) \in \R \times \R_{\geq 0}$,
satisfy at least one of the following two statements
\begin{enumerate}
    \item Either $|U_\beta| = |M_\beta| = 1$,
    \item or all graphs $G_0, G_1, \dotsc, G_K$ are forests.
\end{enumerate}
Moreover, we always assume that $p \in U_\beta$.
\end{assumption}

As for Corollaries  \ref{cor:global_subgraph_clt_intro} and \ref{cor:local_subgraph_clt_intro}, to hold
exactly as stated in the phase coexistence case, we also require that the graph $G$ itself be a forest.
However, this is only so that we obtain the explicit form of the coefficient $2 \se p^{\se-1} n^{\sv-2}$
in front of $\sigma_n$ or $\varsigma_n$ in the denominator.
In the case where $G$ is not a forest (but $G_0, G_1, \dotsc, G_K$ are), we do still obtain CLTs for
$\sN_G(X)$ and $\sN_G^v(X)$, but the coefficient $2 \se p^{\se-1} n^{\sv-2}$ must be replaced by something
with no simple closed form.
This will be more fully explained at the end of Section \ref{sec:hajek_subgraphs}.

Finally, we address the case where $p \in M_\beta \setminus U_\beta$, which is to say that $p$ is a global
maximum of $L_\beta$, but $L_\beta$ is \emph{not strictly concave} at $p$.
In this \emph{critical} case, we do not expect Gaussian fluctuations for the edge count.
Indeed, this case was considered in \cite{mukherjee2023statistics} for the $2$-star ERGM, and there it was
found that the fluctuations of the edge count approach a \emph{different} limiting distribution, with density
proportional to
\begin{equation}
    \Exp{- \fr{\zeta^2}{2} - \fr{\zeta^4}{24}}
\end{equation}
for $\zeta \in \R$.
Additionally, it was conjectured in \cite{bianchi2024limit} that in the edge-triangle model, the fluctuations
of the edge count at the critical point should have limiting density proportional to
$\Exp{-\fr{81}{64} \zeta^4}$.

Interestingly, in the case of the degree $\deg_v(X)$, one still obtains a Gaussian limit at criticality
in the $2$-star ERGM \cite[Theorem 1.6]{mukherjee2023statistics}.
This is the only ERGM for which rigorous results about the critical case have been obtained,
but one might expect similar behavior for more general ERGMs.
However, we do not consider the critical case in the present work as the requisite mixing inputs are not
yet available.
\subsection{Proof ideas}
\label{sec:intro_iop}

At the broadest level, the strategy of the proof follows as in \cite{fang2024normal}, using their general
result which is an application of Stein's method \cite{stein1972bound} to nonlinear exponential families.
The statement of this general result is too long to present here, but the interested reader should
check Theorem \ref{thm:stein} in Section \ref{sec:stein_statement} below,
which is a restatement of \cite[Theorem 2.1]{fang2024normal}.
Additional context for Stein's method is also provided after the statement of Theorem \ref{thm:stein}.
As we do not provide any new modification or use of Stein's method itself, we do not presently provide
an overview of the method, instead referring to the excellent survey \cite{ross2011fundamentals}, as well as
\cite{chatterjee2008new}, which the work of \cite{fang2024normal} is primarily based on.

For any mean-zero function $f$ of a random variable $X$ distributed according to an exponential tilt of a product
measure, as is the case with a ferromagnetic ERGM,
Theorem \ref{thm:stein} gives explicit expressions for error terms bounding $\d(f(X),Z)$,
where $Z \sim \cN(0,1)$ and $\d$ is either $\dW$ or $\dK$.
These expressions can be found in Section \ref{sec:stein_statement}, just before the statement of
Theorem \ref{thm:stein}, and they involve the function $f$ as well as the exponential tilting function.
Using Theorem \ref{thm:stein}, deriving a CLT such as Theorem \ref{thm:global_clt_intro} reduces
to simply bounding all of these error terms appropriately.

The main technical contribution of the present work is the introduction of a new technique, based on concentration
inequalities derived from mixing, for bounding the error terms given by Theorem \ref{thm:stein}.
We remark that this idea is quite flexible, and we expect a similar technique to work for a variety of spin systems
other than ERGMs, and for a variety of observables.
Some domain-specific work is needed in simplifying the error terms, but since they can mostly be written as expectations
of products of factors which have means close to zero, we expect concentration inequalities to be useful regardless.

\subsubsection{H\'ajek projections}
\label{sec:intro_iop_hajek}

When deriving a CLT for $\sE(X)$ under a ferromagnetic ERGM, most of the error terms given in
Theorem \ref{thm:stein} can be bounded routinely.
However, one of the error terms involves the variance of
\begin{equation}
\label{eq:hajekdifference}
    \sN_G(X) - 2 \se p^{\se-1} n^{\sv-2} \sE(X),
\end{equation}
for certain graphs $G = (\cV,\cE)$, where we set $\sv = |\cV|$ and $\se = |\cE|$.
It turns out that in order to conclude a CLT for the edge count as in Theorem \ref{thm:global_clt_intro},
one must show that the variance of \eqref{eq:hajekdifference} is much smaller than $n^{2\sv-2}$.
As the order of the variance of the second term  $2 \se p^{\se-1} n^{\sv-2} \sE(X)$ is itself $n^{2\sv-2}$
(although this fact is only rigorously clear \emph{after} proving the CLT for $\sE(X)$), this bound on
the variance of \eqref{eq:hajekdifference} also allows one to conclude a CLT for $\sN_G(X)$ (such as
Corollary \ref{cor:global_subgraph_clt_intro}) from a CLT for $\sE(X)$ (such as Theorem \ref{thm:global_clt_intro}).

It's worth commenting on the form of the prefactor $2 \se p^{\se-1} n^{\sv-2}$ appearing in
\eqref{eq:hajekdifference}.
Heuristically, this is what one would expect for $\sN_G(x,e)$, which denotes the number of homomorphisms of $G$ in $x$
\emph{which use the edge $e$}, assuming that $e$ is present in $x$, and that $x$ is relatively homogeneous with
density $p$.
Indeed, there are $2 \se$ edges $\phi \in \cE$ which could be mapped to $e$ (counting different orientations
separately).
For each choice of $\phi$, it remains to map the $\sv-2$ vertices of $G$ into $[n]$ somewhere,
meaning there are $2 \se n^{\sv-2}$ possible maps to consider.
Then, for such a map to be a homomorphism, the images of all remaining $\se-1$ edges of $G$ must be present in
$x$, and this happens with probability approximately $p^{\se-1}$ by our homogeneity and density assumptions.

This intuition leads one to expect that $2 \se p^{\se-1} n^{\sv-2} \sE(X)$ is a good proxy for $\sN_G(X)$.
Indeed, in the Erd\H{o}s--R\'enyi case $Y \sim \cG(n,p)$, the quantity $2 \se p^{\se-1} n^{\sv-2} \sE(Y)$
is the \emph{H\'ajek projection} of $\sN_G(Y)$ onto the variable $\sE(Y)$.
This is the $\sE(Y)$-measurable variable which best approximates $\sN_G(Y)$ in $L^2$.
We remark that this idea has previously been called a \emph{Hoeffding decomposition} by other authors
\cite{sambale2020logarithmic,fang2024normal}; however, the two terminologies seem to refer to the same thing,
and in the present work we use the term H\'ajek projection.

While $2 \se p^{\se-1} n^{\sv-2} \sE(X)$ may not be exactly the same as the H\'ajek projection of
$\sN_G(X)$ onto $\sE(X)$ for ERGM-distributed $X$, with non-independent edges,
we nonetheless do still obtain quite a good approximation.
The key reason for this is that the count $\sN_G(X,e)$ of homomorphisms using the edge $e$ concentrates
tightly around $2 \se p^{\se-1} n^{\sv-2}$, simultaneously for all $e \in \edgeset$, which is our notation for the
set of all possible edges of $X$.
This means that each edge accounts for close to the correct number of subgraphs.

\subsubsection{Concentration via mixing}
\label{sec:intro_iop_concentration}

To turn the idea of H\'ajek projections into a rigorous and useful bound, we apply various concentration
inequalities derived from mixing.
The intuition for such inequalities is that if a Markov chain $(X_t)$ mixes rapidly to its stationary
distribution $X$, and if at each step of the chain $f(X_t)$ does not change by much, then the random
variable $f(X)$ must be concentrated.

The rapid mixing of the ERGM Glauber dynamics chain is provided by \cite{bhamidi2008mixing} in the subcritical
case, and more recently \cite{bresler2024metastable} showed \emph{metastable mixing} for the Glauber dynamics
within a metastable well from a \emph{warm start}.
In the context of Glauber dynamics where edges are updated one by one,
``$f(X_t)$ not changing by much'' may be witnessed by a bound on the Lipschitz
vector of $f$, which will be defined in Section \ref{sec:review_fluctuations_concentration} below.
The first instance of a concentration inequality for Lipschitz observables of the ERGM derived from mixing
is due to \cite{ganguly2024sub}, which used the rapid mixing result of \cite{bhamidi2008mixing}
as input to a technique developed by \cite{chatterjee2005concentration}.
Later, \cite{winstein2025concentration} modified this technique and used it with the metastable mixing
result of \cite{bresler2024metastable} to extend the result of \cite{ganguly2024sub} to metastable
wells in the supercritical regime.

In addition to these results for Lipschitz observables, we rely on a general concentration inequality due to
\cite{barbour2022long} which will be introduced in Section \ref{sec:inputs_concentration}.
Like the modification of the technique of \cite{chatterjee2005concentration} provided in \cite{winstein2025concentration},
this result of \cite{barbour2022long} allows us to restrict to certain \emph{large subsets} of the space of all graphs,
wherein nice behavior is observed.
We thus embark on a two-step strategy: we first show that the quantities $\sN_G(X,e)$, which are exactly
the \emph{changes} of $\sN_G(X)$ when an edge $e$ is flipped under Glauber dynamics, themselves concentrate,
as mentioned above.
Then we restrict to a subset of the state space where this concentration holds, and apply a second concentration
inequality which \emph{uses the concentration of the changes} in the bound.
Such a strategy is necessary here because simply calculating the Lipschitz vector of the quantity
\eqref{eq:hajekdifference}
does not lead to a sufficient bound, and a more careful analysis of the changes along the trajectory is needed.

\subsubsection{Key propositions}
\label{sec:intro_iop_props}

The strategy outlined above leads to the following result, which shows that the fluctuations
of $\sN_G(X)$ match closely with those of $2 \se p^{\se-1} n^{\sv-2} \sE(X)$.
In the following result, either the ERGM exhibits phase uniqueness, in which case we allow $G = (\cV,\cE)$
to be an arbitrary graph and set $\sv = |\cV|$ and $\se = |\cE|$, or, in the phase coexistence case, we
set $X$ to be an ERGM sample conditioned on a metastable well around $p \in U_\beta$, and we also assume
that $G$ and the specifying graphs $G_0, G_1, \dotsc, G_K$ are forests.

\begin{proposition}
\label{prop:global_hajek_intro}
Under Assumption \ref{assumption}, define
\begin{equation}
\label{eq:hatngdef}
    \hsN_G(x) \coloneqq \sN_G(x) - 2 \se p^{\se-1} n^{\sv-2} \sE(x).
\end{equation}
Then, for all $\zeta \in (0,1)$, there is some $c > 0$ such that 
\begin{equation}
    \P \left[
        \left|
            \hsN_G(X) - \E\left[\hsN_G(X)\right]
        \right| > n^{\sv-1.5+\zeta}
    \right] \leq \Exp{-cn^\zeta}.
\end{equation}
In particular, for all $\xi > 0$, we have
\begin{equation}
    \Var\left[\hsN_G(X)\right] \lesssim n^{2\sv-3+\xi}.
\end{equation}
\end{proposition}

One should compare this result with \cite[Theorem 3.2]{sambale2020logarithmic} which was restricted
to the Dobrushin (very high-temperature) regime, whereas Proposition \ref{prop:global_hajek_intro}
holds for \emph{all} parameter regimes.
Additionally, \cite[Theorem 3.2]{sambale2020logarithmic} is specific to triangle counts, but it
can easily be generalized to all subgraph counts, as was done in \cite[Lemma 5.1]{fang2024normal}.
Let us just consider the triangle case $G = \triangle$,
where \cite[Theorem 3.2]{sambale2020logarithmic} states
\begin{equation}
    \P\left[ \left| \hsN_\triangle(X) - \E\left[\hsN_\triangle(X)\right] \right| > \lambda \right]
    \leq 2 \Exp{- c \cdot \min \left\{\left(\fr{\lambda}{n^{\fr{3}{2}}}\right)^{\fr{2}{3}},
    \fr{\lambda}{n^{\fr{3}{2}}} \right\}},
\end{equation}
for some constant $c > 0$.
Note that the first expression in the minimum will always be the one chosen when
$\lambda \gg n^{\fr{3}{2}}$, which is the only case where the above inequality is useful.
So, plugging in $\lambda = n^{\fr{3}{2}+\zeta}$ yields an upper bound of $\Exp{-c n^{\fr{2}{3} \zeta}}$.
On the other hand, Proposition \ref{prop:global_hajek_intro} achieves an upper bound of
$\Exp{-c n^\zeta}$ in this case, which is much smaller.
Moreover, the validity of Proposition \ref{prop:global_hajek_intro} for all $\zeta < 1$ is essentially
tight in the supercritical regime, since there may be \emph{small metastable regions} with atypical behavior,
of probability $e^{-O(n)}$ \emph{within} large metastable wells; see
\cite[Theorem 3.3]{bresler2024metastable} for more on this phenomenon.

We remark that, as was brought to our attention by Xiao Fang, the tools used to prove
\cite[Theorem 3.2]{sambale2020logarithmic} (namely the so-called \emph{higher-order concentration inequalities})
have actually been extended to the full subcritical regime in \cite[Appendix C]{fang2025conditionalcentrallimittheorems}.
However, beyond the Dobrushin regime they achieve a worse rate than \cite{sambale2020logarithmic},
meaning that Proposition \ref{prop:global_hajek_intro} still provides a better bound.
As for the supercritical regime, follow-up work by the author \cite{winstein2026wasserstein} provides
another application, using machinery which extends that in the present work,
which fits into the framework of higher-order concentration inequalities.
In that situation, the quantity $\partial_e \sN_G$ is of interest, and in this context
edge correction is not enough, as may be seen even in the Erd\H{o}s--R\'enyi model by a simple calculation.
In \cite{winstein2026wasserstein}, higher-order (namely quadratic) correction terms are necessary.
Furthermore, follow-up work in progress will address the general higher-order concentration phenomenon
in the supercritical regime via another method, yielding effective bounds but without optimal tail estimates.

We should also remark that even in the phase coexistence case we may also treat non-forest graphs $G$
and $G_0, G_1, \dotsc, G_K$ in Proposition \ref{prop:global_hajek_intro}, but we obtain
a slightly different result.
Namely, the coefficient $2 \se p^{\se-1} n^{\sv-2}$ must be replaced by $\E[\sN_G(X,e)]$, which may not
have a convenient closed form.
The assumption that these graphs are forests is essentially only used to approximate this expectation by
$2 \se p^{\se-1} n^{\sv-2}$, and it comes mainly from the application of Proposition \ref{prop:Eproduct}
below, which states that expectations of products of edge variables may be approximated by the product
of their expectations (without using the FKG inequality);
see Section \ref{sec:review_fluctuations_consequences} for more information.
The specific form of this coefficient is not used to conclude the CLT for $\sN_G(X)$ from the CLT for $\sE(X)$.
However, we \emph{do} need this specific form in the proof of the CLT for $\sE(X)$ itself, as using
Theorem \ref{thm:stein} requires that various quantities cancel precisely.

In any case, a similar strategy also leads to a local version of Proposition \ref{prop:global_hajek_intro},
which is used to prove Theorem \ref{thm:local_clt_intro}.
In order to state the local version, for $\rho \in \cV$ and $v \in [n]$, we denote by $\sN_G^{\rho \to v}(x)$
the number of homomorphisms of $G$ in $x$ which map $\rho$ to $v$.
We also denote by $\sd_\rho$ the degree of $\rho$ in $G$.

\begin{proposition}
\label{prop:local_hajek_intro}
Under Assumption \ref{assumption},
define
\begin{equation}
\label{eq:hatngrhovdef}
    \hsN_G^{\rho \to v}(x) \coloneqq \sN_G^{\rho \to v}(x) - \sd_\rho p^{\se-1} n^{\sv-2} \deg_v(x).
\end{equation}
Then, for all $\zeta \in (0, \fr{3}{4})$, there is some $c > 0$ such that
\begin{equation}
    \P \left[
        \left|
            \hsN_G^{\rho \to v}(X) - \E \left[ \hsN_G^{\rho \to v}(X) \right]
        \right| > n^{\sv - 1.75 + \zeta}
    \right] \leq \Exp{- c n^{\fr{4}{3} \zeta}}.
\end{equation}
In particular, for all $\xi > 0$ we have
\begin{equation}
    \Var \left[
        \hsN_G^{\rho \to v}(X)
    \right] \lesssim n^{2\sv - 3.5 + \xi}.
\end{equation}
\end{proposition}

We remark that this local case is somewhat more technical than the global case, and requires
an extra input in the form of a bound on the \emph{homogeneity of the mixing}.
More precisely, we will bound the expected symmetric difference, in the
neighborhood of a vertex, between two monotonically-coupled Glauber dynamics chains which started
in agreement around that vertex.

In order to state this result precisely, let us introduce a bit of notation.
First, for a vertex $v \in [n]$, we define the \emph{local Hamming distance}
\begin{equation}
\label{eq:dlocdef}
    \dloc(x,x') \coloneqq \sum_{e \ni v} |x(e) - x'(e)|,
\end{equation}
which is the number of edges around $v$ at which the graphs $x$ and $x'$ differ.
Let us denote by $(X_t^x)$ the ERGM Glauber dynamics started at $X_0^x = x$,
and conditioned on the metastable well around $p \in U_\beta$.
This will be defined more precisely in Section \ref{sec:review_fluctuations_dynamical} below.
Additionally, for a graph $x$ and an edge $e \in \edgeset$ (the edge set
of the complete graph on $[n]$), we set $x^{\oplus e}$ to be
the same graph but with the status of $e$ (whether it is present or absent) swapped.
Finally, the set $\Lambda^*$ will be defined in Proposition \ref{prop:goodset} below, but for now
it should just be thought of as some set of graphs which is large under the ERGM measure
conditioned on the metastable well, i.e.\ $\P[X \in \Lambda^*] \geq 1 - e^{-\Omega(n)}$.

\begin{proposition}
\label{prop:dloc_uniformbound_intro}
Let $(X_t^x)$ denote the conditioned ERGM Glauber dynamics as above for some $p \in U_\beta$
(we do not require the ``phase-uniqueness-or-forest'' assumption for this result).
There is some $C > 0$ such that the following holds.
Let $v \in [n]$ and $e \in \edgeset$ with $v \notin e$.
Suppose that $x \in \Lambda^*$, and set $x' = x^{\oplus e}$.
Then for all $t \geq 0$, we have
\begin{equation}
    \E\left[\dloc(X_t^x,X_t^{x'})\right] \leq \fr{C}{n}.
\end{equation}
\end{proposition}

In some sense, Proposition \ref{prop:dloc_uniformbound_intro}
shows that the chain cannot become ``too unmixed'' in a small local region while it is mixing globally,
i.e.\ the \emph{global} Hamming distance $\dh(X_t^x, X_t^{x'})$ (which is just the sum of $|x(e) - x'(e)|$
over \emph{all} $e \in \edgeset$) is decreasing exponentially on average.
This local result is derived from previous rapid mixing results via a recurrence relation;
see Section \ref{sec:inputs_dloc} for more specifics.
\subsection{Outline of the paper}
\label{sec:intro_outline}

We begin in Section \ref{sec:review} with a survey of some background material which is crucial for studying ERGMs,
both at the macroscopic scale in Section \ref{sec:review_leadingorder}, and at the microscopic fluctuation
scale in Section \ref{sec:review_fluctuations}.
This includes key results like the large deviations principle of \cite{chatterjee2013estimating}
as well as the concentration inequality for Lipschitz observables of \cite{winstein2025concentration}.

In Section \ref{sec:inputs}, we introduce the main tools which we will use to prove Propositions
\ref{prop:global_hajek_intro} and \ref{prop:local_hajek_intro}.
First and foremost, Section \ref{sec:inputs_concentration} contains the general concentration inequality
due to  \cite{barbour2022long} derived from the rapid mixing of a Markov chain.
We collect necessary mixing inputs from previous works in Section \ref{sec:inputs_mixing},
and use them together with the result of \cite{barbour2022long}
to prove a basic concentration inequality in Section \ref{sec:inputs_basic}.
In addition, we prove Proposition \ref{prop:dloc_uniformbound_intro} in Section \ref{sec:inputs_dloc},
providing a sort of \emph{homogeneity} of the mixing.

In Section \ref{sec:hajek} we prove Propositions \ref{prop:global_hajek_intro} and \ref{prop:local_hajek_intro},
in Sections \ref{sec:hajek_global} and \ref{sec:hajek_local} respectively; these propositions show that
the fluctuations of subgraph counts are controlled by those of the overall edge count.
We also prove Corollaries \ref{cor:global_subgraph_clt_intro} and \ref{cor:local_subgraph_clt_intro},
the CLTs for subgraph counts, in Section \ref{sec:hajek_subgraphs} under the assumption
of Theorems \ref{thm:global_clt_intro} and \ref{thm:local_clt_intro}, the CLTs for edge count $\sE(X)$
and vertex degree $\deg_v(X)$ respectively.

Finally, in Section \ref{sec:stein} we prove Theorems \ref{thm:global_clt_intro} and \ref{thm:local_clt_intro},
in Sections \ref{sec:stein_global} and \ref{sec:stein_local} respectively.
This makes use of the aforementioned application of Stein's method due to \cite{fang2024normal}, which
is introduced in Section \ref{sec:stein_statement}.
We will need to slightly modify the Hamiltonian for this to work, and that will be discussed in
Section \ref{sec:stein_hamiltonian}.
\subsection{Acknowledgements}
\label{sec:intro_acknowledgements}

I was partially supported by the NSF Graduate Research Fellowship grant DGE 2146752.
I would like to thank my advisor, Shirshendu Ganguly, for suggesting this project.
I would also like to thank Xiao Fang, Victor Ginsburg, Kaihao Jing, Sumit Mukherjee,
Nathan Ross, and Rikhav Shah for helpful feedback on an earlier draft of this article.
\section{Review of the exponential random graph model}
\label{sec:review}

We now give a brief review of some basic tools from previous works which are vital for understanding
ERGMs and will be useful for the proofs to follow.
First, in Section \ref{sec:review_leadingorder}, we present results about the macroscopic leading-order behavior of ERGMs.
Then, in Section \ref{sec:review_fluctuations}, we turn to the microscopic fluctuation behavior.

\subsection{Leading-order behavior}
\label{sec:review_leadingorder}

To understand the leading-order behavior of samples from an ERGM, it is helpful to have an infinitary
object to which one may compare finite graphs.
This is the purview of \emph{graph limit theory}, and in particular we will need to introduce \emph{graphons}.
Presently we provide only a brief overview in Section \ref{sec:review_leadingorder_graphons},
referring the reader to \cite{lovasz2012large} for a thorough presentation of the subject.
Then, in Section \ref{sec:review_leadingorder_ldp}, we state the large deviations principle of
\cite{chatterjee2013estimating} and discuss some of its consequences.

\subsubsection{Graphons and the cut distance}
\label{sec:review_leadingorder_graphons}

In short, a graphon is a symmetric measurable function $W : [0,1]^2 \to [0,1]$.
This should be thought of as a generalization of an adjacency matrix: specifically, one may construct a
graphon from any $n$-vertex graph $x$ with adjacency matrix $(a_{ij})_{i,j=1}^n$ by setting
\begin{equation}
    W_x(s,t) = \sum_{i,j=1}^n a_{ij} \ind{\fr{i-1}{n} < s < \fr{i}{n}} \ind{\fr{j-1}{n} < t < \fr{j}{n}}.
\end{equation}
In pictorial language, we convert every entry of the adjacency matrix to a small subsquare of $[0,1]^2$ with
the same value under $W_x$.
Note however that this assignment is not injective; for instance, every complete graph has the same
graphon representation (ignoring Lebesgue-measure-zero subsets of $[0,1]^2$).

The natural distance on the space of graphons is the \emph{cut distance} $\db$, defined as follows:
\begin{equation}
\label{eq:dbdef}
    \db(W,W') \coloneqq \inf_{\substack{\sigma : [0,1] \to [0,1] \\ \text{measure-preserving} \\ \text{bijection}}}
    \sup_{S,T \sse [0,1]} \left| \int_S \int_T \; W(\sigma(s),\sigma(t)) - W'(s,t) \; \,ds \,dt \right|.
\end{equation}
In particular, if there is a measure-preserving bijection $\sigma : [0,1] \to [0,1]$ such that
$W(\sigma(s), \sigma(t)) = W'(s,t)$ a.e., then $\db(W,W') = 0$ and we consider these graphons to be isomorphic
so that we obtain a metric space.

The formula \eqref{eq:dbdef} is not important for most of the present article, but one should be aware that
convergence under $\db$ is equivalent to convergence of all subgraph densities:
\begin{equation}
\label{eq:equivalence}
    \db(W_n,W) \to 0
    \qquad \Longleftrightarrow \qquad
    \st(G,W_n) \to \st(G,W) \quad \text{for all finite graphs } G.
\end{equation}
Here, we define the subgraph density of a graphon as follows: for a graph $G = (\cV, \cE)$,
\begin{equation}
\label{eq:tgwdef}
    \st(G,W) \coloneqq \int_{(s_\rho) \in [0,1]^\cV}
    \prod_{\{ \rho_1, \rho_2 \} \in \cE} W(s_{\rho_1}, s_{\rho_2})
    \; \prod_{\rho \in \cV} ds_\rho.
\end{equation}
Note that this extends the notion of subgraph density for graphs $x$ on $n$ vertices: $\st(G,x) = \st(G,W_x)$.
We also remark that the equivalence \eqref{eq:equivalence} may be made quantitative
\cite{lovasz2006limits,borgs2008convergent}, although this is not relevant for our purposes in this article.

\subsubsection{Large deviations principle}
\label{sec:review_leadingorder_ldp}

As mentioned in the introduction, in general an ERGM is defined by the selection of a Hamiltonian $\sH$,
which is a continuous function on the space of graphons; we evaluate such a function at an $n$-vertex graph $x$
by setting $\sH(x) = \sH(W_x)$.
Given a general Hamiltonian $\sH$, the $n$-vertex random graph $\X$ has the ERGM distribution if \eqref{eq:gibbsmeasure}
holds.
Such random graphs were considered by \cite{chatterjee2013estimating}, who proved a
\emph{large deviations principle} explaining their first-order behavior, i.e.\ how they behave when considered
as graphons under the cut distance.
To state their result, let us define the \emph{graphon entropy functional}
\begin{equation}
\label{eq:entropy}
    I(W) \coloneqq \int_{[0,1]^2} I(W(s,t)) \,ds \,dt,
\end{equation}
where $I(q)$ for numbers $q \in [0,1]$ was defined just after the definition of $L_\beta(q)$ in \eqref{eq:Ldef}.
Additionally, define $\cM_\sH$ to be the set of graphons which maximize the functional $\sH(W) - I(W)$,
which is the graphon version of $L_\beta$.
We are now ready to state the main result of \cite{chatterjee2013estimating}, which states that $n$-vertex
samples $X$ from an ERGM are close to $\cM_\sH$; we denote by $\db(W,\cM_\sH)$ the infimal cut distance
from $W$ to any graphon in $\cM_\sH$.

\begin{theorem}[Theorem 3.2 of \cite{chatterjee2013estimating}]
\label{thm:LDP_general}
For any $\eta > 0$, there are constants $c(\eta), C(\eta) > 0$ such that
\begin{equation}
    \P \left[ \db(\X, \cM_\sH) > \eta \right] \leq C(\eta) \Exp{-c(\eta) n^2}.
\end{equation}
\end{theorem}

Recall from Section \ref{sec:intro_setup}, in particular \eqref{eq:hdef}, that in the present work we
consider Hamiltonians that can be written as linear combinations of subgraph densities
\begin{equation}
    \sH(W) = \sum_{j=0}^K \beta_j \, \st(G_j,W),
\end{equation}
for fixed graphs $G_0, G_1, \dotsc, G_K$ (where $G_0$ is always a single edge)
and parameters $\beta = (\beta_0, \beta_1, \dotsc, \beta_K) \in \R \times \R_{\geq 0}^K$.
In this setting, specifically under the ferromagnetic assumption, the set $\cM_\sH$ simplifies greatly.
Namely, the optimizers of the functional $\sH(W) - I(W)$ are all \emph{constant graphons} $W_q$,
defined for $q \in [0,1]$ by $W_q(s,t) = q$.
When restricted to constant graphons, we have
\begin{equation}
    \sH(W_q) - I(W_q) = L_\beta(q),
\end{equation}
recalling the definition of $L_\beta(q)$ from \eqref{eq:Ldef}
and noting that $\st(G,W_q) = q^{\se}$ for a graph $G$ with $\se$ edges.
Moreover, as also defined in Section \ref{sec:intro_setup}, we denote by $M_\beta$ the finite set of global
maximizers of $L_\beta$ in $(0,1)$.

\begin{theorem}[Theorem 4.1 of \cite{chatterjee2013estimating}]
\label{thm:LDP_ferromagnetic} 
In the above context, $\cM_\sH = \{ W_p : p \in M_\beta \}$.
\end{theorem}

Note that constant graphons $W_q$ are limits of Erd\H{o}s--R\'enyi graphs, in the sense that if
$Y \sim \cG(n,q)$, then $\db(Y, W_q) \pto 0$ as $n \to \infty$.
Thus, if $L_\beta$ has a unique maximizer, i.e.\ $M_\beta = \{p\}$,
then Theorems \ref{thm:LDP_general} and \ref{thm:LDP_ferromagnetic}
(and the equivalence \eqref{eq:equivalence}) show that a sample $\X$ from an ERGM with
Hamiltonian \eqref{eq:hdef} behave similarly under $\db$ as a sample $Y$ from $\cG(n,p)$, in the sense that
they obey the same \emph{laws of large numbers} for homomorphism counts.
Namely,
\begin{equation}
    \fr{\sN_G(\X)}{n^\sv} \pto p^\se
    \qquad \text{and} \qquad
    \fr{\sN_G(Y)}{n^\sv} \pto p^\se
\end{equation}
for any graph $G$ with $\sv$ vertices and $\se$ edges. 

If $|M_\beta| > 1$, the ERGM may be decomposed as a mixture of measures which satisfy such laws of large numbers,
one for each $p \in M_\beta$.
It is easiest to define such a mixture by simply conditioning on a cut-distance ball.
To this end, for any $p \in [0,1]$ and $\eta > 0$ we define
\begin{equation}
\label{eq:balldef}
    \pball \coloneqq \{ W : \db(W,W_p) < \eta \},
\end{equation}
and we define the measure $\condmeas$ by
\begin{equation}
\label{eq:condmeas}
    X \sim \condmeas \qquad \Leftrightarrow \qquad \P[X = x] \propto e^{n^2 \sH(x)} \ind{W_x \in \pball}.
\end{equation}
With the equivalence \eqref{eq:equivalence}, Theorems \ref{thm:LDP_general} and \ref{thm:LDP_ferromagnetic}
imply that for all $p \in M_\beta$ and all small enough $\eta$, samples $X \sim \condmeas$ satisfy the
same laws of large numbers for homomorphism counts as $Y \sim \cG(n,p)$.
Moreover, these theorems imply the existence of coefficients $c_p^\eta(n)$ such that the mixture
\begin{equation}
\label{eq:mixture}
    \sum_{p \in M_\beta} c_p^\eta(n) \, \condmeas
\end{equation}
has total variation distance $\Exp{-\Omega(n^2)}$ to the full ERGM measure itself.
Note that these theorems do not prevent the coefficients $c_p^\eta(n)$ from decaying with $n$,
but they may not decay faster than $\Exp{-\Omega(n^2)}$.

This is essentially the end of the story about the macroscopic behavior of ERGMs in our context, apart from
the important question of the relative weights $c_p^\eta(n)$ of measures in the mixture \eqref{eq:mixture},
which remains open in general.
Some progress towards this result was provided by \cite{mukherjee2023statistics} which proves that
$|M_\beta| = 2$ in the phase coexistence case of the $2$-star ERGM, and that $c_p^\eta(n) \to \fr{1}{2}$
for both optimal densities $p \in M_\beta$.
Additionally, \cite{bianchi2024limit} proved that $|M_\beta| = 2$ for the edge-triangle model in the
phase coexistence case and gave conjectural formulas for the limits of $c_p^\eta(n)$
which are nonzero, i.e.\ they conjecture that both optimal densities $p \in M_\beta$ have comparable weight
under the full edge-triangle measure.

In any case, the above laws of large numbers say almost nothing about the microscopic fluctuations.
To proceed with our goals of proving quantitative CLTs for the edge count and subgraph/homomorphism counts, we will need
more tools to understand this microscopic behavior.
\subsection{Microscopic fluctuations}
\label{sec:review_fluctuations}

The tools we will need for understanding the microscopic fluctuations of ERGM depend on the
\emph{dynamical} behavior of the model under \emph{Glauber dynamics}.
This is a Markov chain $(\X_t)$ on $n$-vertex graphs which has the ERGM measure as its stationary distribution.
As with a single sample, we use a tilde to indicate that we are working with the unconditioned measure;
later, $(X_t)$ will denote the Glauber dynamics with respect to the conditioned measure $\condmeas$.
We first describe its behavior in Section \ref{sec:review_fluctuations_dynamical}, then introduce a
concentration inequality derived from this behavior in Section \ref{sec:review_fluctuations_concentration},
and finally state some consequences of this concentration inequality in
Section \ref{sec:review_fluctuations_consequences}.

\subsubsection{Dynamical behavior}
\label{sec:review_fluctuations_dynamical}

Under Glauber dynamics, at each time step a single potential edge $e \in \edgeset$ of $\X_t$ is chosen uniformly at
random to be resampled.
Once the edge $e$ has been chosen, $\X_{t+1}$ will be either $\X_t^{+e}$ or $\X_t^{-e}$; here, $x^{\pm e}$
denotes the same graph as $x$ but with the edge $e$ fixed to be either \emph{present} or \emph{absent},
respectively.

In order that the ERGM measure be reversible with respect to these dynamics, we must have the following
probabilities:
\begin{equation}
	\Pcond{\X_{t+1} = \X_t^{\pm e}}{e \text{ was chosen}} = 
	\frac{e^{n^2 \sH(\X_t^{\pm e})}}{e^{n^2 \sH(\X_t^{-e})} + e^{n^2 \sH(\X_t^{+e})}}.
\end{equation}
Factoring out $e^{n^2 \sH(\X_t^{-e})}$ from the numerator and denominator of this ratio with
the $\pm$ set to $+$, we obtain
\begin{equation}
\label{eq:glauberratio}
	\Pcond{\X_{t+1} = \X_t^{+e}}{e \text{ was chosen}} = 
	\frac{e^{n^2 \partial_e \sH(\X_t)}}{1 + e^{n^2 \partial_e \sH(\X_t)}},
\end{equation}
where in general we introduce the notation
\begin{equation}
\label{eq:partialdef}
	\partial_e f(x) \coloneqq f(x^{+e}) - f(x^{-e}).
\end{equation}

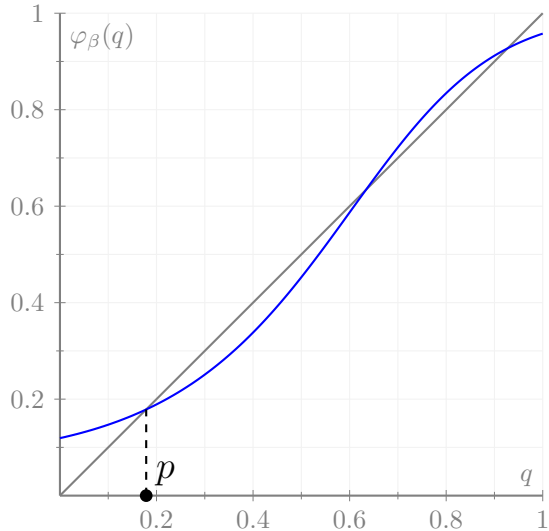
\begin{figure}
\centering
\begin{tikzpicture}
  \begin{axis}[
    width=8cm,
    height=8cm,
    axis lines=middle,
    xlabel={$q$},
    ylabel={$\phi_\beta(q)$},
    domain=0:1,
    samples=500,
    thick,
    grid=both,
    minor tick num=1,
    grid style={gray!10},
    enlargelimits=false,
    clip=false,
    scaled ticks=false,
    axis line style={-,gray},
    label style={gray},
    tick label style={gray},
  ]

  \addplot[gray, thick] {x};
  \addplot[blue, thick] {phi(x)};

  \addplot[dashed, black] coordinates {(\pstar,0) (\pstar,{phi(\pstar)})};
  \addplot[mark=*, black] coordinates {(\pstar,0)};

  \node[above right, black] at (axis cs:\pstar,0) {\Large $p$};
  \end{axis}
\end{tikzpicture}
\caption{Plot of $\phi_\beta(q)$, as defined in \eqref{eq:phidef1}, for the ERGM with the
same specifying graphs and parameters as in Figure \ref{fig:lbeta}.
Compare with that figure, which depicts a plot of $L_\beta(q)$, and observe
that the attracting fixed points of $\phi_\beta$ are the same as the local maxima of $L_\beta$.
In particular, the global maximum $p$ of $L_\beta$ is an attracting fixed point,
but there may be other attracting fixed points of $\phi_\beta$ which are not \emph{global} maxima of $L_\beta$.}
\label{fig:phibeta}
\end{figure}

Analysis of the ratio \eqref{eq:glauberratio} leads to the following heuristic which will be
explained more rigorously in Section \ref{sec:inputs_mixing} below.
If the density of $\X_t$ is approximately $q$ in a sense made precise in that section, then
when an edge $e$ is chosen to be resampled, the probability that $\X_{t+1} = \X_t^{+e}$ is approximately
\begin{equation}
\label{eq:phidef1}
	\phi_\beta(q) \coloneqq \fr{e^{2 \sH'(q)}}{1 + e^{2 \sH'(q)}},
\end{equation}
where we write $\sH'(q)$ for the derivative in $q$ of $\sH(q) = \sH(W_q)$ (recalling from \eqref{eq:hdef}
that this depends on $\beta$ as well).
See Figure \ref{fig:phibeta} for a sample plot of $\phi_\beta(q)$.
Thus, heuristically, over long periods of time the density of $\X_t$ evolves akin to the dynamical system
on $[0,1]$ obtained by repeated application of $\phi_\beta$.
In particular, the density should drift towards an \emph{attracting fixed point} of $\phi_\beta$,
which is a value $p \in [0,1]$ such that
\begin{equation}
\label{eq:attractingfixedpoint}
	\phi_\beta(p) = p
	\qquad \text{and} \qquad
	\phi_\beta'(p) < 1.
\end{equation}
As can easily be checked, the condition \eqref{eq:attractingfixedpoint} is equivalent to
\begin{equation}
	L_\beta'(p) = 0
	\qquad \text{and} \qquad
	L_\beta''(p) < 0,
\end{equation}
where $L_\beta$ was defined in \eqref{eq:Ldef} and explained further in
Section \ref{sec:review_leadingorder} above.
In other words $p$ is an attracting fixed point for $\phi_\beta$ if and only if it is a strictly concave
\emph{local} maximizer for $L_\beta$.

Recall from Section \ref{sec:review_leadingorder} above that the density of samples from the ERGM measure
are close to $p$ for some \emph{global} maximizer of $L_\beta$.
However, there might be local maximizers (i.e.\ attracting fixed points of $\phi_\beta$)
which are not global maximizers of $L_\beta$.
The existence of multiple attracting fixed points of $\phi_\beta$ leads to \emph{bottlenecks} in the Glauber
dynamics, and it may take an extremely long time for the graph to transition between different
\emph{metastable wells} $\pball$ (defined in \eqref{eq:balldef}) for different attracting fixed points $p$.

For this reason, as already alluded to multiple times, we separate the parameters into multiple regimes of consideration.
First, the \emph{subcritical regime} of parameters consists of all $\beta$ for which $\phi_\beta$ has a
\emph{unique} fixed point which is attracting.
Next, the \emph{supercritical regime} consists of $\beta$ for which $\phi_\beta$ has \emph{multiple}
fixed points and each one which is a global maximizer of $L_\beta$ is attracting.
In all other cases, we say that $\beta$ is \emph{critical}.
The \emph{Dobrushin regime} is a subset of the subcritical regime, defined by the condition $\sH''(1) < 2$.
This condition turns out to imply \emph{uniform contraction} with respect to the Hamming distance under the monotone coupling,
which immediately yields rapid mixing via a simple path coupling argument \cite{bubley1997path}.

In \cite{bhamidi2008mixing} it was proved that the Glauber dynamics mixes rapidly in the \emph{full} subcritical
phase, even beyond the Dobrushin regime, and they also proved that Glauber dynamics mixes exponentially slowly
in the supercritical case from worst-case starting configurations, due to the bottleneck behavior mentioned above.
More recently, however, \cite{bresler2024metastable} showed that in the supercritical case Glauber dynamics exhibits
\emph{metastable mixing}, i.e.\ rapid mixing \emph{within} a metastable well $\pball$ with a \emph{warm start}
$X_0 \sim \cG(n,p)$ (here we use $(X_t)$ for the \emph{conditioned} Glauber dynamics, which will be
defined precisely in Section \ref{sec:inputs_mixing}).
In \cite{bresler2024metastable} this is only stated for metastable wells $\pball$ around attracting fixed points
$p$ of $\phi_\beta$ which are also \emph{global} maxima of $L_\beta$, i.e.\ \emph{large} metastable wells.
We remind the reader that $U_\beta \sse M_\beta$ is our notation for the set of such values of $p$.
While it is plausible that metastable mixing also occurs for metastable wells around attracting fixed points
which are only \emph{local} maxima of $L_\beta$, we restrict our attention to $p \in U_\beta$ for the present work.

\subsubsection{Gaussian concentration for Lipschitz observables}
\label{sec:review_fluctuations_concentration}

The mixing results of \cite{bhamidi2008mixing,bresler2024metastable} lead to a variety of important static
consequences for ERGM measures and the mixture measures $\condmeas$ defined in \eqref{eq:condmeas},
which are conditioned on the cut-distance ball $\pball$ defined in \eqref{eq:balldef}.
Chief among them is a Gaussian concentration inequality for Lipschitz observables of $X \sim \condmeas$.
This was first proved in the subcritical regime by \cite{ganguly2024sub}, using the rapid mixing result of
\cite{bhamidi2008mixing} as input to a technique for proving such concentration inequalities originally
developed by \cite{chatterjee2005concentration}.
Note that in the subcritical case, since we have phase uniqueness, the result may be stated simply for the
full ERGM measure as opposed to the conditioned measure $\condmeas$.
Later, this concentration inequality was extended by \cite{winstein2025concentration} to the measures
$\condmeas$ in the supercritical regime,
using the metastable mixing of \cite{bresler2024metastable} as input for a modification of
the technique of \cite{chatterjee2005concentration}, as well as a result of \cite{barbour2022long}.

In order to state this result, we introduce the notation of Lipschitz vectors of functions on $n$-vertex graphs.
We say that $\cL \in \R^{\binom{[n]}{2}}$ is a \emph{Lipschitz vector} of a real-valued function $f$ on
$n$-vertex graphs if $|f(x^{+e}) - f(x^{-e})| \leq \L_e$ for all $e \in \binom{[n]}{2}$ and all $n$-vertex
graphs $x$ (recall that $x^{\pm e}$ was defined at the beginning of Section \ref{sec:review_fluctuations_dynamical}).
Furthermore, in the following result, we fix some $p \in U_\beta$.

\begin{theorem}[Theorem 5.1 of \cite{winstein2025concentration}]
\label{thm:lipschitzconcentration}
Let $p \in U_\beta$ (we do not require the full Assumption \ref{assumption}).
There is some $\eta > 0$ as well as constants $C, c > 0$
such that the following holds for $X \sim \condmeas$.
For any function $f$ with Lipschitz vector $\cL$ and $\lambda \geq 0$ with
$\lambda \leq c \| \L \|_1$,
\begin{equation}
    \P[| f(X) - \E[f(X)] | > \lambda]
    \leq 2 \Exp{ - c \cdot \max \left\{
      \fr{\lambda^2}{n^2 \| \L \|_\infty^2},
      \fr{\lambda^2}{\| \L \|_1 \| \L \|_\infty} - \Exp{\fr{C \lambda}{\| \L \|_\infty} - c n}
      \right\} }
      + e^{-c n}.
\end{equation}
\end{theorem}

A few remarks are in order.
First, the initial expression in the maximum gives concentration with a
variance proxy of $n^2 \| \cL \|_\infty^2$, applicable to all Lipschitz observables.
This is optimal in the case of ``global'' observables such as the overall edge count
(which should have variance of order $n^2$), but is suboptimal in the case of ``local''
observables such as the degree of a vertex (which should have variance of order $n$).
The second expression in the maximum gives a better bound in the case of local observables,
but cannot apply to global observables due to the subtracted exponential term; this is why
both expressions are needed.

Next, the external $e^{-c n}$ term is necessary due to the potential presence of small metastable regions
with atypical behavior \emph{within} the large metastable well $\pball$, which was already mentioned
in Section \ref{sec:intro_iop_props}.
This will not affect our proof as all we require is an upper bound that is exponentially small in certain
small powers of $n$.
For more information about this behavior, we again refer to \cite[Theorem 3.3]{bresler2024metastable}.

Finally, although Theorem \ref{thm:lipschitzconcentration} gives a quick and easy method of proving
concentration in many cases which will be relevant to us, the Lipschitz vector often does not capture
the true scale of fluctuations as it essentially considers a worst-case change in $f$.
Instead, to make use of the average-case change when an edge is flipped, we will use a different
concentration inequality derived from a general result of \cite{barbour2022long}, which will be
discussed in Section \ref{sec:inputs_concentration} below.
The use of this inequality requires much more careful analysis, however, and even still it
yields worse bounds than Theorem \ref{thm:lipschitzconcentration} in some cases.
As such, we will use both Theorem \ref{thm:lipschitzconcentration} as well as this alternate method
at different points during the proof.

\subsubsection{Consequences of concentration}
\label{sec:review_fluctuations_consequences}

Theorem \ref{thm:lipschitzconcentration} leads to a variety of consequences which will be useful in what
follows.
In all remaining results of this section, we fix some $p \in U_\beta$ and let $X \sim \condmeas$ for some
$\eta > 0$ small enough for the result of Theorem \ref{thm:lipschitzconcentration} to hold.

First, we have a bound on the covariance between edges in $X$.
For an edge $e \in \edgeset$, we introduce the notation $X(e)$, used throughout the remainder of the article,
for the \emph{edge indicator}, i.e.\ $X(e) = 1$ if the edge $e$ is present in $X$, and $X(e) = 0$ otherwise.

\begin{lemma}[Lemma 6.1 of \cite{winstein2025concentration}]
\label{lem:cov}
Let $p \in U_\beta$ (we do not require the full Assumption \ref{assumption}).
For any distinct $e, e' \in \edgeset$, we have
\begin{equation}
	\left| \Cov[X(e), X(e')] \right| \lesssim \fr{1}{n}.
\end{equation}
Moreover, if $e$ and $e'$ do not share a vertex, then the right-hand side
above can be improved to $\fr{1}{n^2}$.
\end{lemma}

This can be extended to multilinear moments of edge indicator variables.
In the phase uniqueness case (i.e.\ $M_\beta = U_\beta = \{p\}$), since $\condmeas$ is essentially the same as
the full ERGM measure, this follows from a result of \cite{newman1980normal} which uses the FKG inequality.
However, in the phase coexistence case this tool is absent, as the FKG inequality only holds for the full
ERGM measure itself and is quite brittle in this regard.
By an alternative strategy, we may still proceed but only if the multilinear moment consists of edge
variables which form a \emph{forest}, i.e.\ there may not be any cycles.
It remains an interesting open question whether this condition may be removed in the phase coexistence case.
Such a result would immediately expand the validity of every result in the present article, removing
the ``phase-coexistence or forest'' assumption throughout.

\begin{proposition}[Proposition 6.2 of \cite{winstein2025concentration}]
\label{prop:Eproduct}
For any fixed $k$, let $e_1, \dotsc, e_k \in \binom{[n]}{2}$ be distinct potential edges and suppose that
either of the following two conditions hold:
\begin{enumerate}[label=(\alph*)]
	\item The density $p$ is the unique global maximizer of $L_\beta$, i.e.\ $M_\beta = U_\beta = \{p\}$.
	\item The edges $e_1, \dotsc, e_k$ form a forest in $K_n$, the complete graph on $[n]$.
\end{enumerate}
Then we have
\begin{equation}
	\left| \E\left[ \prod_{j=1}^k X(e_j) \right]
	- \E[X(e)]^k \right| \lesssim \fr{1}{n},
\end{equation}
where $e$ is an arbitrary edge in $\binom{[n]}{2}$.
\end{proposition}

Another important consequence of Theorem \ref{thm:lipschitzconcentration}, which also uses
Proposition \ref{prop:Eproduct} as an input, is a bound on the marginal probability that an edge
$e \in \edgeset$ is present in $X$, or in other words $\E[X(e)]$.

\begin{proposition}[Proposition 6.5 of \cite{winstein2025concentration}]
\label{prop:marginal}  
If $M_\beta = U_\beta = \{p\}$ or all graphs $G_0, G_1, \dotsc, G_K$ in the ERGM specification are forests
(i.e.\ we make Assumption \ref{assumption}), then
\begin{equation}
    \left|\E[X(e)] - p \right| \lesssim \sqrt{\fr{\log n}{n}}
\end{equation}
for any $e \in \edgeset$.
\end{proposition}

We remark that the error bound above was recently improved to $\fr{1}{n}$ in the subcritical parameter regime
by \cite[Proposition 1.1]{fang2025conditionalcentrallimittheorems}, resolving
\cite[Conjecture 2.1]{winstein2025concentration} in the subcritical regime.
Moreover, we expect that a similar bound also holds in the supercritical regime, without the assumption
that the graphs $G_0, G_1, \dotsc, G_K$ are forests in the phase coexistence case, since this assumption
comes from the use of Proposition \ref{prop:Eproduct}, which we do not believe to be necessary for such a result in general,
e.g.\ by following a strategy similar to \cite{reinert2019approximating} with the help of tools developed by
\cite{winstein2025concentration}.

Nevertheless, we will use Proposition \ref{prop:marginal} as stated since it holds throughout all parameter regimes,
and we already have to use Proposition \ref{prop:Eproduct} elsewhere in our proofs, necessitating the
``phase-uniqueness or forest'' assumption anyway.
Additionally, the stronger bound provided by \cite[Proposition 1.1]{fang2025conditionalcentrallimittheorems}
would not improve the quantitative error bounds we obtain in Theorems \ref{thm:global_clt_intro}
or \ref{thm:local_clt_intro} (although as mentioned it is not clear if these bounds are sharp).
\section{Mixing and concentration}
\label{sec:inputs}

In this section we collect various technical inputs needed for the proof in the sequel.
As mentioned in Section \ref{sec:review_fluctuations_concentration} above, we will need a general
concentration inequality from \cite{barbour2022long}, and this is discussed in Section \ref{sec:inputs_concentration}.
To apply this concentration inequality to ferromagnetic ERGMs, we will need some details of the rapid mixing
and metastable mixing results of \cite{bhamidi2008mixing} and \cite{bresler2024metastable} respectively;
these are presented in Section \ref{sec:inputs_mixing}.
In Section \ref{sec:inputs_basic}, we provide a basic application of these mixing results with the result
of \cite{barbour2022long}, which yields a basic concentration inequality which will nonetheless be useful
at a certain point in the proof of Theorem \ref{thm:global_clt_intro}, and may be viewed as a warm-up to
the broader proof strategy.
Finally, for our local results like Theorem \ref{thm:local_clt_intro},
Corollary \ref{cor:local_subgraph_clt_intro}, and Proposition \ref{prop:local_hajek_intro}, we need an
additional technical input, bounding the discrepancy within a neighborhood of a single vertex between
two coupled copies of the Glauber dynamics chain.
This provides a form of \emph{homogeneity} for the mixing, showing that no small local area can become
``too unmixed'' while the overall chain mixes.
To the best of our knowledge, such a result has not appeared in the literature before and so we provide
a complete proof in Section \ref{sec:inputs_dloc}.

\subsection{Concentration via mixing}
\label{sec:inputs_concentration}

As already hinted at in Sections \ref{sec:intro_iop} and \ref{sec:review_fluctuations_concentration},
it is possible to derive concentration inequalities for a variable $f(X)$ if there is a Markov chain $(X_t)$
for which the distribution of $X$ is stationary, such that the chain mixes rapidly and $f(X_t)$ does not change
by much at each step.
Intuitively, this is simply because different ``runs'' of the algorithm sampling $X$ by running the chain
(from a fixed initialization) up to the mixing time $T$ must yield similar values of $f(X_T)$.
The idea of analyzing a sampling algorithm to obtain concentration results was already present in the classical
Azuma/Hoeffding inequality for martingales, and was made more explicit in the famous OSSS inequality for Boolean
functions \cite{o2005every}.

Applying this idea to Markov chains
appears to have first been considered by \cite{chatterjee2005concentration}, which is used by \cite{ganguly2024sub}
in the subcritical regime for ferromagnetic ERGMs, and which is modified by \cite{winstein2025concentration}
to extend the results of \cite{ganguly2024sub} to the supercritical regime and prove
Theorem \ref{thm:lipschitzconcentration} above.
In a more recent and separate development due to \cite{barbour2022long}, a different proof technique leads to
a subtly different statement than what appears in \cite{chatterjee2005concentration} or
\cite{winstein2025concentration}, with various benefits and drawbacks.

In what follows, we consider a discrete-time ergodic Markov chain with a finite state space $\Omega$ and
transition kernel $P$.
We use $(X_t^x)$ to denote the chain started at $X_0^x = x$.
Additionally, we use $X$ to denote a sample from the stationary distribution of the chain.
The following theorem from \cite{barbour2022long} has been updated to match our notation.
We remark that the finiteness of the state space was not assumed in \cite{barbour2022long}, only that
$\Omega$ is discrete; the only additional requirement in the case of infinite state space is that
$\E[f(X_t^x)]$ exists for all $t \in \N$ and $x \in \Omega$.

\begin{theorem}[Theorem 2.1 of \cite{barbour2022long}]
\label{thm:barbour}
Let $\Lambda \sse \Omega$ and $f : \Omega \to \R$ such that the following two conditions hold for some
positive numbers $V$ and $\Delta$, as well as some time $T \in \N$.
\begin{enumerate}[label=(\roman*)]
    \item
    \label{cond:variance}
    A bound on a proxy for the variance of $f$: for all $x \in \Lambda$,
    \begin{equation}
        \sum_{x' \in \Omega} P(x,x') \sum_{t=0}^{T-1} \left(
            \E[f(X_t^x)] - \E[f(X_t^{x'})] \right)^2 \leq V.
    \end{equation}
    \item
    \label{cond:smalljumps}
    A bound on the future expected difference, measured by $f$, between chains started at two adjacent initial states:
    for all $x \in \Lambda$, and $x' \in \Omega$ with $P(x,x') > 0$, and for all $t \geq 0$,
    \begin{equation}
        \left| \E[f(X_t^x)] - \E[f(X_t^{x'})] \right| \leq \Delta.
    \end{equation}
\end{enumerate}
Then, for any $x \in \Lambda$, we have
\begin{equation}
\label{eq:barbourconclusion}
    \P \left[
        X_t^x \in \Lambda \text{ for } 0 \leq t < T,
        \text{ and }
        |f(X_T^x) - \E[f(X_T^x)]| \geq \lambda
    \right] \leq 2 \Exp{- \fr{\lambda^2}{2 V + \fr{4}{3} \Delta \lambda}}.
\end{equation}
\end{theorem}

This result provides concentration for the chain at some large enough time $T$, but we would like to
have concentration under the stationary distribution itself.
Moreover, in the conclusion there appears the extraneous condition that $X_t^x \in \Lambda$ for all
$0 \leq t < T$, which we would like to avoid carrying through the proofs that follow.
Both of these issues may be remedied if the chain mixes rapidly enough and if the extraneous condition
is very likely.
As such, we introduce some extra hypotheses and state a straightforward corollary which allows
for more seamless application of Theorem \ref{thm:barbour} in the sequel.
In this corollary and below, we use $\dtv$ to denote the \emph{total variation distance} between
probability distributions or random variables.

\begin{corollary}
\label{cor:barbour}
Suppose there is a set $\Lambda$ and a function $f$ such that the conditions \ref{cond:variance} and \ref{cond:smalljumps} in
the statement of Theorem \ref{thm:barbour} hold for positive numbers $V$ and $\Delta$, as well as time $T$.
Suppose also that the following three additional conditions hold for some nonnegative numbers
$\epsilon, \delta$, and $M$, and some $z \in \Lambda$.
\begin{enumerate}[label=(\roman*)]
    \setcounter{enumi}{2}
    \item
    \label{cond:setlarge}
    The chain has a high probability of staying in $\Lambda$ for time $T$ when started from $z$:
    \begin{equation}
        \P[X_t^z \in \Lambda \text{ for } 0 \leq t < T] \geq 1 - \epsilon.
    \end{equation}

    \item
    \label{cond:mixing}
    The chain mixes rapidly with respect to $f$ when started from $z$:
    \begin{equation}
        \dtv(f(X_T^z), f(X)) \leq \delta.
    \end{equation}

    \item
    \label{cond:fbound}
    The function $f$ is bounded close to its mean: for all $x \in \Omega$,
    \begin{equation}
        |f(x) - \E[f(X)]| \leq M.
    \end{equation}
\end{enumerate}
Then for all $\lambda \geq 2 \delta M$ we have
\begin{equation}
\label{eq:barbourcor}
    \P[|f(X) - \E[f(X)]| > \lambda] \leq 2 \exp \left(
        - \fr{(\lambda-2\delta M)^2}{2V + \fr{4}{3} \Delta (\lambda - 2\delta M)}
    \right) + \epsilon + \delta.
\end{equation}
\end{corollary}

\begin{proof}[Proof of Corollary \ref{cor:barbour}]
First of all, conditions \ref{cond:mixing} and \ref{cond:fbound} imply that
\begin{equation}
    \left| \E[f(X_T^z)] - \E[f(X)] \right| \leq 2 \delta M,
\end{equation}
so if we have $\left| f(X_T^z) - \E[f(X)] \right| > \lambda$ then we must also have
\begin{equation}
\label{eq:lambdaminus2deltam}
    \left|
        f(X_T^z) - \E[f(X_T^z)]
    \right| > \lambda - 2 \delta M.
\end{equation}
By conditions \ref{cond:setlarge} and \ref{cond:mixing}, we may couple $X$ and $(X_t^z)$ such that the event
where $f(X) \neq f(X_T^z)$ or $X_t^z \notin \Lambda$ for some $t$ with $0 \leq t < T$ happens with probability
at most $\epsilon + \delta$.
Thus \eqref{eq:barbourcor} follows from \eqref{eq:barbourconclusion} with $\lambda$ replaced by $\lambda - 2 \delta M$.
\end{proof}

One should compare the statement of Corollary \ref{cor:barbour} with \cite[Theorem 4.2]{winstein2025concentration},
which was used to prove Theorem \ref{thm:lipschitzconcentration}, and is a modification of
\cite[Theorem 3.3]{chatterjee2005concentration} that has five hypotheses similar to Corollary \ref{cor:barbour},
as well as a similar conclusion.
The main advantage that \cite[Theorem 4.2]{winstein2025concentration} has over Corollary \ref{cor:barbour}
is that condition \ref{cond:variance} is replaced by the potentially milder condition
\begin{equation}
    \sum_{x' \in \Omega} P(x,x') |f(x) - f(x')| \sum_{t=0}^{T-1} \left| \E[f(X_t^x)] - \E[(X_t^{x'})] \right|
    \leq V.
\end{equation}
This is precisely what allows the scale of concentration in Theorem \ref{thm:lipschitzconcentration} to be
$\sqrt{\| \L \|_1 \| \L \|_\infty}$ rather than the potentially larger scale $\sqrt{\| \L \|_\infty^2 \binom{n}{2}}$, which is
all that can be guaranteed in general by Corollary \ref{cor:barbour},
as discussed in \cite[Section 7.2]{winstein2025concentration} and briefly mentioned in
Section \ref{sec:review_fluctuations_concentration} above.
Note also that a concentration inequality with rate $\| \cL \|_\infty n$ will be proved in Section \ref{sec:inputs_basic} below,
which should be seen as a warm-up to the use of \cite{barbour2022long} in the more technical proofs to follow.
The deficiency in this rate comes down to the fact that it is not always the case that
\begin{equation}
\label{eq:distancegrows}
    \left| \E[f(X_t^x)] - \E[f(X_t^{x'})] \right| \leq |f(x) - f(x')|,
\end{equation}
even for adjacent $x,x'$.
This is particularly apparent when the function $f$ does not depend on certain edge variables, meaning that
the right-hand side above is zero in many cases.

However, \cite[Theorem 4.2]{winstein2025concentration} has its own deficiency which is an extra term inside
the exponential of the upper bound in the conclusion.
This term takes the form
\begin{equation}
    \epsilon \cdot e^{\lambda M / V},
\end{equation}
which means that the bound is destroyed as soon as $\lambda \gg \fr{V}{M} \log \fr{1}{\epsilon}$.
This is what leads to the extra term inside the exponential in Theorem \ref{thm:lipschitzconcentration}.
Separately, this extra term makes \cite[Theorem 4.2]{winstein2025concentration} unsuitable for our present applications.

Instead, we will apply Corollary \ref{cor:barbour}, but we will need to prove an approximate version of
the inequality \eqref{eq:distancegrows} for the particular choices of $f$ that we will consider in the proof of
the \emph{local} results, i.e.\ Theorem \ref{thm:local_clt_intro}, Corollary \ref{cor:local_subgraph_clt_intro}, 
and Proposition \ref{prop:local_hajek_intro}.
Indeed, without some version of this inequality the concentration afforded by Corollary \ref{cor:barbour} would not be sufficient.
The approximate version of \eqref{eq:distancegrows} is provided by a bound on the discrepancy within a neighborhood
of a vertex $v$, between two coupled Glauber chains which start in agreement around $v$.
This is stated in Proposition \ref{prop:dloc_uniformbound} in Section \ref{sec:inputs_dloc} below,
and provides a sort of ``homogeneity'' to the mixing, as already mentioned at the end of
Section \ref{sec:intro_iop_props}.
\subsection{Mixing inputs}
\label{sec:inputs_mixing}

In order to apply Corollary \ref{cor:barbour}, we need some technical inputs regarding the mixing properties
of Glauber dynamics.
The rapid mixing of the Glauber chain in the subcritical regime was proved by \cite{bhamidi2008mixing},
and there are two main steps.
First the density of the chain evolves with some drift towards the unique attracting fixed point $p$ of
$\phi_\beta$, and secondly the chain fully mixes among graphs with density close enough to $p$ in a
precise sense.
As we will mainly consider conditioned measures $\condmeas$, the first step of the argument is not so
relevant for our purposes.
However the second step will be quite important, and the fact that it holds for \emph{all} attracting fixed
points of $\phi_\beta$ allows it to be used beyond the subcritical regime.

To precisely state the lemma we need, let us specify what we mean by graphs with density close enough to $p$.
Rather than just having macroscopic subgraph counts be close to the corresponding features of the constant graphon
$W_p$, we need certain more \emph{local} subgraph counts to also have the correct behavior.
Specifically, for a graph $G = (\cV,\cE)$ with $\sv = |\cV|$ and $\se = |\cE| > 1$, following \cite{bhamidi2008mixing},
let us define
\begin{equation}
\label{eq:rdef}
    \sr_G(x,e) \coloneqq \left(
        \fr{\sN_G(x,e)}{2 \se n^{\sv-2}}
    \right)^{\fr{1}{\se-1}},
\end{equation}
where $\sN_G(x,e) = \partial_e \sN_G(x)$ is the number of homomorphisms of $G$ in $x^{+e}$ which use the
edge $e$, as mentioned in Section \ref{sec:intro_iop_hajek} (recall also the notation $\partial_e$ introduced
in \eqref{eq:partialdef}).
Notice that if $Y \sim \cG(n,p)$ then $\sr_G(Y,e) \pto p$ as $n \to \infty$.
Moreover, we may rewrite the ERGM Glauber dynamics transition probabilities \eqref{eq:glauberratio} as
\begin{align}
    \Pcond{\X_{t+1} = \X_t^{+e}}{e \text{ was chosen}} &= \phi \left(
        n^2 \sum_{j=0}^K \beta_j \fr{\sN_{G_j}(x,e)}{n^{\sv_j}}
    \right) \\
    &= \phi \left(2 \sum_{j=0}^K \beta_j \se_j \sr_{G_j}(\X_t, e)^{\se_j-1}\right),
\end{align}
where $\phi(z) = \fr{e^z}{1+e^z}$ and the graphs $G_j = (\cV_j, \cE_j)$ defining the ERGM Hamiltonian
satisfy $\sv_j = |\cV_j|$ and $\se_j = |\cE_j|$.
Further, for the $j=0$ term which corresponds to $G_0$ being a single edge, 
we interpret $\sr_{G_0}(\X_t,e)^{\se_j-1}$ as $1$, which is consistent with
the definition \eqref{eq:rdef}.
Thus, if every $\sr_{G_j}(\X_t, e)$ is close enough to $q$, then the probability that $\X_{t+1} = \X_t^{+e}$
if $e$ is resampled is approximately
\begin{equation}
\label{eq:phidef2}
    \phi_\beta(q) = \phi \left( 2 \sum_{j=0}^K \beta_j \se_j q^{\se_j-1} \right).
\end{equation}
Now, since
\begin{equation}
    \sH(q) = \sH(W_q) = \sum_{j=0}^K \beta_j q^{\se_j},
\end{equation}
we may easily see that \eqref{eq:phidef2} agrees with the previous definition \eqref{eq:phidef1} for $\phi_\beta$.
In any case, if $p$ is an attracting fixed point of $\phi_\beta$, then the chain will remain with density near
$p$ for a long period of time if it is started as such (in the sense that all $\sr_{G_j}(\X_0,e)$ are
approximately $p$).

Let us define the set of such ``good'' starting configurations, relative to the attracting fixed point
$p$ of $\phi_\beta$ and with some error $\eps > 0$, as follows:
\begin{equation}
\label{eq:gammadef}
    \gampe \coloneqq \left\{ x : \left| \sr_G(x,e) - p \right| < \eps \text{ for all } e \in \edgeset
    \text{ and } G = (\cV, \cE) \text{ with } |\cV| \leq \max_{0 \leq j \leq K} \sv_j \text{ and } |\cE| > 1 \right\}.
\end{equation}
This set already appeared in \cite{bhamidi2008mixing}, but we borrow the notation from
\cite{bresler2024metastable} instead.
Two features of $\gampe$ will be relevant for our present purposes.
First is the fact that the Glauber dynamics experiences \emph{contraction} under the monotone coupling
when started within $\gampe$.
This is the coupling wherein both chains choose the same edge $e$ to update at each step, and the
choices of whether to include or remove that edge, according to the probabilities \eqref{eq:glauberratio},
are coupled monotonically, for instance by using a background uniform random number between $0$ and $1$.
This ensures that any monotonicity relation between distinct chains will be preserved under the evolution.
The second relevant fact about $\gampe$ is that it is large under $\condmeas$, recalling this notation
for the ERGM measure conditioned on $\pball$ from \eqref{eq:condmeas}.
These facts were proved in \cite{bhamidi2008mixing} and \cite{bresler2024metastable} respectively,
and they are stated in Lemmas \ref{lem:contraction} and \ref{lem:cavitymethod} below.

In Lemma \ref{lem:contraction} we consider the standard Glauber dynamics with respect to the full ERGM measure,
which we denote by $(\tilde{X}_t)$; also recall from Sections \ref{sec:intro_iop_props} and \ref{sec:inputs_concentration}
the notation of $(\tilde{X}_t^x)$ for the chain started at $\tilde{X}_0^x = x$.
Also, in this lemma and the remainder of the article, $\dh(x,x')$ denotes the \emph{Hamming distance} between
graphs $x$ and $x'$, i.e.
\begin{equation}
\label{eq:hammingdef}
    \dh(x,x') \coloneqq \sum_{e \in \edgeset} |x(e) - x'(e)|,
\end{equation}
which is the number of edges $e \in \edgeset$ for which $x(e) \neq x'(e)$, recalling the
notation of edge indicators from Section \ref{sec:review_fluctuations_consequences}.
This lemma provides \emph{one-step contraction} with respect to the Hamming distance under the Glauber
dynamics started at adjacent graphs in $\gampe$, when the chains are monotonically coupled.

\begin{lemma}[Lemma 18 of \cite{bhamidi2008mixing}]
\label{lem:contraction}
For any $p \in U_\beta$, there are some $\eps, \kappa > 0$ such that,
whenever $x,x' \in \gampe$ have $\dh(x,x') = 1$, then
\begin{equation}
    \E\left[\dh(\tilde{X}_1^x, \tilde{X}_1^{x'})\right] \leq 1 - \fr{\kappa}{n^2}.
\end{equation}
\end{lemma}

Next, the following lemma provides the fact that $\gampe$ is large under $\condmeas$.

\begin{lemma}[Lemma 4.5 of \cite{bresler2024metastable}]
\label{lem:cavitymethod}
For any $p \in U_\beta$ and $\eps > 0$, there are some $\eta, \gamma > 0$ such that,
when $X \sim \condmeas$, we have
\begin{equation}
    \P[X \in \gampe] > 1 - e^{-\gamma n}.
\end{equation}
\end{lemma}

In the remainder of the article, we will consider Glauber dynamics with respect to $\condmeas$,
under which the chain will not make any moves that leave $\pball$, and we denote this
conditioned Glauber dynamics by $(X_t)$, or by $(X_t^x)$ when it starts at $x \in \pball$.
Conditioning the chain to lie in $\pball$ introduces some additional challenges as the intersection
$\gampe \cap \pball$ may not be \emph{connected} under the nearest-neighbor relation where graphs $x$
and $x'$ are adjacent if $\dh(x,x') = 1$.
This poses a problem for applying the path coupling technique of \cite{bubley1997path}
to the conditioned Glauber dynamics using
the contraction within $\gampe$ provided by Lemma \ref{lem:contraction}.
Moreover, if two chains are monotonically coupled and one lies above the other,
any time one of them attempts to \emph{leave} the metastable
well $\pball$, then the ordering between the two chains may be broken.

Both of these issues are remedied by the following proposition which exhibits a large \emph{connected}
set $\Lambda^*$ which has desirable properties, including the fact that a chain started near $\Lambda^*$
stays in both $\gampe$ and $\phalfball$.
Staying in $\gampe$ allows us to use contraction under the monotone coupling, and staying in
$\phalfball$ ensures that the monotonicity relations will not be destroyed, since a chain may not even attempt
to leave $\pball$ if it is currently in $\phalfball$, at least for large enough $n$.

\begin{proposition}[Proposition 5.9 of \cite{winstein2025concentration}]
\label{prop:goodset}
Let $p \in U_\beta$ (we do not require the full Assumption \eqref{assumption}).
For all small enough $\eps > 0$, there are some $\eta, \gamma > 0$ such that for all large enough $n$ there is a set
$\Lambda^* \sse \pball$ satisfying the following properties, where $X \sim \condmeas$ and $(X_t^x)$ is the
Glauber dynamics conditioned on $\pball$ and started at $x$.
\begin{enumerate}[label=(\Roman*)]
    \item
    \label{property:lambdalarge}
    The set $\Lambda^*$ is large, i.e.\ we have
    \begin{equation}
        \P[X \in \Lambda^*] \geq 1 - e^{-\gamma n}.
    \end{equation}
    \item
    \label{property:lambdaplusstaygood}
    The chain started near $\Lambda^*$ stays in the good region.
    Namely, whenever $\dh(x,\Lambda^*) \leq 1$,
    \begin{equation}
        \P\left[X_t^x \in \gampe \cap \phalfball \text{ for all } t < e^{\gamma n}\right] \geq 1 - e^{-\gamma n}.
    \end{equation}
    \item
    \label{property:lambdadiambound}
    The set $\Lambda^*$ is connected under the nearest-neighbor relation.
    Moreover, for any $x,x' \in \Lambda^*$, we have
    \begin{equation}
        \dls(x,x') \leq 2 n^2,
    \end{equation}
    where $\dls(x,x')$ denotes the distance within $\Lambda^*$, i.e.\ the length of the shortest nearest-neighbor
    path between $x,x'$ which stays inside $\Lambda^*$.
\end{enumerate}
\end{proposition}

Proposition \ref{prop:goodset} allows us to apply the path coupling technique of \cite{bubley1997path}
to obtain certain markers of rapid mixing for the Glauber dynamics with respect to $\condmeas$.
Lemma \ref{lem:mixing} below makes this precise, providing various mixing results which will be used in the
proofs to follow.
This lemma is a modification of \cite[Lemma 5.12]{winstein2025concentration}, replacing one of the
facts which appeared in the statement of that lemma by a more fundamental fact which appeared in the proof.

\begin{lemma}
\label{lem:mixing}
Let $p \in U_\beta$ (we do not require the full Assumption \eqref{assumption}).
Let $\Lambda^*$ and $\eta,\gamma > 0$ be, respectively, the set and constants given in
Proposition \ref{prop:goodset}.
Additionally, let $\kappa$ be the constant in Lemma \ref{lem:contraction}.
Suppose $x \in \Lambda^*$ and $x' \in \pball$ with $\dh(x,x') = 1$.
Under the monotone coupling of the $\condmeas$ Glauber dynamics, for all $t \leq e^{\gamma n}$,
\begin{equation}
    \E\left[\dh(X_t^x,X_t^{x'})\right] \leq \left(1-\fr{\kappa}{n^2}\right)^t + e^{-\Omega(n)}.
\end{equation}
In particular, for all $t \geq n^3$,
\begin{equation}
    \P\left[X_t^x \neq X_t^{x'}\right] \leq e^{-\Omega(n)}.
\end{equation}
Moreover, for any $z \in \Lambda^*$ and $t \geq n^3$, we have
\begin{equation}
    \dtv(X_t^z, X) \leq e^{-\Omega(n)}.
\end{equation}
\end{lemma}

The last part of Lemma \ref{lem:mixing} guarantees that condition \ref{cond:mixing} of
Corollary \ref{cor:barbour} holds with $\delta = e^{-\Omega(n)}$ for any $z \in \Lambda^*$, as long as
$T \geq n^3$.
In the sequel we will always take $T = n^3$.
\subsection{Basic concentration inequality}
\label{sec:inputs_basic}

In this section we apply the mixing inputs of Section \ref{sec:inputs_mixing} to derive a basic
concentration inequality for Lipschitz observables of the ERGM
using Corollary \ref{cor:barbour}.
This follows the same reasoning as part of the proof of Theorem \ref{thm:lipschitzconcentration},
and the result is simply a version of part of that theorem.
However, we include the following streamlined derivation as it
follows a similar idea to what we will eventually use for the proofs of
Propositions \ref{prop:global_hajek_intro} and \ref{prop:local_hajek_intro}, and so
this section can be thought of as a warm-up.
We will also make use of the following result in the proof of Theorem \ref{thm:global_clt_intro},
though we could have used Theorem \ref{thm:lipschitzconcentration} instead.
As usual, we assume that $X \sim \condmeas$ for some $p \in U_\beta$ and some $\eta > 0$ such that
Proposition \ref{prop:goodset} and Lemma \ref{lem:mixing} hold.

\begin{lemma}
\label{lem:basic_concentration}
Let $p \in U_\beta$ (we do not require the full Assumption \ref{assumption}).
There are constants $c_1, c_2, c_3, c_4 > 0$ such that the following holds.
Let $f$ be a function on $n$-vertex graphs with Lipschitz vector $\cL$,
and let $\lambda > 0$ such that
$c_1 n^2 e^{-c_2 n} \leq \fr{\lambda}{\| \cL \|_\infty} \leq c_3 n^{\fr{3}{2}}$.
Then
\begin{equation}
    \P\left[ |f(X) - \E[f(X)]| > \lambda \right] \leq 
    2 \Exp{- \fr{c_4 \lambda^2}{\| \cL \|_\infty^2 n^2} }.
\end{equation}
\end{lemma}

\begin{proof}[Proof of Lemma \ref{lem:basic_concentration}]
We will apply Corollary \ref{cor:barbour} with the choice $\Lambda = \Lambda^*$, and with $T = n^3$.
Thus we need to check the conditions \ref{cond:variance}, \ref{cond:smalljumps},
\ref{cond:setlarge}, \ref{cond:mixing}, and \ref{cond:fbound} and provide appropriate bounds
on the numbers $V, \Delta, \epsilon, \delta$, and $M$ appearing in those conditions.

First, notice that for any $n$-vertex graphs $x$ and $x'$, we have
\begin{equation}
    |f(x) - f(x')| \leq \| \cL \|_\infty \cdot \dh(x,x'),
\end{equation}
since we may change one edge at a time to get from $x$ to $x'$, at each step incurring a difference
of at most $\| \cL \|_\infty$.
Thus, by Lemma \ref{lem:mixing}, for all $t \leq e^{\gamma n}$ we have
\begin{equation}
\label{eq:mainboundbasic}
    \left| \E[f(X_t^x)] - \E[f(X_t^{x'})] \right| \leq \| \cL \|_\infty \cdot \left(
        \left(1 - \fr{\kappa}{n^2} \right)^t + e^{-\Omega(n)}
    \right),
\end{equation}
where $x \in \Lambda^*$ and $\dh(x,x') \leq 1$, and we recall that $(X_t^x)$ denotes the conditioned
Glauber dynamics starting at $x$.
Therefore, to validate condition \ref{cond:variance}, for $x \in \Lambda^*$ we obtain
\begin{align}
    \sum_{x' \in \pball} P(x,x') \sum_{t=0}^{T-1} \left(
        \E[f(X_t^x)] - \E[f(X_t^{x'})]
    \right)^2 
    &\leq \sum_{e \in \edgeset} \fr{1}{\binom{n}{2}} \sum_{t=0}^{n^3-1} \left(
        \E[f(X_t^x)] - \E[f(X_t^{x^{\oplus e}})]
    \right)^2 \\
    &\lesssim \| \cL \|_\infty^2 \sum_{t=0}^{n^3-1} \left( 1 - \fr{\kappa}{n^2} \right)^{2t}
        + \| \cL \|_\infty^2 e^{-\Omega(n)} \\
    &\lesssim \| \cL \|_\infty^2 n^2,
\end{align}
where we used \eqref{eq:mainboundbasic} at the second inequality.
Note that in the first sum, when we write $x' \in \pball$, it is understood that we mean the sum
over all $n$-vertex graphs whose graphon representations lie in $\pball$.
The first inequality above is also worth commenting on: it follows because at each step of
Glauber dynamics, an edge is chosen \emph{uniformly}, and it may or may not be changed;
recall our notation $x^{\oplus e}$ for the graph $x$ with the status of edge $e$ flipped,
i.e.\ $x^{\oplus e}(e) = 1-x(e)$ and $x^{\oplus e}(e') = x(e')$ for all $e' \neq e$.
All this aside, we obtain condition \ref{cond:variance} with $V \lesssim \| \cL \|_\infty^2 n^2$.

We may also use \eqref{eq:mainboundbasic} to verify condition \ref{cond:smalljumps}.
For $t \leq e^{\gamma n}$, \eqref{eq:mainboundbasic} immediately implies that
\begin{equation}
\label{eq:Deltabound}
    \left| \E[f(X_t^x)] - \E[f(X_t^{x'})] \right| \lesssim \| \cL \|_\infty.
\end{equation}
For larger values of $t$, we must apply the second part of Lemma \ref{lem:mixing},
which states that for $t \geq n^3$ we have
\begin{equation}
    \P \left[ X_t^x \neq X_t^{x'} \right] \leq e^{-\Omega(n)},
\end{equation}
which implies \eqref{eq:Deltabound} for $t \geq n^3$ as well.
In other words, we obtain condition \ref{cond:smalljumps} with $\Delta \lesssim \| \cL \|_\infty$.

Next we turn to condition \ref{cond:setlarge}.
Let us consider the \emph{stationary} chain $(X_t)$, i.e.\ the chain with $X_0 \sim \condmeas$.
By Property \ref{property:lambdalarge} of $\Lambda^*$ in Proposition \ref{prop:goodset}, for each $t$ we have
\begin{equation}
    \P[X_t \in \Lambda^*] \geq 1 - e^{-\gamma n}.
\end{equation}
Thus, by a union bound,
\begin{equation}
    \P[X_t \in \Lambda^* \text{ for } 0 \leq t < n^3 ] \geq 1 - e^{-\Omega(n)}
\end{equation}
as well.
In particular there is some $z \in \Lambda^*$ for which we also have
\begin{equation}
    \P[X_t^z \in \Lambda^* \text{ for } 0 \leq t < n^3] \geq 1 - e^{-\Omega(n)}.
\end{equation}
In other words, we obtain condition \ref{cond:setlarge} with $\epsilon = e^{-\Omega(n)}$
for this choice of $z$.
Since $z \in \Lambda^*$, we also immediately obtain condition \ref{cond:mixing} with
$\delta = e^{-\Omega(n)}$ by the last part of Lemma \ref{lem:mixing}.
Finally, since we may get from any $n$-vertex graph to any other by changing at most $O(n^2)$
edges, incurring a change of at most $\| \cL \|_\infty$ at each step, we obtain the
boundedness condition \ref{cond:fbound} with $M \lesssim \| \cL \|_\infty n^2$.
All together, we have
\begin{equation}
    V \leq C \| \cL \|_\infty^2 n^2, \qquad
    \Delta \leq C \| \cL \|_\infty, \qquad
    \epsilon \leq e^{-c n}, \qquad
    \delta \leq e^{-c n}, \qquad \text{and} \qquad
    M \leq C \| \cL \|_\infty n^2,
\end{equation}
for some $C c > 0$.
Therefore, the conclusion of Corollary \ref{cor:barbour} is
\begin{equation}
    \P \left[
        |f(X) - \E[f(X)]| > \lambda
    \right] \leq 2 \Exp{
        - \fr{\left( \lambda - 2 C \| \cL \|_\infty n^2 e^{- c n} \right)^2}
        {2 C \| \cL \|_\infty^2 n^2 + \fr{4}{3} C \| \cL \|_\infty \left( \lambda - 2 C \| \cL \|_\infty n^2 e^{-c n} \right)}
    } + 2e^{-cn}.
\end{equation}
Now if $\lambda > 3 C \| \cL \|_\infty n^2 e^{-cn}$, then
$\lambda - 2 C \| \cL \|_\infty n^2 e^{-cn} > \fr{1}{3} \lambda$, and since
$\fr{z^2}{a + b z}$ is an increasing function of $z \geq 0$ whenever $a,b \geq 0$, we obtain
\begin{equation}
\label{eq:sumofexpos}
    \P \left[
        |f(X) - \E[f(X)]| > \lambda
    \right] \leq 2 \Exp{
        - \fr{\fr{1}{9} \lambda^2}
        {2 C \| \cL \|_\infty^2 n^2 + \fr{4}{9} C \| \cL \|_\infty \lambda}
    } + 2e^{-cn}.
\end{equation}
Next, if $\lambda \lesssim \| \cL \|_\infty n^{\fr{3}{2}}$, then 
\begin{equation}
    \fr{\fr{1}{9} \lambda^2}{2 C \| \cL \|_\infty^2 n^2 + \fr{4}{9} C \| \cL \|_\infty \lambda}
    \lesssim n,
\end{equation}
meaning that the first exponential in \eqref{eq:sumofexpos} is the larger term.
Finally, for such $\lambda$, the first term in the denominator inside the first exponential dominates
the second term, leading to the desired conclusion.
\end{proof}

Let us compare this result with Theorem \ref{thm:lipschitzconcentration} above.
As discussed below the statement of that theorem, the result we obtained here with the most simple
application of Corollary \ref{cor:barbour} only obtains a variance proxy of order $\| \cL \|_\infty^2 n^2$,
whereas Theorem \ref{thm:lipschitzconcentration} obtains $\| \cL \|_1  \| \cL \|_\infty$, which is
smaller in general (note that these Lipschitz vectors have $\binom{n}{2}$ entries), especially for
\emph{local} observables which don't depend strongly on too many edges.
However, Theorem \ref{thm:lipschitzconcentration} also has an exponential term \emph{inside} the main
exponential, which destroys the bound when considering \emph{global} observables such as the overall
edge count.
Luckily, for such global observables, the quantity $\| \cL \|_\infty^2 n^2$ is often close to the
true variance, meaning that by using both Theorem \ref{thm:lipschitzconcentration} and Lemma
\ref{lem:basic_concentration}, we can adequately cover many local \emph{and} global Lipschitz observables.

Despite this, in order to prove our main results, in particular Propositions \ref{prop:global_hajek_intro}
and \ref{prop:local_hajek_intro}, even the two of these bounds combined will not suffice.
Instead, we need to embark upon a more careful analysis which goes beyond the simple calculation of
the Lipschitz vectors of the relevant quantities.
Nevertheless, Theorem \ref{thm:lipschitzconcentration} and Lemma \ref{lem:basic_concentration}
will be useful along the way.
\subsection{Homogeneity of mixing}
\label{sec:inputs_dloc}

Lemma \ref{lem:mixing} provides exponential contraction of the overall Hamming distance under Glauber dynamics with
the monotone coupling between different chains,
which is a crucial input for applying Corollary \ref{cor:barbour} to prove both Propositions \ref{prop:global_hajek_intro}
and \ref{prop:local_hajek_intro}, leading to Theorems \ref{thm:global_clt_intro} and \ref{thm:local_clt_intro}.
This strategy was hinted at just above in Section \ref{sec:inputs_basic}.
However, as mentioned briefly at the end of Section \ref{sec:intro_iop} as well as the end of
Section \ref{sec:inputs_concentration}, in order to prove the \emph{local} results
we will also require control on the \emph{local} Hamming distance.
More precisely, recall from the definition \eqref{eq:dlocdef} that
\begin{equation}
    \dloc(x,x') = \sum_{e \ni v} |x(e) - x'(e)|,
\end{equation}
which is the number of edges adjacent to $v$ at which $x$ and $x'$ differ.

In this section, we prove Proposition \ref{prop:dloc_uniformbound_intro}, restated below for the reader's
convenience.
This proposition gives a bound on
$\dloc(X_t^x, X_t^{x'})$ under the monotone coupling, in the case where $\dloc(x,x') = 0$ but $x$ and $x'$
may differ at some edge away from $v$.
This leads to an approximate version of \eqref{eq:distancegrows}, showing that there can't be a neighborhood
of a vertex where the chain becomes ``unmixed'' locally while continuing to mix globally.
In other words, the mixing must happen somewhat homogeneously.

To state the proposition, recall from Section \ref{sec:inputs_mixing} that $(X_t^x)$ denotes the conditioned
Glauber dynamics, which has $\condmeas$ as its stationary distribution, and which starts at $X_0^x = x \in \pball$,
recalling the notation for $\pball$ and $\condmeas$ from \eqref{eq:balldef} and \eqref{eq:condmeas}.
The lemma statement should be interpreted as holding for some $\eta > 0$.
Recall also the definition of the monotone coupling from Section \ref{sec:inputs_mixing}, and our convention
that all chains are monotonically coupled.
Finally, recall the good connected set $\Lambda^*$ introduced in Proposition \ref{prop:goodset}.

\begin{proposition}[Restatement of Proposition \ref{prop:dloc_uniformbound_intro}]
\label{prop:dloc_uniformbound}
Let $p \in U_\beta$ (we do not require the full Assumption \ref{assumption}).
There is some constant $C > 0$ such that the following holds.
Let $v \in [n]$ and $e \in \edgeset$ with $v \notin e$.
Suppose that $x \in \Lambda^*$, and set $x' = x^{\oplus e}$.
Then for all $t \geq 0$, we have
\begin{equation}
    \E\left[\dloc(X_t^x,X_t^{x'})\right] \leq \fr{C}{n}.
\end{equation}
\end{proposition}

To prove this bound, we first prove a recurrence inequality in Lemma \ref{lem:dloc_recurrence} below, and then
solve this recurrence to obtain the bound in Proposition \ref{prop:dloc_uniformbound_intro}/\ref{prop:dloc_uniformbound}.
This lemma statement includes the sets $\gampe$ and $\phalfball$, and should be interpreted as the
existence of $\eta > 0$ for any $\eps > 0$, such that the statement holds, the same as in Lemma \ref{lem:cavitymethod}.

\begin{lemma}
\label{lem:dloc_recurrence}
Let $p \in U_\beta$ (we do not require the full Assumption \ref{assumption}).
Let $v \in [n]$, and 
suppose that $x, x' \in \gampe \cap \phalfball$ with $x \preceq x'$.
Then under the monotone coupling,
\begin{equation}
    \E \left[\dloc(X_1^x, X_1^{x'})\right]
    \leq \left(1 - \fr{\kappa}{n^2}\right) \dloc(x,x')
    + \fr{C}{n^3} \dh(x,x').
\end{equation}
\end{lemma}

\begin{proof}[Proof of Lemma \ref{lem:dloc_recurrence}]
First write $x = x_0 \preceq x_1 \preceq \dotsb \preceq x_k = x'$, where $k = \dh(x,x')$, such that we have
$\dh(x_{i-1},x_i) = 1$ for each $i = 1, \dotsc, k$.
Note that $\gampe$ and $\phalfball$ are \emph{intervals} in the Hamming cube, so each $x_i \in \gampe \cap \phalfball$.
The interval property for $\gampe$ holds since it is an intersection of pullbacks of intervals in $\R$ under
monotone functions $\sr_G(-,e)$, and for $\phalfball$ the interval property is provided by
\cite[Lemma 7.3]{bresler2024metastable}, for instance.

In any case, we have
\begin{equation}
\label{eq:dloc_sum}
    \E\left[\dloc(X_1^x,X_1^{x'})\right] \leq \sum_{i=1}^k \E\left[ \dloc(X_1^{x_{i-1}}, X_1^{x_i}) \right],
\end{equation}
since $\dloc$ satisfies the triangle inequality (although it is not a metric).
Now, for exactly $\dloc(x,x')$ indices $i$, the edge at which $x_{i-1}$ and $x_i$ differ is adjacent to $v$.
For these indices $i$, we have
\begin{equation}
    \E\left[\dloc(X_1^{x_{i-1}}, X_1^{x_i})\right] \leq \E\left[\dh(X_1^{x_{i-1}}, X_1^{x_i})\right]
    \leq 1 - \fr{\kappa}{n^2}
\end{equation}
by Lemma \ref{lem:contraction}, since all $x_i$ are in $\phalfball$ so the \emph{conditioned} Glauber steps
are the same as the \emph{unconditioned} Glauber steps.

For all the (at most $\dh(x,x')$) other indices $i$, the edge at which $x_{i-1}$ and $x_i$ differ is \emph{not}
adjacent to $v$.
Let us fix such an index $i$ and call the corresponding edge $e_i$, so that $v \notin e_i$.
Then for any $e$, by the definition of the monotone coupling and recalling the transition probabilities
\eqref{eq:glauberratio}, we have
\begin{align}
    \P\left[X_1^{x_{i-1}}(e) \neq X_1^{x_i}(e)\right] &= \fr{1}{\binom{n}{2}} \left|
        \phi(n^2 \partial_e \sH(x_{i-1})) - \phi(n^2 \partial_e \sH(x_i)) \right| \\
        &\lesssim \fr{1}{n^2} \left| n^2 \partial_e \sH(x_{i-1}) - n^2 \partial_e \sH(x_i)\right| \\
        &\leq \fr{1}{n^2} \sum_{j=1}^K \fr{\beta_j}{n^{\sv_j-2}} \partial_e \sN_{G_j}(x_i,e_i),
\end{align}
recalling also the notation $\partial_e$ defined in \eqref{eq:partialdef}, and using the fact that
the function $\phi(z) = \fr{e^z}{1+e^z}$ is $1$-Lipschitz.

Now for $e_i \neq e$ the differential $\partial_e \sN_{G_j}(x_i,e_i)$ is equal to $\sN_{G_j}(x,e_i,e)$,
which is the number of homomorphisms of $G_j$ in $x^{+e_i+e}$ using both $e_i$ and $e$.
Recalling the definition \eqref{eq:dlocdef} of $\dloc$, we need to bound this for $e$ such that $v \in e$.
There are two cases: either $e$ and $e_i$ do not share a vertex, or they do.
There are exactly \emph{two} choices of $e$ with $v \in e$ and $e \cap e_i \neq \emptyset$, and for these edges we have
\begin{equation}
    \sN_{G_j}(x_i,e_i,e) \lesssim n^{\sv_j-3}
\end{equation}
since $3$ vertices are used up by $e_i$ and $e$ together.
For all of the other (at most $n$) edges $e$ adjacent to $v$ with $e \cap e_i = \emptyset$,
since $e_i$ and $e$ together use up $4$ vertices, we have
\begin{equation}
    \sN_{G_j}(x_i,e_i,e) \lesssim n^{\sv_j-4}.
\end{equation}
Putting this all together, there are two edges $e$ with $v \in e$ for which
\begin{equation}
    \P\left[X_1^{x_{i-1}}(e) \neq X_1^{x_i}(e)\right] \lesssim n^{-3},
\end{equation}
and there are at most $n$ edges $e$ with $v \in e$ for which
\begin{equation}
    \P\left[X_1^{x_{i-1}}(e) \neq X_1^{x_i}(e)\right] \lesssim n^{-4}.
\end{equation}
This yields
\begin{equation}
    \E\left[\dloc(X_1^{x_{i-1}}, X_1^{x_i})\right] = \sum_{e \ni v} \P[X_1^{x_i-1}(e) \neq X_1^{x_i}(e)] \leq C n^{-3}
\end{equation}
for the indices $i$ where the edge $e_i$ at which $x_{i-1}$ and $x_i$ differ is not adjacent to $v$, for some
constant $C$.

So we have bounded each term in \eqref{eq:dloc_sum}, and we obtain
\begin{equation}
    \E\left[\dloc(X_1^x, X_1^{x'})\right] \leq \left(1 - \fr{\kappa}{n^2}\right) \dloc(x,x')
    + \fr{C}{n^3} \dh(x,x')
\end{equation}
as required.
\end{proof}

With the recurrence inequality Lemma \ref{lem:dloc_recurrence} in hand, we may finish the proof of the
desired bound by considering the solution to a general recurrence of this form.

\begin{proof}[Proof of Proposition \ref{prop:dloc_uniformbound_intro}/\ref{prop:dloc_uniformbound}]
Let $\Lambda^*$ and $\gamma$ be the set and the constant in Proposition \ref{prop:goodset},
so that
\begin{equation}
    \P\left[X_s^x \in \gampe \cap \phalfball \text{ for all } s < e^{\gamma n}\right] \geq 1 - e^{-\gamma n}
\end{equation}
whenever $\dh(x,\Lambda^*) \leq 1$; this is property \ref{property:lambdaplusstaygood}
of that proposition.
Let us denote by $\cA_t$ the event that this holds for both $x$ and $x'$ simultaneously under the monotone coupling
up to time $t$, for any fixed $t \leq e^{\gamma n}$.
More precisely, set
\begin{equation}
    \cA_t = \left\{ X_s^x, X_s^{x'} \in \gampe \cap \phalfball \text{ for all } s < t \right\}
\end{equation}
and note that $\P[\cA_t] \geq 1 - 2 e^{-\gamma n}$ by a union bound, recalling that $x \in \Lambda^*$
and that $x' = x^{\oplus e}$ for some $e \not \ni v$, so $\dh(x',\Lambda^*) \leq 1$ as well.

On the event $\cA_t$, under the monotone coupling, since neither chain has had a chance to attempt to escape
$\pball$, the initial ordering (either $x \preceq x'$ or vice versa) is retained throughout the evolution
until time $t-1$.
Therefore we can apply Lemma \ref{lem:dloc_recurrence} to the last time step as follows:
\begin{align}
    \E\left[\dloc(X_t^x, X_t^{x'})\right]
    &= \Econd{\dloc(X_t^x, X_t^{x'})}{\cA_t} \cdot \P[\cA_t] + \Econd{\dloc(X_t^x, X_t^{x'})}{\cA_t^\c} \cdot \P[\cA_t^\c] \\
    &\leq \Econd{\left(1 - \fr{\kappa}{n^2}\right) \dloc(X_{t-1}^x, X_{t-1}^{x'})
        + \fr{C}{n^3} \dh(X_{t-1}^x,X_{t-1}^{x'})}{\cA_t} \cdot \P[\cA_t]
        + e^{-\Omega(n)} \\
    &\leq \left(1 - \fr{\kappa}{n^2} \right) \E\left[\dloc(X_{t-1}^x,X_{t-1}^{x'})\right]
        + \fr{C}{n^3} \E\left[\dh(X_{t-1}^x,X_{t-1}^{x'})\right]
        + e^{-\Omega(n)}.
\end{align}
In first inequality above, we applied Lemma \ref{lem:dloc_recurrence} via the tower property by conditioning on
$X_{t-1}^x, X_{t-1}^{x'}$ and using the fact that these are in $\gampe \cap \phalfball$ and are
ordered under $\cA_t$.
In the second inequality, we simply used the fact that $\E[A \ind*{\cA_t}] \leq \E[A]$ for any random variable $A$.

Now, by Lemma \ref{lem:mixing}, $\E\left[\dh(X_{t-1}^x, X_{t-1}^{x'})\right]$ is bounded above by some constant for
all $t \leq e^{\gamma n}$, and so by adjusting the constant $C$ to absorb this and the extra $e^{-\Omega(n)}$
term, we obtain
\begin{equation}
    \E\left[\dloc(X_t^x, X_t^{x'})\right] \leq \left(1 - \fr{\kappa}{n^2}\right) \E\left[\dloc(X_{t-1}^x, X_{t-1}^{x'})\right]
    + \fr{C}{n^3}.
\end{equation}
This means that, for all $t \leq e^{\gamma n}$, we have
\begin{equation}
    \E\left[\dloc(X_t^x, X_t^{x'})\right] \leq a_t
\end{equation}
where $(a_t)_{t=0}^\infty$ is a sequence satisfying $a_0 = 0$ and
\begin{equation}
    a_t = \left(1 - \fr{\kappa}{n^2} \right) a_{t-1} + \fr{C}{n^3}.
\end{equation}
In general, the solution to the initial-value recurrence problem
\begin{equation}
    \begin{cases}
        a_t = b \cdot a_{t-1} + c \\
        a_0 = 0
    \end{cases}
\end{equation}
with $b \neq 1$ is given uniquely by
\begin{equation}
\label{eq:recurrencesolution}
    a_t = c \fr{1 - b^t}{1 - b},
\end{equation}
since $a_t$ is simply the sum of $t$ terms in a geometric sequence.
Another way to see this is that \eqref{eq:recurrencesolution}
satisfies both equations of the initial-value recurrence problem
and the solution of such a problem is unique.
Plugging in $b = 1 - \fr{\kappa}{n^2}$ and $c = \fr{C}{n^3}$, we obtain
\begin{align}
    a_t &= \fr{C}{n^3} \fr{n^2}{\kappa} \left(1 - \left(1 - \fr{\kappa}{n^2}\right)^t \right) \\
    &\leq \fr{C}{\kappa n}
\end{align}
as desired, letting the constant of the proposition statement be $\fr{C}{\kappa}$.

This finishes the proof in the case where $t \leq e^{\gamma n}$, and the case where
$t > e^{\gamma n}$ follows immediately from the total-variation distance bound in
Lemma \ref{lem:mixing}, using the fact that $\dloc(a,b) \leq n$ for all $a,b$.
\end{proof}

We remark that a stronger bound which depends on $t$ may be derived using the same
strategy, and this may prove useful in other situations.
However, for our present purposes, the stated bound will suffice.
\section{H\'ajek projections}
\label{sec:hajek}

In this section we come to two of the main contributions of the present article, Propositions \ref{prop:global_hajek_intro}
and \ref{prop:local_hajek_intro}.
We remind the reader that Proposition \ref{prop:global_hajek_intro} states quantitatively that the fluctuations of
$\sN_G(X)$ are driven by those of $\sE(X)$, and Proposition \ref{prop:local_hajek_intro} is the analogous statement
for local homomorphism counts (i.e.\ the number of homomorphisms using a particular vertex).
We prove Proposition \ref{prop:global_hajek_intro} in Section \ref{sec:hajek_global}, and
Proposition \ref{prop:local_hajek_intro} in Section \ref{sec:hajek_local} below.
Finally, in Section \ref{sec:hajek_subgraphs} we prove Corollaries \ref{cor:global_subgraph_clt_intro}
and \ref{cor:local_subgraph_clt_intro} (the CLTs for global and local subgraph counts $\sN_G(X)$
and $\sN_G^v(X)$ respectively),
assuming the results of Theorem \ref{thm:global_clt_intro} (the CLT for edge count $\sE(X)$) and
Theorem \ref{thm:local_clt_intro} (the CLT for vertex degree $\deg_v(X)$) respectively.

Throughout this section, we fix a finite simple graph $G = (\cV, \cE)$ and set $\sv = |\cV|$ and $\se = |\cE|$.
For $\rho \in \cV$, we also denote by $\sd_\rho$ the degree of $\rho$ in $G$.
In what follows, we always make one of the following two assumptions: either we have $M_\beta = U_\beta = \{p\}$
(recalling these notations from the beginning of Section \ref{sec:intro_results}, for instance),
i.e.\ we are in a non-critical phase uniqueness case, or otherwise we fix some $p \in U_\beta$ and assume that
the graph $G$, as well as the specifying graphs $G_0, G_1, \dotsc, G_K$ of the ERGM,
are all \emph{forests}, meaning that they contain no cycles (this is Assumption \ref{assumption}).
This assumption is necessary only because we use Proposition \ref{prop:Eproduct}
(including implicit use via Proposition \ref{prop:marginal}).

The assumption that the graphs $G, G_0, G_1, \dotsc, G_K$ are forests in the phase coexistence case
can be dropped in exchange for slightly worse conclusions;
this will be commented on at the end of each subsection, and will be relevant for proving a CLT 
for $\sN_G(X)$ in the phase coexistence case when $G$ is not a forest
(but the specifying graphs $G_0, G_1, \dotsc, G_K$ of the ERGM are).
Note that we \emph{do not} make any other assumptions on $G$, such as being balanced
\cite[Page 50, Definition 1]{alon2016probabilistic} or even connected.

We will always denote by $X$ a sample from $\condmeas$, recalling this notation from \eqref{eq:condmeas}.
Every result should be interpreted as holding for some fixed $\eta > 0$, which is determined by some arbitrary
(but small enough) $\eps > 0$ so that Lemma \ref{lem:cavitymethod}, which states that $\gampe$ is
large under $\condmeas$, holds.
Furthermore, we let $(X_t^x)$ denote the Glauber dynamics with respect
to $\condmeas$ (i.e.\ the dynamics which does not step outside of $\pball$), starting at $X_0^x = x$,
and we use $P$ to denote the transition kernel of this Markov chain.

\subsection{Global subgraph counts}
\label{sec:hajek_global}

In this section we prove Proposition \ref{prop:global_hajek_intro}, which has been restated
below for the reader's convenience.
Keep in mind the ``phase-uniqueness or forest'' assumption \ref{assumption} which we make throughout
all of Section \ref{sec:hajek}, except for in brief discussions at the end of each subsection.

\begin{proposition}[Restatement of Proposition \ref{prop:global_hajek_intro}]
\label{prop:global_hajek}
Let us make Assumption \ref{assumption}.
As in \eqref{eq:hatngdef}, let us write
\begin{equation}
    \hsN_G(x) = \sN_G(x) - 2 \se p^{\se-1} n^{\sv-2} \sE(x).
\end{equation}
Then, for all $\zeta \in (0,1)$, there is some $c > 0$ such that
\begin{equation}
    \P\left[\left|\hsN_G(X) - \E\left[\hsN_G(X)\right]\right| > n^{\sv-1.5+\zeta}\right] \leq \Exp{-cn^\zeta}.
\end{equation}
In particular, for all $\xi > 0$ we have
\begin{equation}
    \Var\left[\hsN_G(X)\right] \lesssim n^{2\sv-3+\xi}.
\end{equation}
\end{proposition}

As described in Section \ref{sec:intro_iop_concentration}, the broad strategy for proving
Proposition \ref{prop:global_hajek_intro}/\ref{prop:global_hajek} is to first bound the fluctuations of the
\emph{changes} in $\sN_G(X_t)$ under Glauber dynamics, and then use this as
an input for the concentration inequality Corollary \ref{cor:barbour}.
If an edge $e$ is flipped, taking $x^{-e}$ to $x^{+e}$ or vice versa, the change in $\sN_G(x)$ is
$\partial_e \sN_G(x)$, recalling the notation introduced in \eqref{eq:partialdef}.
This is the same as $\sN_G(x,e)$, which is the number of homomorphisms of $G$ in $x^{+e}$ which use the edge $e$,
as introduced in Section \ref{sec:intro_iop}.
The following lemma provides the concentration of $\sN_G(X,e)$ which is required as input for the broader strategy.

\begin{lemma}
\label{lem:global_step_concentration}
Let us make Assumption \ref{assumption}.
For $\alpha \in (0,\fr{1}{2})$ there is some $c > 0$ such that
\begin{equation}
    \P \left[ | \sN_G(X,e) - 2 \se p^{\se-1} n^{\sv-2} | > n^{\sv - 2.5 + \alpha} \right]
    \leq \Exp{- c n^{2\alpha}}.
\end{equation}
\end{lemma}

\begin{proof}[Proof of Lemma \ref{lem:global_step_concentration}]
We will first invoke Theorem \ref{thm:lipschitzconcentration} to obtain concentration of
$\sN_G(X,e)$ around \emph{its mean}, and then estimate the mean.
To apply Theorem \ref{thm:lipschitzconcentration}, we need to construct a Lipschitz
vector $\L$ for the function $x \mapsto \sN_G(x,e)$.
First, we can take $\L_e = 0$ since $\sN_G(x,e)$ does not depend on $x(e)$, being the count of copies of
$G$ in $x^{+e}$ (the edge $e$ is always added).
Next, for any $e' \in \binom{[n]}{2}$ with $e' \neq e$, we have
\begin{equation}
    \sN_G(x^{+e'},e) - \sN_G(x^{-e'},e) = \partial_{e'} \sN_G(x,e) = \sN_G(x,e,e'),
\end{equation}
which is the number of labeled copies of $G$ in $x^{+e+e'}$ using both edges $e$ and $e'$.
If $e$ and $e'$ share a vertex then there are $\lesssim n^{\sv-3}$ of these, and if $e$ and $e'$ do not
share a vertex then there are $\lesssim n^{\sv-4}$.
In other words we can take $\L_{e'} \lesssim n^{\sv-3}$ for at most $n$ edges $e'$, and we can take
$\L_{e'} \lesssim n^{\sv-4}$ for the remaining at most $n^2$ edges.
This yields
\begin{equation}
    \| \L \|_\infty \lesssim n^{\sv-3} \qquad \text{and} \qquad
    \| \L \|_1 \lesssim n^{\sv - 3} \cdot n + n^{\sv - 4} \cdot n^2 \lesssim n^{\sv-2}.
\end{equation}
Now we apply Theorem \ref{thm:lipschitzconcentration}, slightly adjusting the constants to account
for the $\lesssim$'s above.
Plugging in $\lambda = n^{\sv-2.5+\alpha}$, which is less than $c_1 n^{\sv-2}$ for $n$ large enough,
we obtain
\begin{align}
    \P\left[|\sN_G(X,e) - \E[\sN_G(X,e)]| > n^{\sv-2.5+\alpha} \right]
    &\leq 2 \Exp{- c_2 \fr{n^{2\sv-5+2\alpha}}{n^{\sv-2}n^{\sv-3}} + \Exp{c_3 \fr{n^{\sv-2.5+\alpha}}{n^{\sv-3}} - c_4 n} } + e^{-c_5 n}. \\
    &= 2 \Exp{- c_2 n^{2\alpha} + \Exp{c_3 n^{0.5+\alpha} - c_4 n}} + e^{-c_5 n}.
\end{align}
Since $\alpha < \fr{1}{2}$, the inner exponential term tends to $0$, and also
$e^{-\Omega(n^{2\alpha})} \gg e^{-\Omega(n)}$, so we can combine terms and change the constants to yield the existence
of some $c > 0$ for which
\begin{equation}
    \P \left[ | \sN_G(X,e) - \E[\sN_G(X,e)] | > n^{\sv - 2.5 + \alpha} \right]
    \leq \Exp{- c n^{2\alpha}},
\end{equation}
which is the desired rate of concentration.

It just remains to estimate $\E[\sN_G(X,e)]$; for this, we invoke Propositions \ref{prop:Eproduct} and \ref{prop:marginal},
which together imply
\begin{equation}
\label{eq:expectedproduct}
    \E[\sN_G(X,e)] = 2 \se p^{\se-1} n^{\sv-2} + O(n^{\sv-2.5} \sqrt{\log n}).
\end{equation}
To see why \eqref{eq:expectedproduct} holds, note that all but $O(n^{\sv-3})$ of the homomorphisms counted
in $\sN_G(X,e)$ are injective, simply because there are only $O(n^{\sv-3})$ non-injective maps $\cV \to [n]$
which contain $e$ in their image.
To define an injective map that uses $e$, one must choose which edge $\phi \in \cE$ to map to
$e \in \edgeset$, and the orientation, giving $2 \se$ choices.
The other vertex variables are free, giving $(n-2)(n-3) \dotsb (n-\sv+1) = n^{\sv-2} + O(n^{\sv-3})$ 
choices of overall injective map.
For each such map $m$ which sends $\phi$ to $e$, the probability that it is a homomorphism in $X^{+e}$ is
\begin{equation}
\label{eq:eprod}
    \E \left[ \prod_{\phi' = \{\rho_1, \rho_2\} \in \cE \setminus \{\phi\}}
    X(\{m(\rho_1), m(\rho_2)\}) \right].
\end{equation}
Now, if we are in the phase uniqueness case then Proposition \ref{prop:Eproduct} immediately applies,
and if we are in the case where $G$ is a forest, then since the map $m$ is injective, the images of 
edges under $m$ also form a forest in the complete graph on $[n]$.
So, in both cases we may apply Proposition \ref{prop:Eproduct} to conclude that \eqref{eq:eprod}
equals
\begin{equation}
\label{eq:prode}
    \E[X(e)]^{\se-1} + O \left( \fr{1}{n} \right).
\end{equation}
Further, since we assume that either we are in the phase uniqueness case or that all specifying
graphs $G_0, G_1, \dotsc, G_K$ are forests, we may also apply Proposition \ref{prop:marginal} to conclude
that \eqref{eq:prode} equals
\begin{equation}
    p^{\se-1} + O \left( n^{-\fr{1}{2}} \sqrt{\log n} \right).
\end{equation}
Combined with the above reasoning, this proves \eqref{eq:expectedproduct}.
The conclusion of the lemma thus follows from the fact that $n^{\sv-2.5+\alpha} \gg n^{\sv-2.5}\sqrt{\log n}$.
\end{proof}

We remark that the verification of \eqref{eq:expectedproduct} is the \emph{only} place in the proof of
Proposition \ref{prop:global_hajek_intro}/\ref{prop:global_hajek} where the ``phase-uniqueness or forest''
assumption is used.
Without this assumption, a similar proposition could indeed be proved, simply by replacing the factor of
$2 \se p^{\se-1} n^{\sv-2}$ in the statement with $\E[\sN_G(X,e)]$.
However, we will need the explicit form of this H\'ajek projection factor later, in the proof of
Theorem \ref{thm:global_clt_intro}, as certain quantities will need to multiply and cancel nicely.
Moreover, this form is what leads to the explicit representation of the proxy $\sigma_n^2$ given in 
\eqref{eq:sigmadef_global_intro} for the variance of $\sE(X)$.

Now, as previously stated, we will use Corollary \ref{cor:barbour} to prove
Proposition \ref{prop:global_hajek_intro}/\ref{prop:global_hajek}.
This corollary requires the selection of a ``nice set'' $\Lambda$ where the dynamics remains with high probability,
and where the behavior of the chain can be controlled.
We would like $\Lambda$ to be a set where every variable $\sN_G(X,e)$ is close to $2 \se p^{\se-1} n^{\sv-2}$
for all $e \in \edgeset$.
By a union bound over the polynomially-many edges $e$, Lemma \ref{lem:global_step_concentration} shows that this
condition is very \emph{likely}, but this alone is not a priori enough to use as input for Corollary \ref{cor:barbour}.
The following lemma allows us to ``upgrade'' sets $\Pi$ which are simply \emph{large} to sets $\Lambda$ which
are \emph{nice enough} to be used for our purposes, and it will be helpful for both this section and the next.
The proof is similar to part of the proof of Proposition \ref{prop:goodset} (which was originally formulated as
\cite[Proposition 5.9]{winstein2025concentration}).

\begin{lemma}
\label{lem:goodset_from_bigset}
Let $p \in U_\beta$ (we do not require the full Assumption \ref{assumption}).
Suppose that $\Pi \sse \pball$ with $\P[X \in \Pi] \geq 1 - \Exp{-\Omega(n^\xi)}$ for some $\xi \in (0,1)$.
Then there is a set $\Lambda \sse \Lambda^*$ (where $\Lambda^*$ is the set given in
Proposition \ref{prop:goodset}) such that
\begin{equation}
    \P[X_t^x \in \Pi \text{ for all } 0 \leq t \leq n^3] \geq 1 - \Exp{-\Omega(n^\xi)}
\end{equation}
for all $x$ with $\dh(x,\Lambda) \leq 1$.
Moreover, there is some $z \in \Lambda$ for which
\begin{equation}
    \P[X_t^z \in \Lambda \text{ for all } 0 \leq t \leq n^3] \geq 1 - \Exp{-\Omega(n^\xi)}.
\end{equation}
\end{lemma}

\begin{proof}[Proof of Lemma \ref{lem:goodset_from_bigset}]
Consider the stationary chain $(X_t)$.
By a union bound over $0 \leq t \leq n^3$, we have
\begin{equation}
    \P[X_t \in \Pi \text{ for } 0 \leq t \leq n^3] \geq 1 - \Exp{-\Omega(n^\xi)}.
\end{equation}
So by the tower property of conditional expectation (conditioning on $X_0$ in the stationary chain)
this implies that, for some $c > 0$, the set
\begin{equation}
    \tilde{\Lambda} \coloneqq \left\{
        x \in \pball : \P\left[X_t^x \in \Pi \text{ for } 0 \leq t \leq n^3\right]
        \geq 1 - \Exp{- c n^\xi} \right\}
\end{equation}
also satisfies
\begin{equation}
    \P\left[X \in \tilde{\Lambda} \right] \geq 1 - \Exp{-\Omega(n^\xi)}.
\end{equation}
Now set
\begin{equation}
    \Lambda^+ = \Lambda^* \cap \tilde{\Lambda},
\end{equation}
where $\Lambda^*$ is as in Proposition \ref{prop:goodset}.
Since $\P[X \in \Lambda^*] \geq 1 - \Exp{-\Omega(n)}$ and $\xi < 1$, we also have
\begin{equation}
\label{eq:lambdaminus}
    \P[X \in \Lambda^+] \geq 1 - \Exp{-\Omega(n^\xi)},
\end{equation}
by a union bound.
Finally, we set
\begin{equation}
    \Lambda = \left\{
        x \in \Lambda^+ : x' \in \Lambda^+ \text{ for all } x' \in \pball
        \text{ with } \dh(x,x') = 1
    \right\}.
\end{equation}
By construction, this set $\Lambda$ satisfies the first property in the lemma statement.
Finally we claim that
\begin{equation}
\label{eq:toprove}
    \P[X \in \Lambda] \geq 1 - \Exp{-\Omega(n^\xi)}.
\end{equation}
Given this for now, applying the strategy just employed but with $\Pi$ replaced by $\Lambda$,
we find that the set of starting
states $z$ satisfying the second statement in the lemma is large and thus in particular it is not empty,
which finishes the proof modulo \eqref{eq:toprove}.

Now \eqref{eq:toprove} follows by an easy modification of the proof of \cite[Lemma 5.11]{winstein2025concentration},
but we present a sketch for completeness.
The key ingredient is the existence of some constant $\delta > 0$ such that, for all
$x \in \gampe \cap \phalfball$, every neighbor
$x'$ of $x$ (meaning $\dh(x,x') \leq 1$) satisfies
\begin{equation}
\label{eq:steplb}
    \P \left[ X_1^x = x' \right] \geq \frac{\delta}{n^2}.
\end{equation}
This is stated as \cite[Equation (15)]{winstein2025concentration}.
Now because $\xi \in (0,1)$, letting $\Xi = \Lambda^+ \cap \gampe \cap \phalfball$,
a union bound of \eqref{eq:lambdaminus} (using the stationarity of the chain $(X_t)$) implies that
\begin{equation}
    \P\left[ X_0, X_1 \in \Xi \right]  \geq 1 - \Exp{-\Omega(n^\xi)}.
\end{equation}
Now because of \eqref{eq:steplb}, we have
\begin{equation}
    \P \left[ X \notin \Lambda \right]
    \leq \P \left[ X_0 \notin \Xi \right] + \frac{n^2}{\delta} \P\left[ X_0 \in \Xi, X_1 \notin \Xi \right],
\end{equation}
because if $X_0 \in \Xi$ and one of the neighbors of $X_0$ is not in $\Xi$,
then the stationary chain will have at least $\frac{\delta}{n^2}$ chance to find it as $X_1$.
This finishes the proof.
\end{proof}

With Lemmas \ref{lem:global_step_concentration} and \ref{lem:goodset_from_bigset} in hand,
we are now ready to finish the main proof of this subsection.

\begin{proof}[Proof of Proposition \ref{prop:global_hajek_intro}/\ref{prop:global_hajek}]
For $\alpha \in (0,\fr{1}{2})$ to be determined later, set
\begin{equation}
    \Pi^\alpha = \left\{ x \in \pball : \left| \sN_G(x,e) - 2 \se p^{\se-1} n^{\sv-2} \right| \leq
    n^{\sv-2.5+\alpha} \text{ for all } e \in \binom{[n]}{2} \right\}.
\end{equation}
By Lemma \ref{lem:global_step_concentration} and a union bound over all $\binom{n}{2}$ edges, we have
\begin{equation}
    \P[X \in \Pi^\alpha] \geq 1 - \Exp{-\Omega(n^{2\alpha})}.
\end{equation}
Now let $\Lambda^\alpha \sse \Lambda^*$ be the set given by applying Lemma \ref{lem:goodset_from_bigset} to
$\Pi^\alpha$ (where we take $\xi = 2\alpha \in (0,1)$), and let $z \in \Lambda^\alpha$ be also as in
Lemma \ref{lem:goodset_from_bigset}.
We will apply Corollary \ref{cor:barbour} with $f = \hsN_G$ and $\Lambda = \Lambda^\alpha$, with this choice of $z$,
and with $T = n^3$.

First, condition \ref{cond:setlarge} holds with $\epsilon = \Exp{-\Omega(n^{2\alpha})}$ by the second part of
Lemma \ref{lem:goodset_from_bigset}.
Next, condition \ref{cond:mixing} follows, with $\delta = \Exp{-\Omega(n)}$,
from Lemma \ref{lem:mixing}, since $z \in \Lambda^*$.
Additionally, we can take $M \lesssim n^\sv$ in condition \ref{cond:fbound}, and so it just remains to consider
conditions \ref{cond:variance} and \ref{cond:smalljumps} and then analyze the resulting bound.

Let us first examine condition \ref{cond:smalljumps}.
We will use the fact that $\Pi^\alpha$ is an interval in the sense that if
$a \preceq b \preceq c$ and $a, c \in \Pi^\alpha$, then $b \in \Pi^\alpha$ as well.
This holds because $\Pi^\alpha$ is the intersection of pullbacks of intervals in $\R$ under monotone functions
$\sN_G(-,e)$.
By the defining feature of $\Pi^\alpha$, this means that if $a, b \in \Pi^\alpha$ and $a \preceq b$,
then
\begin{equation}
\label{eq:ngdiffhamming}
    \left| \hsN_G(a) - \hsN_G(b) \right|
    \leq n^{\sv-2.5+\alpha} \cdot \dh(a, b),
\end{equation}
since we can flip one edge at a time to get from $a$ to $b$, all while remaining in $\Pi^\alpha$.
Indeed, when edge $e$ is flipped, the quantity $2 \se p^{\se-1} n^{\sv-2} \sE(x)$
changes by exactly $2 \se p^{\se-1} n^{\sv-2}$,
and $\sN_G(x)$ changes by $\sN_G(x,e)$, which is controlled by the definition of $\Pi^\alpha$
and is no more than $n^{\sv-2.5+\alpha}$ away from $2 \se p^{\se-1} n^{\sv-2}$.

Now, for any $x \in \Lambda^\alpha$ and $x'$ with $P(x,x') > 0$, both $(X_t^x)$ and $(X_t^{x'})$
will remain in $\Pi^\alpha$ for all $t \leq n^3$ with probability at least $1 - \Exp{-\Omega(n^{2\alpha})}$
by Lemma \ref{lem:goodset_from_bigset}.
Furthermore, since $\Lambda^\alpha \sse \Lambda^*$, by Proposition \ref{prop:goodset} neither
chain will leave $\phalfball$ in time $n^3$ either, with probability at least $1 - \Exp{-\Omega(n)}$.
This is important because it ensures that the ordering between the chains is preserved.
Since either $x \preceq x'$ or $x \succeq x'$ (as $\dh(x,x') \leq 1$ since $P(x,x') > 0$), the same ordering holds for
$X_t^x$ and $X_t^{x'}$ as long as both trajectories remain in $\phalfball$ until time $t$,
as in this case neither chain will have had the opportunity to attempt to leave $\pball$, which is the only
way that the ordering can be interrupted, as discussed in Section \ref{sec:inputs_mixing} (just above
the statement of Proposition \ref{prop:goodset}).

All this means that for $t \leq n^3$, $x \in \Lambda^\alpha$, and $x'$ with $P(x,x') > 0$, we have
\begin{equation}
    \left|
        \E\left[\hsN_G(X_t^x) - \hsN_G(X_t^{x'})\right]
    \right| \leq n^{\sv-2.5+\alpha} \cdot \E \left[ \dh(X_t^x,X_t^{x'}) \right] + \Exp{-\Omega(n^{2\alpha})}.
\end{equation}
And, by Lemma \ref{lem:mixing}, since $\Lambda^\alpha \sse \Lambda^*$ we have
\begin{equation}
    \E \left[ \dh(X_t^x,X_t^{x'}) \right] \leq \left(1 - \fr{\kappa}{n^2}\right)^t + \Exp{-\Omega(n)}.
\end{equation}
Finally, by the last part of Lemma \ref{lem:mixing} (the total variation distance bound), for all $t \geq n^3$ we have
\begin{equation}
    \left|
        \E\left[\hsN_G(X_t^x) - \hsN_G(X_t^{x'})\right]
    \right| \leq \Exp{-\Omega(n)},
\end{equation}
here using the fact that $|\hsN_G(x)| \leq M$ for all $x \in \pball$ with $M \lesssim n^\sv$.
Therefore there is some constant $C$ for which
\begin{equation}
    \left|
        \E\left[\hsN_G(X_t^x) - \hsN_G(X_t^{x'})\right]
    \right| \leq C n^{\sv - 2.5 + \alpha}
\end{equation}
for \emph{all} $t \in \N$, whenever $x \in \Lambda^\alpha$ and $P(x,x') > 0$.
In other words, we obtain condition \ref{cond:smalljumps} with $\Delta \lesssim n^{\sv - 2.5 + \alpha}$.

Finally, we examine condition \ref{cond:variance}.
For any $x \in \Lambda^\alpha$, with probability at least $1 - \Exp{-\Omega(n^{2\alpha})}$, each of the
chains $X_t^{x'}$ for $x'$ with $\dh(x,x') \leq 1$ remain in the set
$\Pi^\alpha \cap \gampe \cap \phalfball$ for all times $0 \leq t < n^3$, using Proposition \ref{prop:goodset}
which gives the properties of $\Lambda^*$, as well as Lemma \ref{lem:goodset_from_bigset}.
So by similar reasoning as before, applying the interval property of $\Pi^\alpha$ and Lemma \ref{lem:mixing},
we find that for all $x \in \Lambda^\alpha$,
\begin{align}
    &\sum_{x' \in \pball} P(x,x') \sum_{t=0}^{n^3-1} \left(
        \E\left[\hsN_G(X_t^x)\right] - \E\left[\hsN_G(X_t^{x'})\right]\right)^2 \\
    &\qquad \qquad \leq \fr{1}{\binom{n}{2}} \sum_{e \in \edgeset} \sum_{t=0}^{n^3-1}
    \left( \E\left[\hsN_G(X_t^x)\right] - \E\left[\hsN_G(X_t^{x^{\oplus e}})\right] \right)^2 \\
    &\qquad \qquad \leq \fr{1}{\binom{n}{2}} \sum_{e \in \binom{[n]}{2}} \sum_{t=0}^{n^3-1}
    \left(
        n^{\sv-2.5+\alpha} \cdot \E\left[\dh(X_t^x,X_t^{x^{\oplus e}})\right]
    \right)^2 + \Exp{-\Omega(n^{2\alpha})} \\
    &\qquad \qquad \lesssim n^{2\sv-5+2\alpha} \sum_{t=0}^{n^3-1} \left(
        1 - \fr{\kappa}{n^2}
    \right)^{2t} + \Exp{-\Omega(n^{2\alpha})} \\
    &\qquad \qquad \lesssim n^{2\sv-3+2\alpha},
\end{align}
meaning that we can take $V \lesssim n^{2\sv-3+2\alpha}$.
Recall our notation $x^{\oplus e}$ for the graph $x$ with the state of $e$ flipped.
Additionally, as in the proof of Lemma \ref{lem:basic_concentration},
it should be understood that in the sum over $x' \in \pball$, we only consider $n$-vertex graphs
whose graphon representations lie in $\pball$.

To summarize, we can apply Corollary \ref{cor:barbour} with
\begin{equation}
    V = C n^{2\sv-3+2\alpha}, \quad
    \Delta = C n^{\sv-2.5+\alpha}, \quad
    \epsilon = \Exp{-cn^{2\alpha}}, \quad
    \delta = \Exp{-cn}, \quad
    M = C n^\sv,
\end{equation}
for some $C, c > 0$.
The conclusion of Corollary \ref{cor:barbour} thus yields, for all $\lambda \geq 2 C n^\sv e^{-cn}$,
\begin{equation}
    \P\left[
        \left|
            \hsN_G(X) - \E[\hsN_G(X)]
        \right| > \lambda
    \right] \leq 2 \Exp{- \fr{\left(\lambda - 2 C n^\sv e^{-cn}\right)^2}{2 C n^{2\sv-3+2\alpha}
    + \fr{4}{3} C n^{\sv-2.5+\alpha} (\lambda - 2 C n^\sv e^{-cn})}}
    + e^{-cn^{2\alpha}} + e^{-cn}.
\end{equation}
Now let us plug in $\lambda = n^{\sv-1.5+\zeta}$.
Note that the second term in the denominator inside the exponential is always dominated by the first one, since
\begin{equation}
    2\sv - 3 + 2\alpha \geq 2\sv - 4 + \alpha + \zeta,
\end{equation}
which holds as $\alpha > 0$ and $\zeta < 1$.
So, combining terms and using the fact that $2\alpha < 1$ and that the
$C n^\sv e^{-cn}$ terms are exponentially small, we find that
\begin{equation}
    \P\left[
        \left|
            \hsN_G(X) - \E[\hsN_G(X)]
        \right| > n^{\sv-1.5+\zeta}
    \right] \leq \Exp{- cn^{2\zeta-2\alpha}}
    + \Exp{-cn^{2\alpha}};
\end{equation}
here we have changed the value of $c$.
Setting $\alpha = \fr{\zeta}{2}$ gives the optimal bound here, which is $\Exp{- cn^\zeta}$,
finishing the proof.
\end{proof}

As previously mentioned, the only place where we used the fact that $G$ and $G_0, G_1, \dotsc, G_K$ are forests
in the phase coexistence case is in approximating $\E[\sN_G(X,e)]$ by $2 \se p^{\se-1} n^{\sv-2}$
via Propositions \ref{prop:Eproduct} and \ref{prop:marginal}.
Thus, the same proof, without this step, yields the same statement as in
Proposition \ref{prop:global_hajek_intro}/\ref{prop:global_hajek}, but with $\hsN_G(X)$ replaced by
\begin{equation}
\label{eq:worsenhat}
    \sN_G(X) - \E[\sN_G(X,e)] \cdot \sE(X).
\end{equation}
This will not be useful for proving the CLT for $\sE(X)$ itself, because the proof relies on certain cancellations
between various algebraic formulas.
Nevertheless, it does allow us to deduce a CLT (albeit with nonexplicit fluctuation rate) for
arbitrary $\sN_G(X)$ once we have one for $\sE(X)$; see Section \ref{sec:hajek_subgraphs} for more details.
\subsection{Local subgraph counts}
\label{sec:hajek_local}

We next proceed with the proof of Proposition \ref{prop:local_hajek_intro}, the local version of
Proposition \ref{prop:global_hajek_intro}/\ref{prop:global_hajek}, whose statement has been reproduced below
for the reader's convenience.
Throughout this section, we fix some $\rho \in \cV$ and $v \in [n]$, and we remind the reader that, as first
introduced in Section \ref{sec:intro_iop}, the notation $\sN_G^{\rho \to v}(x)$ denotes the number of
homomorphisms of $G$ in $x$ which map $\rho$ to $v$.

\begin{proposition}[Restatement of Proposition \ref{prop:local_hajek_intro}]
\label{prop:local_hajek}
Let us make Assumption \ref{assumption}.
As in \eqref{eq:hatngrhovdef}, let us write
\begin{equation}
    \hsN_G^{\rho \to v}(x) = \sN_G^{\rho \to v}(x) - \sd_\rho p^{\se-1} n^{\sv-2} \deg_v(x).
\end{equation}
Then for all $\zeta \in (0,\fr{3}{4})$, there is some $c > 0$ such that
\begin{equation}
    \P \left[ \left| \hsN_G^{\rho \to v}(X) - \E \left[ \hsN_G^{\rho \to v}(X) \right]\right| > n^{\sv-1.75+\zeta} \right]
    \leq \Exp{-c n^{\fr{4}{3} \zeta}}.
\end{equation}
In particular, for all $\xi > 0$ we have
\begin{equation}
    \Var \left[
        \hsN_G^{\rho \to v}(X)
    \right] \lesssim n^{2\sv - 3.5 + \xi}.
\end{equation}
\end{proposition}

The proof follows the same general strategy as for Proposition \ref{prop:global_hajek_intro}/\ref{prop:global_hajek},
but there are some more technical points encountered along the way due to the local nature of the proof.
In particular, to make Corollary \ref{cor:barbour} work as desired, we need to use Proposition \ref{prop:dloc_uniformbound_intro},
which bounds $\dloc(X_t^x, X_t^{x'})$, the discrepancy within a neighborhood around $v$ between two monotonically-coupled
Glauber chains started at $x$ and $x'$ which differ at an edge \emph{not adjacent to} $v$.
In essence, this gives an approximate version of \eqref{eq:distancegrows}, as will be seen in the proof details.

In any case, we begin with the following local version of Lemma \ref{lem:global_step_concentration}, showing that
the changes $\partial_e \sN_G^{\rho \to v}(X)$ under Glauber dynamics are themselves concentrated around the correct
value for $e \ni v$.
Note that the change $\partial_e \sN_G^{\rho \to v}(x)$ is the same as $\sN_G^{\rho \to v}(x,e)$, the number of 
homomorphisms of $G$ in $x^{+e}$ which map $\rho$ to $v$ and use the edge $e$.
As was the case in the proof of Proposition \ref{prop:global_hajek_intro}/\ref{prop:global_hajek}, the only place
where the ``phase-uniqueness or forest'' assumption is used is in approximating the mean
$\E[\sN_G^{\rho \to v}(X,e)]$.

\begin{lemma}
\label{lem:local_step_concentration}
Let us make Assumption \ref{assumption}.
For $\alpha \in (0,\fr{1}{2})$ there is some $c > 0$ such that, for all $e \in \edgeset$ with $v \in e$,
\begin{equation}
    \P \left[ | \sN^{\rho \to v}_G(X,e) - \sd_\rho p^{\se-1} n^{\sv-2} | > n^{\sv - 2.5 + \alpha} \right]
    \leq \Exp{- c n^{2\alpha}}.
\end{equation}
\end{lemma}

\begin{proof}[Proof of Lemma \ref{lem:local_step_concentration}]
As in the proof of Lemma \ref{lem:global_step_concentration}, we will first invoke
Theorem \ref{thm:lipschitzconcentration} to obtain concentration for $\sN_G^{\rho \to v}(X,e)$ around its mean,
and then estimate the mean.
So we need to construct a Lipschitz vector $\L$
for the function $x \mapsto \sN_G^{\rho \to v}(x,e)$.
Again, we can take $\L_e = 0$; next for any other $e' \in \edgeset$, we have
\begin{equation}
    \sN_G^{\rho \to v} (x^{+e'}, e) - \sN_G^{\rho \to v} (x^{-e'}, e) = \sN_G^{\rho \to v}(x,e,e'),
\end{equation}
which is the number of homomorphisms $G \to x^{+e+e'}$ which map $\rho$ to $v$ and contain both $e$ and $e'$ as
images of edges in $G$.
If $e$ and $e'$ do not share a vertex, this is $\lesssim n^{\sv-4}$ and if they do share a vertex then it's
$\lesssim n^{\sv-3}$ (recall that $v \in e$).
This means we can actually take $\L$ to be the same Lipschitz vector as obtained in the
proof of Lemma \ref{lem:global_step_concentration}, up to constants, and thus we obtain the same concentration for
$\sN_G^{\rho \to v}(X,e)$ \emph{around its mean} as for $\sN_G(X,e)$.

It just remains to estimate the mean; again we invoke Propositions \ref{prop:Eproduct} and \ref{prop:marginal}
and obtain
\begin{equation}
    \E \left[ \sN_G^{\rho \to v} (X, e) \right] = \sd_\rho p^{\se-1} n^{\sv-2} + O(n^{\sv-2.5} \sqrt{\log n}),
\end{equation}
via the same argument as in Lemma \ref{lem:global_step_concentration}, although this time instead of $2\se$
options for an oriented edge to map onto $e$, we only have $\sd_\rho$ options since we already know that $\rho$
maps to $v$ and we just need to pick one neighbor of $\rho$ to map to the other vertex in $e$.
In any case, this leads to the desired conclusion as before.
\end{proof}

With Lemma \ref{lem:local_step_concentration} in hand, we may again use Corollary \ref{cor:barbour}
to conclude the proof.
We will also use Lemma \ref{lem:goodset_from_bigset} as in the previous section, and moreover will make use
of Proposition \ref{prop:dloc_uniformbound_intro}, bounding the local Hamming distance around $v$ in monotonically-coupled 
Glauber chains.

\begin{proof}[Proof of Proposition \ref{prop:local_hajek_intro}/\ref{prop:local_hajek}]
Let $\alpha \in (0,\fr{1}{2})$ be determined later and define 
\begin{equation}
    \Pi^\alpha = \left\{
        x \in \pball: 
        \left| \sN_G^{\rho \to v}(x,e) - \sd_\rho p^{\se-1} n^{\sv-2} \right| \leq n^{\sv-2.5+\alpha}
        \text{ for all } e \in \edgeset \text{ with } v \in e
    \right\}.
\end{equation}
By Lemma \ref{lem:local_step_concentration}, the triangle inequality, and a union bound, we have
$\P[X \in \Pi^\alpha] \geq 1 - \Exp{-\Omega(n^{2\alpha})}$.
Let $\Lambda^\alpha$ be the set guaranteed by Lemma \ref{lem:goodset_from_bigset} applied with $\Pi = \Pi^\alpha$.
We will apply Corollary \ref{cor:barbour} with $f = \hsN_G^{\rho \to v}$ and $\Lambda = \Lambda^\alpha$, and with $T = n^3$.

As in the proof of Proposition \ref{prop:global_hajek}, we can take $\epsilon = \Exp{-\Omega(n^{2\alpha})}$,
$\delta = \Exp{-\Omega(n)}$, and $M \lesssim n^\sv$, and it remains to check conditions \ref{cond:variance}
and \ref{cond:smalljumps}.

For this, first note that if $a, b \in \Pi^\alpha$ differ only in an edge $e$ with $v \in e$,
then
\begin{equation}
\label{eq:diffatv}
    \left|\hsN_G^{\rho \to v}(a) - \hsN_G^{\rho \to v}(b)\right| \leq n^{\sv-2.5+\alpha},
\end{equation}
by the definition of $\Pi^\alpha$, since
\begin{equation}
    \left| \hsN_G^{\rho \to v}(a) - \hsN_G^{\rho \to v}(b)\right|
    = \left| \sN_G^{\rho \to v}(a,e) - \sd_\rho p^{\se-1} n^{\sv-2} \right|.
\end{equation}
On the other hand, if $a,b$ differ only in an edge $e$ with $v \notin e$, then we have
\begin{equation}
\label{eq:diffawayv}
    \left| \hsN_G^{\rho \to v}(a) - \hsN_G^{\rho \to v}(b)\right|
    = \left| \sN_G^{\rho \to v} (a,e) \right| \leq C n^{\sv-3},
\end{equation}
for some $C$, since there are already three vertices accounted for ($v$, and the two other vertices of $e$) in any
copy of $G$.
Now again the set $\Pi^\alpha$ is an interval in the same sense as in the proof of Proposition \ref{prop:global_hajek}.
So the facts \eqref{eq:diffatv} and \eqref{eq:diffawayv} imply that for any $a, b \in \Pi^\alpha$ with $a \preceq b$
or $b \preceq a$, we have
\begin{equation}
\label{eq:local_hajek_mainbound}
    \left| \hsN_G^{\rho \to v} (a) - \hsN_G^{\rho \to v} (b) \right|
    \leq n^{\sv-2.5+\alpha} \dloc(a,b) + C n^{\sv-3} \dh(a,b),
\end{equation}
where we recall from \eqref{eq:dlocdef} the definition of the ``local Hamming distance at $v$'' as follows:
\begin{equation}
    \dloc(a,b) = \sum_{e \ni v} |a(e)-b(e)|.
\end{equation}
We will use the inequality \eqref{eq:local_hajek_mainbound} to verify both conditions \ref{cond:variance}
and \ref{cond:smalljumps}.

First, for condition \ref{cond:smalljumps}, since $\dloc(a,b) \leq \dh(a,b)$ and $n^{\sv-2.5+\alpha} \gg n^{\sv-3}$,
we have
\begin{equation}
\label{eq:local_hajek_mainbound_simplified}
    \left| \hsN_G^{\rho \to v}(a) - \hsN_G^{\rho \to v}(b) \right| \leq C n^{\sv-2.5+\alpha} \dh(a,b)
\end{equation}
for some $C$, whenever $a,b \in \Pi^\alpha$ and $a \preceq b$ or vice versa.
Up to a constant, this is the same inequality \eqref{eq:ngdiffhamming} which held in the case of $\hsN_G$ in the proof of
Proposition \ref{prop:global_hajek}, and so the same reasoning carries through, allowing us to conclude
that condition \ref{cond:smalljumps} holds with $\Delta \lesssim n^{\sv-2.5+\alpha}$.

Next, for condition \ref{cond:variance}, as mentioned in Sections \ref{sec:inputs_concentration}
and \ref{sec:inputs_dloc}, we will need to prove a version of the inequality \eqref{eq:distancegrows},
using Proposition \ref{prop:dloc_uniformbound_intro}.
First, for all $x \in \Lambda^\alpha$ we have
\begin{align}
    &\sum_{x' \in \pball} P(x,x') \sum_{t=0}^{n^3-1} \left(
        \E\left[\hsN_G^{\rho \to v}(X_t^x)\right] - \E\left[\hsN_G^{\rho \to v}(X_t^{x'})\right]\right)^2 \\
    &\qquad \qquad \leq \fr{1}{\binom{n}{2}} \sum_{e \in \edgeset} \sum_{t=0}^{n^3-1}
    \left( \E\left[\hsN_G^{\rho \to v}(X_t^x)\right] - \E\left[\hsN_G^{\rho \to v}(X_t^{x^{\oplus e}})\right] \right)^2,
    \label{eq:local_hajek_sumovere}
\end{align}
recalling the notation $x^{\oplus e}$ for the graph $x$ with the status of edge $e$ flipped.
For the terms in \eqref{eq:local_hajek_sumovere} corresponding to $e \ni v$, we may use
\eqref{eq:local_hajek_mainbound_simplified} and the fact that $X_t^x$ and $X_t^{x^{\oplus e}}$ will remain ordered
and within $\Pi^\alpha \cap \gampe \cap \phalfball$ for all time up to $n^3-1$ with probability at least
$1 - \Exp{-\Omega(n^{2\alpha})}$ in order to obtain, for all $e \ni v$,
\begin{align}
    &\sum_{t=0}^{n^3-1} \left( \E\left[\hsN_G^{\rho \to v}(X_t^x)\right] - \E\left[\hsN_G^{\rho \to v}(X_t^{x^{\oplus e}})\right] \right)^2 \\
    &\qquad \qquad \lesssim n^{2\sv-5+2\alpha} \sum_{t=0}^{n^3-1} \left( \E\left[\dh(X_t^x,X_t^{x^{\oplus e}})\right]\right)^2 + \Exp{-\Omega(n^{2\alpha})} \\
    &\qquad \qquad \lesssim n^{2\sv-3+2\alpha}.
    \label{eq:ubwithv}
\end{align}
Here we have invoked Lemma \ref{lem:mixing}.

On the other hand, for the terms corresponding to $e \not\ni v$, we have
\begin{align}
    &\left|\E\left[\hsN_G^{\rho \to v}(X_t^x)\right] - \E\left[\hsN_G^{\rho \to v}(X_t^{x^{\oplus e}})\right] \right| \\
    &\qquad \qquad\leq n^{\sv-2.5+\alpha} \cdot \E \left[ \dloc(X_t^x, X_t^{x^{\oplus e}}) \right]
    + C n^{\sv-3} \cdot \E \left[ \dh(X_t^x, X_t^{x^{\oplus e}}) \right] + \Exp{-\Omega(n^{2\alpha})}.
\end{align}
Now, by Proposition \ref{prop:dloc_uniformbound_intro}, since $x$ and $x^{\oplus e}$ differ only at some edge $e \not \ni v$
and also $x \in \Lambda^*$, we have
\begin{equation}
    \E \left[ \dloc(X_t^x, X_t^{x^{\oplus e}}) \right] \lesssim \fr{1}{n}
\end{equation}
Since additionally $\E\left[\dh(X_t^x, X_t^{x^{\oplus e}})\right] \lesssim 1$ and $\alpha < \fr{1}{2}$, we obtain
\begin{equation}
\label{eq:approxdistancegrows}
    \left|\E\left[\hsN_G^{\rho \to v}(X_t^x)\right] - \E\left[\hsN_G^{\rho \to v}(X_t^{x^{\oplus e}})\right] \right|
    \lesssim n^{\sv-3}
\end{equation}
in this case.
This is the aforementioned version of \eqref{eq:distancegrows} we need, since the difference at $t=0$ is
$\lesssim n^{\sv-3}$ as well.
With \eqref{eq:approxdistancegrows}, we obtain
\begin{align}
    &\sum_{t=0}^{n^3-1} \left( \E\left[\hsN_G^{\rho \to v}(X_t^x)\right] - \E\left[\hsN_G^{\rho \to v}(X_t^{x^{\oplus e}})\right] \right)^2 \\
    &\qquad \qquad \lesssim n^{\sv-3} \sum_{t=0}^{n^3-1} \left| \E\left[\hsN_G^{\rho \to v}(X_t^x)\right] - \E\left[\hsN_G^{\rho \to v}(X_t^{x^{\oplus e}})\right] \right| \\
    &\qquad \qquad \lesssim n^{\sv-3} \sum_{t=0}^{n^3-1} n^{\sv-2.5+\alpha} \cdot \E \left[ \dh(X_t^x, X_t^{x^{\oplus e}})\right] + \Exp{-\Omega(n^{2\alpha})} \\
    &\qquad \qquad \lesssim n^{2\sv-3.5+\alpha},
    \label{eq:ubwithoutv}
\end{align}
for all $e \not \ni v$.

Now, applying \eqref{eq:ubwithv} and \eqref{eq:ubwithoutv} in \eqref{eq:local_hajek_sumovere}, we find that
\begin{align}
    &\sum_{x' \in \pball} P(x,x') \sum_{t=0}^{n^3-1} \left(
        \E\left[\hsN_G^{\rho \to v}(X_t^x)\right] - \E\left[\hsN_G^{\rho \to v}(X_t^{x'})\right]\right)^2 \\
    &\qquad \qquad \lesssim \fr{1}{n^2} \left( n \cdot n^{2\sv-3+2\alpha} + n^2 \cdot n^{2\sv-3.5+\alpha} \right) \\
    &\qquad \qquad \lesssim n^{2\sv-3.5+\alpha},
\end{align}
using again the fact that $\alpha < \fr{1}{2}$.
In other words, condition \ref{cond:variance} holds with $V \lesssim n^{2\sv-3.5+\alpha}$.

To summarize, all conditions of Corollary \ref{cor:barbour} hold with the following choices
\begin{equation}
    V = C n^{2\sv-3.5+\alpha}, \quad
    \Delta = C n^{\sv-2.5+\alpha}, \quad
    \epsilon = \Exp{-cn^{2\alpha}}, \quad
    \delta = \Exp{-cn}, \quad
    M = C n^\sv,
\end{equation}
for some constants $C$ and $c$.
The conclusion of Corollary \ref{cor:barbour} thus yields, for $\lambda \geq 2 C n^\sv e^{-cn}$,
\begin{equation}
    \P\left[
        \left|
            \hsN_G^{\rho \to v}(X) - \E\left[\hsN_G^{\rho \to v}(X)\right]
        \right| > \lambda
    \right] \leq 2 \Exp{- \fr{\left(\lambda - 2 C n^\sv e^{-cn}\right)^2}{2 C n^{2\sv-3.5+\alpha}
    + \fr{4}{3} C n^{\sv-2.5+\alpha} (\lambda - 2 C n^\sv e^{-cn})}}
    + e^{-cn^{2\alpha}} + e^{-cn}.
\end{equation}
Plugging in $\lambda = n^{\sv-1.75+\zeta}$, the second term in the denominator is dominated by the first exactly when
\begin{equation}
    2\sv-3.5+\alpha \geq 2\sv-4.25+\alpha+\zeta,
\end{equation}
and so we must take $\zeta \leq 0.75$ for this to hold (this restriction does not impose any additional constraints
on $\alpha$).
In the case where $\zeta \leq 0.75$, we thus find that
\begin{equation}
    \P\left[
        \left|
            \hsN_G^{\rho \to v}(X) - \E\left[\hsN_G^{\rho \to v}(X)\right]
        \right| > n^{\sv-1.75+\zeta}
    \right] \leq \Exp{- c n^{2\zeta-\alpha}} + \Exp{-cn^{2\alpha}},
\end{equation}
where we have changed the value of $c$.
The optimal choice $\alpha$ here is $\alpha = \fr{2}{3} \zeta$, which is in $(0,\fr{1}{2})$
as $\zeta \in (0,\fr{3}{4})$, and this yields an upper bound of $\Exp{- c n^{\fr{4}{3} \zeta}}$,
finishing the proof.
\end{proof}

Next, we consider the count $\sN_G^v(x)$, which is simply the number of homomorphisms of $G$ in $x$
which map \emph{some vertex} of $G$ to $v \in [n]$.
Notice first that
\begin{equation}
    \sN_G^v(x) = \sum_{\rho \in \cV} \sN_G^{\rho \to v}(x) + O(n^{\sv-2}).
\end{equation}
To see why this holds, first note that we would have exact equality without the $O(n^{\sv-2})$ error term if
we required all homomorphisms counted above to be \emph{injective}.
Thus, we obtain this approximate equality since there are $\lesssim n^{\sv-2}$ non-injective homomorphisms
which map \emph{some vertex} of $G$ to $v \in [n]$.
Since $n^{\sv-2} \ll n^{\sv-1.75+\zeta}$ for any $\zeta > 0$ and $\sum_{\rho \in \cV} \sd_\rho = 2 \se$ by
the handshake lemma,
Proposition \ref{prop:local_hajek_intro}/\ref{prop:local_hajek} immediately implies the following corollary
via a union bound.

\begin{corollary}
\label{cor:local_hajek_norho}
Let us make Assumption \ref{assumption}.
Define
\begin{equation}
    \hsN_G^v(x) \coloneqq \sN_G^v(x) - 2 \se p^{\se-1} n^{\sv-2} \deg_v(x).
\end{equation}
Then for all $\zeta \in (0,\fr{3}{4})$, there is some $c > 0$ such that
\begin{equation}
    \P \left[
        \left|
            \hsN_G^v(X) - \E\left[\hsN_G^v(X)\right]
        \right| > n^{\sv-1.75+\zeta}
    \right] \leq \Exp{-c n^{\fr{4}{3}\zeta}}.
\end{equation}
\end{corollary}

In other words, the fluctuations of the overall number of subgraphs at a vertex are also driven by the
fluctuations in the degree of that vertex.

Finally, we mention that by the same reasoning outlined at the end of Section \ref{sec:hajek_global},
we may remove the assumption that $G$ is a forest in the phase coexistence case, at the cost of not approximating
$\E[\sN_G^{\rho \to v}(X,e)]$ or $\E[\sN_G^v(X,e)]$.
In other words, in this case we still have the same statement as Proposition \ref{prop:local_hajek_intro}/\ref{prop:local_hajek}
but with $\hsN_G^{\rho \to v}(X)$ replaced by
\begin{equation}
\label{eq:worsenrhohat}
    \sN_G^{\rho \to v}(X) - \E[\sN_G^{\rho \to v}(X,e)] \cdot \deg_v(X).
\end{equation}
Additionally, we still have the same statement as Corollary \ref{cor:local_hajek_norho} but with
$\hsN_G^v(X)$ replaced by
\begin{equation}
\label{eq:worsenvhat}
    \sN_G^v(X) - \E[\sN_G^v(X,e)] \cdot \deg_v(X).
\end{equation}
As before, these results will not be useful for proving the CLT for $\deg_v(X)$, but they do allow us to
conclude a CLT for $\sN_G^{\rho \to v}(X)$ or $\sN_G^v(X)$ from a CLT for $\deg_v(X)$.
See Section \ref{sec:hajek_subgraphs} for more details.
\subsection{Contingent CLT for subgraph counts}
\label{sec:hajek_subgraphs}

In this section, we finish the proofs of Corollaries \ref{cor:global_subgraph_clt_intro}
and \ref{cor:local_subgraph_clt_intro}, under the assumption of Theorems \ref{thm:global_clt_intro}
and \ref{thm:local_clt_intro}.
It should be noted that the variance bound given in
Proposition \ref{prop:global_hajek_intro}/\ref{prop:global_hajek} implies, via Chebyshev's
inequality, that if $\sE(X)$ has a CLT when the fluctuations are normalized by some quantity $\sigma_n$
which is of order at least $n$, then $\sN_G(X)$ has a CLT (in the non-quantitative sense of convergence in distribution)
when the fluctuations are normalized by $2 \se p^{\se-1} n^{\sv-2} \sigma_n$,
and a similar statement holds for the local subgraph counts.

The following corollary uses the tail bound in Proposition \ref{prop:global_hajek_intro}/\ref{prop:global_hajek}
to upgrade this non-quantitative statement to a bound on the Wasserstein error in the CLT
for $\sN_G(X)$, given the Wasserstein distance in the CLT for $\sE(X)$.
This is a corollary solely of Proposition \ref{prop:global_hajek_intro}/\ref{prop:global_hajek},
but together with Theorem \ref{thm:global_clt_intro} it implies Corollary \ref{cor:global_subgraph_clt_intro}.
For now, we assume that we are either in the phase uniqueness case, or that $G$ and $G_0, G_1, \dotsc, G_K$ are forests.
This assumption can be removed at the price of an expression for the fluctuations which has no closed form;
this will be clarified at the end of this section.

\begin{corollary}
\label{cor:global_subgraph_clt}
Let us make Assumption \ref{assumption}.
Assume that there are sequences $\sigma_n \gtrsim n$ and $\delta_n \to 0$ such that
\begin{equation}
    \dW\left(\frac{\sE(X) - \E[\sE(X)]}{\sigma_n}, Z \right) \leq \delta_n,
\end{equation}
where $Z \sim \cN(0,1)$.
Then for all $\eps > 0$,
\begin{equation}
    \dW\left(\frac{\sN_G(X) - \E[\sN_G(X)]}{2 \se p^{\se-1} n^{\sv-2} \sigma_n}, Z\right)
    \lesssim n^{-\fr{1}{2}+\eps} + \delta_n.
\end{equation}
\end{corollary}

\begin{proof}[Proof of Corollary \ref{cor:global_subgraph_clt}]
We may assume that $\eps < 1$.
It suffices to show that
\begin{equation}
    \E \left[ \left| \frac{\sN_G(X) - \E[\sN_G(X)]}{2 \se p^{\se-1} n^{\sv-2} \sigma_n}
    - \fr{\sE(X) - \E[\sE(X)]}{\sigma_n} \right| \right] \lesssim n^{-\fr{1}{2}+\eps},
\end{equation}
using the ``optimal coupling'' form of the Wasserstein distance, which was mentioned at the beginning of
Section \ref{sec:intro_results}.
Rearranging, this is exactly
\begin{equation}
\label{eq:subgrapheq1}
    \fr{1}{2 \se p^{\se-1} n^{\sv-2} \sigma_n}
    \E\left[ \left| \hsN_G(X) - \E[\hsN_G(X)]\right| \right].
\end{equation}
Now, by Proposition \ref{prop:global_hajek}, since $\hsN_G(x) \lesssim n^\sv$
uniformly for all $x$, we have
\begin{align}
    \E\left[ \left| \hsN_G(X) - \E[\hsN_G(X)]\right| \right]
    &\lesssim n^{\sv-1.5+\eps} + n^\sv \cdot \P \left[ \left| \hsN_G(X) - \E[\hsN_G(X)] \right| > n^{\sv-1.5+\eps} \right] \\
    &\leq n^{\sv-1.5+\eps} + n^\sv \cdot \Exp{- c n^{\eps}} \\
    &\lesssim n^{\sv-1.5+\eps}.
\end{align}
Therefore \eqref{eq:subgrapheq1} is
\begin{equation}
    \lesssim \fr{n^{\sv-1.5+\eps}}{n^{\sv-2} \sigma_n} \lesssim n^{-\fr{1}{2}+\eps},
\end{equation}
using the fact that $\sigma_n \gtrsim n$.
\end{proof}

The same proof, but with an application of Proposition \ref{prop:local_hajek_intro}/\ref{prop:local_hajek}
instead of Proposition \ref{prop:global_hajek_intro}/\ref{prop:global_hajek}, yields the following local version
of Corollary \ref{cor:global_subgraph_clt}.

\begin{corollary}
\label{cor:local_subgraph_clt}
Let us make Assumption \ref{assumption}.
Assume that there are sequences $\varsigma_n \gtrsim \sqrt{n}$ and $\delta_n \to 0$ such that
\begin{equation}
    \dW\left(\frac{\deg_v(X) - \E[\deg_v(X)]}{\varsigma_n}, Z \right) \leq \delta_n,
\end{equation}
where $Z \sim \cN(0,1)$.
Then for all $\eps > 0$,
\begin{equation}
    \dW\left(\frac{\sN_G^{\rho \to v}(X) - \E[\sN_G^{\rho \to v}(X)]}{\sd_\rho p^{\se-1} n^{\sv-2} \varsigma_n}, Z\right)
    \lesssim n^{-\fr{1}{4}+\eps} + \delta_n.
\end{equation}
\end{corollary}

Another local version of Corollary \ref{cor:global_subgraph_clt} is provided by the following, where we instead
apply Corollary \ref{cor:local_hajek_norho} and consider the count $\sN_G^v(x)$ of homomorphisms which map
\emph{some vertex} to $v$.
The following corollary, combined with Theorem \ref{thm:local_clt_intro}, implies
Corollary \ref{cor:local_subgraph_clt_intro}.

\begin{corollary}
\label{cor:local_subgraph_clt_norho}
Under the same assumptions as in Corollary \ref{cor:local_subgraph_clt},
for all $\eps > 0$,
\begin{equation}
    \dW\left(\frac{\sN_G^{v}(X) - \E[\sN_G^v(X)]}{2 \se p^{\se-1} n^{\sv-2} \varsigma_n}, Z\right)
    \lesssim n^{-\fr{1}{4}+\eps} + \delta_n.
\end{equation}
\end{corollary}

Note that in all of these corollaries, we do not provide bounds on the Kolmogorov distance for the subgraph
counts which match the corresponding Kolmogorov distance bounds obtained in Theorems \ref{thm:global_clt_intro}
and \ref{thm:local_clt_intro}.
However, as previously mentioned in Section \ref{sec:intro_results}, one may apply the standard inequality
\begin{equation}
    \dK(f(X),Z) \lesssim \sqrt{\dW(f(X),Z)}
\end{equation}
to obtain quantitative bounds on the Kolmogorov distance.
These are not sharp, but neither are the Kolmogorov distance bounds we obtain in Theorems \ref{thm:global_clt_intro}
or \ref{thm:local_clt_intro}, at least in principle.
Indeed, as already mentioned in Section \ref{sec:intro_results},
simulations from \cite{winstein2025concentration} show that we should expect
\begin{equation}
    \dK\left(\fr{\sE(X) - \E[\sE(X)]}{\tilde{\sigma}_n}, Z \right) \lesssim \fr{1}{n}
    \qquad \text{and} \qquad
    \dK\left(\fr{\deg_v(X) - \E[\deg_v(X)]}{\tilde{\varsigma}_n}, Z\right) \lesssim \fr{1}{\sqrt{n}}
\end{equation}
for some $\tilde{\sigma}_n^2$ and $\tilde{\varsigma}_n^2$, although these \emph{may not be the same} as
the $\sigma_n^2$ and $\varsigma_n^2$ given by \eqref{eq:sigmadef_global_intro} and \eqref{eq:sigmadef_local_intro}.

Finally, as previously mentioned in Section \ref{sec:intro_results_pcc},
we may conclude CLTs for $\sN_G(X)$, $\sN_G^{\rho \to v}(X)$, and $\sN_G^v(X)$
once we have CLTs for $\sE(X)$ and $\deg_v(X)$, even in the phase coexistence regime when $G$ is not a forest.
For this, we must use the versions of Propositions \ref{prop:global_hajek_intro}/\ref{prop:global_hajek} and
\ref{prop:local_hajek_intro}/\ref{prop:local_hajek}, as well as Corollary \ref{cor:local_hajek_norho},
with the relevant quantities $\hsN_G(X)$, $\hsN_G^{\rho \to v}(X)$, and $\hsN_G^v(X)$ replaced by the expressions
\eqref{eq:worsenhat}, \eqref{eq:worsenrhohat}, and \eqref{eq:worsenvhat} respectively.
This replaces the explicit expressions $2 \se p^{\se-1} n^{\sv-2}$ and $\sd_\rho p^{\se-1} n^{\sv-2}$
which appear as coefficients of $\sigma_n$ and $\varsigma_n$ in Corollaries \ref{cor:global_subgraph_clt},
\ref{cor:local_subgraph_clt}, and \ref{cor:local_subgraph_clt_norho}, with $\E[\sN_G(X,e)]$ and
$\E[\sN_G^{\rho \to v}(X,e)]$, respectively.
In other words, we cannot approximate the size of the fluctuations in this case by an explicit expression, but
we can still prove that they are Gaussian in the limit.
\section{Stein's method for nonlinear exponential families}
\label{sec:stein}

Finally we are ready to prove our main results, the central limit theorems for $\sE(X)$ and $\deg_v(X)$.
Note that Theorems \ref{thm:global_clt_intro} and \ref{thm:local_clt_intro} in their full generality as
discussed in Section \ref{sec:intro_results} are implied by the same statements when $X \sim \condmeas$,
by the large deviations principle given in Theorems \ref{thm:LDP_general} and \ref{thm:LDP_ferromagnetic}
and the discussion thereafter; recall the definition of $\condmeas$ from \eqref{eq:condmeas}.

As mentioned in Section \ref{sec:intro_iop}, the proof uses an application of Stein's method for
\emph{nonlinear exponential families} which was recently developed by \cite{fang2024normal}.
The statement of this theorem is provided in Section \ref{sec:stein_statement} below.
We will need a mild modification of the Hamiltonian which will be discussed thereafter, in
Section \ref{sec:stein_hamiltonian}.
Then in Section \ref{sec:stein_global}
the CLT for $\sE(X)$ (Theorem \ref{thm:global_clt_intro}) will be proved,
followed by Section \ref{sec:stein_local} which contains the proof of the CLT for $\deg_v(X)$
(Theorem \ref{thm:local_clt_intro}).

\subsection{Statement of the theorem}
\label{sec:stein_statement}

We warn the reader that our notation differs from that of \cite{fang2024normal}
by \emph{swapping the roles of the letters $X$ and $Y$}.
In \cite{fang2024normal}, $X$ was a sample from $\cG(n,p)$ and $Y$ was a sample from the ERGM,
but we choose the opposite convention.
This is primarily because we have previously denoted the ERGM sample using $X$, and also
because the ERGM is the primary object of study, whereas $\cG(n,p)$ is seen more as
an auxiliary distribution.

The statement is general, and can be stated for arbitrary random variables as well,
but we state things in terms of discrete variables for the present work.
Let $Y = (Y_1, \dotsc, Y_N)$ be a vector of independent discrete random variables,
and let $X = (X_1, \dotsc, X_N)$ be an exponentially tilted random vector, i.e.\
\begin{equation}
\label{eq:tilt}
    \P[X = x] \propto e^{\sg(x)} \P[Y = x],
\end{equation}
for some function $\sg$.
In our applications, $\sg$ will be a suitable modification of the ERGM Hamiltonian $\sH$.
We implicitly assume that $\E[e^{\sg(Y)}] < \infty$.
In addition, \emph{we assume that the samples $X$ and $Y$ are independent}.

Now let $f$ be a function with $\E[f(X)] = 0$.
Although it is not precisely required by our definitions, $f$ should be
appropriately normalized so as to have variance close to $1$.
The theorem of \cite{fang2024normal} gives explicit expressions for error terms
$\delta_0, \delta_1, \delta_1', \delta_2, \delta_3$ which bound
\begin{equation}
    \dW(f(X),Z) \qquad \text{and} \qquad \dK(f(X),Z),
\end{equation}
where $Z \sim \cN(0,1)$ (recall the definitions of these distances from \eqref{eq:wasdist} and \eqref{eq:koldist}).
See the statement of Theorem \ref{thm:stein} for the general structure of these bounds.
Note that $\delta_0$, $\delta_2$ and $\delta_3$ are used in both bounds, but $\delta_1$ is only used
for the Wasserstein distance bound and $\delta_1'$ is only used for the Kolmogorov distance bound.

To present the error terms, let us first define random variables
\begin{align}
    \label{eq:xbrackdef}
    x^{[i]} &\coloneqq (x_1, \dotsc, x_{i-1}, x_i, Y_{i+1}, \dotsc, Y_N), \\
    \label{eq:xparendef}
    x^{(i)} &\coloneqq (x_1, \dotsc, x_{i-1}, Y_i, x_{i+1}, \dotsc, x_N)
\end{align}
for any fixed $x$ and any index $i$.
Next, define
\begin{align}
    \label{eq:Delta1def}
    \Delta_{1,i}(x) &\coloneqq \fr{1}{2} \E \left[
        \left(f(x) - f(x^{(i)})\right)
        \left(f(x^{[i]}) - f(x^{[i-1]})\right) \right] \qquad \text{and} \\
    \label{eq:Delta2def}
    \Delta_{2,i}(x) &\coloneqq \fr{1}{2} \E \left[
        \left(\sg(x) - \sg(x^{(i)})\right)
        \left(f(x^{[i]}) - f(x^{[i-1]})\right) \right],
\end{align}
where $\sg$ is the exponential tilting function,
and set
\begin{equation}
\label{eq:bdef}
    b \coloneqq \E \left[
        \sum_{i=1}^N \Delta_{1,i}(X)
    \right].
\end{equation}
Finally, let $D_i^*(x,y)$ be any symmetric function such that
\begin{equation}
\label{eq:dstarcond}
    D_i^*(x,Y) \geq |f(x^{[i]}) - f(x^{[i-1]})|
\end{equation}
almost surely.

With these preliminary definitions in hand, we can write down the formulas for
the aforementioned error terms $\delta_0, \delta_1, \delta_1', \delta_2, \delta_3$,
as follows, and then state the theorem.
Note that our notation differs from \cite{fang2024normal} by extracting a common
term from their definitions of $\delta_1$ and $\delta_1'$, which we denote by $\delta_0$.
\begin{align}
    \label{eq:delta0def}
    \delta_0 &\coloneqq \sum_{i=1}^N \E \bigg[
        e^{\left|\sg(X) - \sg(X^{(i)})\right|}
        \left( \sg(X) - \sg(X^{(i)}) \right)^2 \\
    &\qquad \qquad \qquad \qquad \times \left(\left|\sg(X) - \sg(X^{(i)})\right|
        + \left|f(X) - f(X^{(i)})\right|\right)
        \left|f(X^{[i]}) - f(X^{[i-1]})\right|
    \bigg], \\
    \label{eq:delta1def}
    \delta_1 &\coloneqq \sum_{i=1}^N \E \left[
        \left( f(X) - f(X^{(i)}) \right)^2
        \left| f(X^{[i]}) - f(X^{[i-1]}) \right|
    \right], \\
    \label{eq:delta1primedef}
    \delta_1' &\coloneqq \sum_{i=1}^N \E \left[
        e^{\left|\sg(X) - \sg(X^{(i)})\right|}
        D_i^*(X,Y)
        \left|f(X) - f(X^{(i)})\right|
        \left|\sg(X) - \sg(X^{(i)})\right|
    \right] \\
    &\quad \qquad +\E\left[
        \left| \sum_{i=1}^N \Econd{D_i^*(X,Y) (f(X) - f(X^{(i)}))}{X} \right|
    \right], \\
    \label{eq:delta2def}
    \delta_2 &\coloneqq \sqrt{
        \Var \left[
            \sum_{i=1}^N \Delta_{1,i}(X)
        \right]},
    \\
    \label{eq:delta3def}
    \delta_3 &\coloneqq \sqrt{
        \Var \left[
            \sum_{i=1}^N \Delta_{2,i}(X)
            - (1-b) f(X)
        \right]}.
\end{align}
Recall that in the definition of $\delta_1'$ the variables $X$ and $Y$ are independent.
With these definitions in hand, we may state the theorem.

\begin{theorem}[Theorem 2.1 of \cite{fang2024normal}]
\label{thm:stein}
Let $X$ be as in \eqref{eq:tilt} and $f : \R^N \to \R$ be some function satisfying $\E[f(X)] = 0$.
Define $b, \delta_0, \delta_1, \delta_1', \delta_2, \delta_3$ as above, and let $Z \sim \cN(0,1)$.
Then we have
\begin{align}
    \dW(f(X),Z) &\leq \fr{C}{|b|} \left(\delta_0 + \delta_1 + \delta_2 + \delta_3\right), \\
    \text{and} \qquad \dK(f(X),Z) &\leq \fr{C}{|b|} \left(\delta_0 + \delta_1' + \delta_2 + \delta_3\right),
\end{align}
for some absolute constant $C$.
\end{theorem}

As mentioned multiple times already, this theorem is proved using \emph{Stein's method}, which is a
powerful technique that allows one to bound distributional distances of the following form, written in
terms of random variables $S$ and $Z$:
\begin{equation}
\label{eq:dcdef}
    \d_\cC(S,Z) \coloneqq \sup \left\{ \left| \E[\psi(S)] - \E[\psi(Z)] \right| : \psi \in \cC \right\};
\end{equation}
here $\cC$ is some class of functions $\psi : \R \to \R$.
For instance, one recovers the Wasserstein distance $\dW$ defined in \eqref{eq:wasdist} by choosing
\begin{equation}
    \cC = \cC_{\mathrm{Was}} \coloneqq \{ \psi : \R \to \R \text{ $1$-Lipschitz} \},
\end{equation}
and one recovers the Kolmogorov distance $\dK$ defined in \eqref{eq:koldist} by choosing
\begin{equation}
    \cC = \cC_{\mathrm{Kol}} \coloneqq \{ \ind*{(-\infty,s]} : s \in \R \}.
\end{equation}
In Stein's method one fixes the distribution of $Z$,
and for any $\psi \in \cC$ one considers the \emph{Stein equation} with unknown function $F$:
\begin{equation}
\label{eq:steineq}
    \cA F (s) = \psi(s) - \E[\psi(Z)].
\end{equation}
Here, $\cA$ is the \emph{characterizing operator} for the distribution $Z$.
In our case, $Z$ is a standard normal random variable and $\cA F (s) = s F(s) - F'(s)$.
Note that with this definition,
\begin{equation}
\label{eq:expzero}
    \E\left[\cA F (Z)\right] = 0
\end{equation}
for all nice enough functions $F$ by Gaussian integration-by-parts, and indeed the condition
\eqref{eq:expzero} is how one identifies the characterizing operator for other distributions in general.
In any case, letting $F_\psi$ denote the solution to \eqref{eq:steineq}, we immediately observe that
\begin{equation}
    \left| \E [ \psi(S) ] - \E[ \psi(Z) ] \right| = \left| \E[\cA F_\psi (S)]\right|,
\end{equation}
meaning that, in our context with $Z \sim \cN(0,1)$, recalling the definition \eqref{eq:dcdef}, we have
\begin{equation}
    \d_\cC(S,Z) = \sup \left\{ \left| \E \left[ S F_\psi(S) - F_\psi'(S) \right] \right| : \psi \in \cC \right\}.
\end{equation}
This is useful because for various classes $\cC$, the functions $F_\psi$ for $\psi \in \cC$ satisfy convenient
properties allowing for the above quantity to be bounded.
For instance, whenever $\psi \in \cC_{\mathrm{Kol}}$, we have $\| F_\psi \|_\infty, \| F_\psi' \|_\infty \leq 1$,
and whenever $\psi \in \cC_{\mathrm{Was}}$, we have $\| F_\psi \|_\infty, \| F_\psi' \|_\infty, \| F_\psi'' \|_\infty \leq 2$
(see e.g.\ \cite[Lemma 2.5]{ross2011fundamentals}).
Thus, in order to bound the Kolmogorov distance or Wasserstein distance between the random variable $f(X)$
and the Gaussian $Z$, it suffices to bound
\begin{equation}
\label{eq:tobound}
    \E \left[ f(X) F(f(X)) - F'(f(X)) \right]
\end{equation}
whenever $F$ satisfies either $\| F \|_\infty, \| F' \|_\infty \leq 1$ or $\| F \|_\infty, \| F' \|_\infty, \| F'' \|_\infty \leq 2$,
respectively.
For more information on Stein's method, we refer the reader to the survey \cite{ross2011fundamentals}.

Theorem \ref{thm:stein} follows this approach and bounds the expression \eqref{eq:tobound} in both cases
mentioned above, in a rather long and technical derivation making use of Taylor's theorem repeatedly,
leading to the expressions for the error terms $\delta_0, \delta_1, \delta_1', \delta_2, \delta_3$ appearing
in \eqref{eq:delta0def}-\eqref{eq:delta3def}, as well as the quantity $b$.
Note that typically, and in particular in our applications, $b$ will be a quantity of order $1$,
so using Theorem \ref{thm:stein} to prove a CLT reduces to bounding the five error terms appropriately.
\subsection{Modified Hamiltonian}
\label{sec:stein_hamiltonian}

In the setup for Theorem \ref{thm:stein} above, we would like to express the distribution of $X$ as
an exponential tilt of a product measure, i.e.\ the distribution of an independent random vector $Y$.
The method seems to work best if $X$ and $Y$ behave similarly on a macroscopic scale.
In our context, where $X$ is a sample from the conditioned measure $\condmeas$ for some $p \in U_\beta$,
the product measure should be $\cG(n,p)$.
Recall that every $p \in U_\beta$ is an attracting fixed point of the function
\begin{equation}
    \phi_\beta(q) = \fr{e^{2\sH'(q)}}{1 + e^{2\sH'(q)}},
\end{equation}
recalling the definition \eqref{eq:phidef1} of $\phi_\beta$,
as well as the discussion of Section \ref{sec:review_fluctuations_dynamical}.
Since $\log(\fr{z}{1-z})$ is the inverse of the function $\fr{e^z}{1+e^z}$, and since $\phi_\beta(p) = p$
for $p \in U_\beta$, we thus have
\begin{equation}
    2 \sH'(p) = \log \left(\fr{p}{1-p}\right),
\end{equation}
meaning that, for any $n$-vertex graph $x$, we have
\begin{equation}
\label{eq:tiltfromgnp}
    e^{n^2 \sH(x)} = e^{n^2 \sH(x) - 2 \sH'(p) \sE(x)} \cdot \left(\fr{p}{1-p}\right)^{\sE(x)}.
\end{equation}
Since $\P[Y = x] \propto (\fr{p}{1-p})^{\sE(x)}$ if $Y \sim \cG(n,p)$, \eqref{eq:tiltfromgnp} means that we may
express the full ERGM distribution as an exponential tilt of $\cG(n,p)$ by the function
\begin{equation}
    n^2 \sH(x) - 2 \sH'(p) \sE(x).
\end{equation}
However, recall that we actually wish to consider the conditioned measure $\condmeas$ for some $\eta > 0$.
One way to do this would be to set the tilting function to $-\infty$ outside of $\pball$.
However, since the error estimates in Theorem \ref{thm:stein} rely on estimates of the tilting function
on graphs which may be outside of $\pball$, this won't quite work for our purposes.

Instead, we will make a more gradual modification to the tilting function which will lead to a measure
very close to $\condmeas$.
Specifically, let us define
\begin{equation}
\label{eq:Rdef}
    \sR(x) \coloneqq n^2 \db(x, \pball),
\end{equation}
where $\db(x,\pball)$ denotes the infimal cut distance from the graphon representing $x$ to a graphon
in $\pball$; in particular, if $x \in \pball$, then $\sR(x) = 0$.
Then we will define
\begin{align}
\label{eq:smallgdef}
    \sg(x) &\coloneqq n^2 \sH(x) - 2 \sH'(p) \sE(x) - \sR(x) \\
\label{eq:smallgexpanded}
    &= \sum_{j=1}^K \beta_j \left( \fr{\sN_G(x)}{n^{\sv_j-2}} - 2 \se_j p^{\se_j-1} \sE(x) \right)
    - \sR(x),
\end{align}
noting that the $j=0$ term in the sum above is zero since $\sN_{G_0}(x) = 2 \sE(x)$.
Now we let $X$ be a random graph distributed according to \eqref{eq:tilt} with $Y \sim \cG(n,p)$, i.e.
\begin{equation}
    \P[X = x] \propto e^{\sg(x)} \P[Y = x] \propto \Exp{n^2 (\sH(x) - \db(x,\pball))}.
\end{equation}
In other words, we have a general ERGM \eqref{eq:gibbsmeasure} with Hamiltonian
$W \mapsto \sH(W) - \db(W, \pball)$, which is a continuous function on the space of graphons.
Thus the general large deviations principle, Theorem \ref{thm:LDP_general} applies, and samples
from this ERGM are typically near a maximizer of $F(W) = \sH(W) - I(W) - \db(W,\pball)$, where $I$ is the
entropy functional of \eqref{eq:entropy}.

We now claim that $F$ has a \emph{unique} maximizer, which is $W_p$.
To prove this, first note that since $\sH(W)$ has the form \eqref{eq:hdef},
namely a linear combination of homomorphism densities
with ferromagnetic parameters $\beta$, Theorem \ref{thm:LDP_ferromagnetic} applies and the maximizers
of the function $G(W) = \sH(W) - I(W)$ (without the distance term)
are all constant graphons $W_q$ with $q \in M_\beta$, which is a finite set.
Thus if $\eta$ is small enough, since $\db(W,\pball)$ is zero for $W \in \pball$ and is positive for
$W$ outside of $\pball$, we have $F(W) \leq G(W)$ for all graphons $W$,
and we have $F(W_p) = G(W_p)$ and $F(W_q) < G(W_q)$ for $q \in M_\beta \setminus \{p\}$.
This proves the claim.

So Theorem \ref{thm:LDP_general} states that $X$ is within $\pball$ with probability $1 - \Exp{-\Omega(n^2)}$.
Also since $\db(x,\pball) = 0$ for $x \in \pball$, the distribution of $X$ conditioned on being within
$\pball$ is \emph{exactly the same} as $\condmeas$.
All together, the distribution of $X$ given by \eqref{eq:tilt} for this choice of $\sg$
has total variation distance $\Exp{-\Omega(n^2)}$ to $\condmeas$.

\begin{remark}
\label{rmk:equivalence}
This means that to prove Theorems \ref{thm:global_clt_intro} and \ref{thm:local_clt_intro},
it suffices to prove the corresponding statements for this $X$.
Moreover, the H\'ajek projection results
Propositions \ref{prop:global_hajek_intro} and \ref{prop:local_hajek_intro} hold for this $X$,
as do the Lipschitz concentration inequality Theorem \ref{thm:lipschitzconcentration} and its consequences,
Propositions \ref{prop:Eproduct} and \ref{prop:marginal}, as well as our ``basic'' Lipschitz
concentration inequality, Lemma \ref{lem:basic_concentration}.
\end{remark}

Before we continue, let us record the following useful facts about $\partial_e \sR$ (recalling this notation
from \eqref{eq:partialdef}) which may easily be verified  from the definition \eqref{eq:dbdef} of the cut
distance and the fact that $X$ is within $\phalfball$ with probability $1 - \exp(-\Omega(n^2))$.

\begin{lemma}
\label{lem:Rfacts}
There are constants $c, C > 0$ such that for all $e \in \edgeset$,
we have $|\partial_e \sR(x)| \leq C$ for all $n$-vertex graphs $x$, and
moreover $\P[\partial_e \sR(X) \neq 0] \leq e^{-c n^2}$.
\end{lemma}
\subsection{Quantitative CLT for edge count}
\label{sec:stein_global}

In this section we prove Theorem \ref{thm:global_clt_intro} which has been restated below for the
reader's convenience.
This is a quantitative central limit theorem for the edge count $\sE(X)$ of the $n$-vertex random graph $X$.

As mentioned in the previous section, for this theorem we will let $X$ have the distribution
\eqref{eq:tilt}, with $Y \sim \cG(n,p)$ for some $p \in U_\beta$,
and $\sg$ given by \eqref{eq:smallgdef}; as also mentioned in the previous
section, this is within $\Exp{-\Omega(n^2)}$ total-variation distance from $\condmeas$
(recalling the definition of this conditioned measure from \eqref{eq:condmeas}),
so Theorem \ref{thm:global_clt_intro} is equivalent to our statement below (which will be reiterated
for the reader's convenience).
Additionally, we make one of the following two assumptions: either $M_\beta = U_\beta = \{p\}$
(i.e.\ we are in the non-critical phase uniqueness case), or \emph{all graphs} $G_0, G_1, \dotsc, G_K$
in the ERGM specification \eqref{eq:hdef} are \emph{forests}.
As was the case for Propositions \ref{prop:global_hajek_intro} and \ref{prop:local_hajek_intro},
the only place we use this ``phase-uniqueness or forest'' assumption is in the application of
Proposition \ref{prop:Eproduct}, which is also implicit in
the use of Proposition \ref{prop:marginal} and Lemma \ref{lem:global_step_concentration}.

Recall that we denote by $G_j = (\cV_j, \cE_j)$ for $j = 0, \dotsc, K$ the graphs in the ERGM specification
\eqref{eq:hdef}, and moreover that we set $\sv_j = |\cV_j|$ and $\se_j = |\cE_j|$.
Since $G_0$ is always a single edge, the $j=0$ term drops out of most sums we consider.

\begin{theorem}[Restatement of Theorem \ref{thm:global_clt_intro}]
\label{thm:global_clt}
Let us make Assumption \ref{assumption}.
As in \eqref{eq:sigmadef_global_intro}, set
\begin{equation}
\label{eq:sigmadef_global}
    \sigma_n^2 = 
    \left( 1 - 2 p (1-p) \sum_{j=1}^K \beta_j \se_j (\se_j-1) p^{\se_j-2} \right)^{-1}
    \times p(1-p) \binom{n}{2},
\end{equation}
Then for any $\eps > 0$, if $X$ is distributed as \eqref{eq:tilt} with $Y \sim \cG(n,p)$ and $\sg$ given in
\eqref{eq:smallgdef}, under the ``phase-uniqueness or forest'' assumption \ref{assumption}, we have
\begin{equation}
    \d \left( \fr{\sE(X) - \E[\sE(X)]}{\sigma_n}, Z \right) \lesssim n^{-\fr{1}{2} + \eps},
\end{equation}
for both $\d = \dW$ and $\d = \dK$.
Recalling Remark \ref{rmk:equivalence}, this is equivalent to the original statement of
Theorem \ref{thm:global_clt_intro} which was stated for $X \sim \condmeas$.
\end{theorem}

The reason for this explicit form of $\sigma_n^2$ will become clear later, in Section \ref{sec:stein_global_delta3}.
Until then, we will only make use of the fact that $\sigma_n = \Theta(n)$.
However, before proceeding, we should check that $\sigma_n^2$ is well-defined, i.e.\ that
\begin{equation}
\label{eq:sigmawelldef}
    2 p (1-p) \sum_{j=1}^K \beta_j \se_j (\se_j - 1) p^{\se_j-2} < 1.
\end{equation}
To see this, note that $p$ is an attracting fixed point of $\phi_\beta$, i.e.\ $\phi_\beta(p) = p$
and $\phi_\beta'(p) < 1$.
So we have
\begin{align}
    1 &> \phi_\beta'(p) \\
    &= \phi_\beta(p) (1 - \phi_\beta(p)) \cdot 2 \sH''(p) \\
    &= p (1-p) \cdot 2 \sum_{j=0}^K \beta_j \se_j (\se_j-1) p^{\se_j-2},
\end{align}
recalling from \eqref{eq:phidef1} that $\phi_\beta(q) = \phi(2 \sH'(q))$ with
$\phi(z) = \fr{e^z}{1+e^z}$, which satisfies $\phi'(z) = \phi(z) (1 - \phi(z))$.

To prove Theorem \ref{thm:global_clt_intro}/\ref{thm:global_clt}, we will set
\begin{equation}
    f(X) = \fr{\sE(X) - \E[\sE(X)]}{\sigma_n}
\end{equation}
and apply Theorem \ref{thm:stein} to this choice of $f$.
In other words, we need to bound the expression $b$ given in \eqref{eq:bdef} from below and bound
the error terms $\delta_0, \delta_1, \delta_1', \delta_2, \delta_3$ given by
\eqref{eq:delta0def}-\eqref{eq:delta3def} from above by appropriate quantities.

\subsubsection{Preliminary calculations}
\label{sec:stein_global_prelim}

We begin by calculating expressions for $b$ given by \eqref{eq:bdef}, as well as
the quantities $\Delta_{1,e}(x)$ and $\Delta_{2,e}(x)$ given by \eqref{eq:Delta1def} and \eqref{eq:Delta2def},
for a fixed graph $x$ and edge $e \in \edgeset$.
Recall also the definitions of $x^{[e]}$ and $x^{(e)}$ from \eqref{eq:xbrackdef} and \eqref{eq:xparendef},
and notice that $x^{[e]}$ depends on some \emph{ordering} of the edges $e \in \edgeset$.
For our present purposes, this ordering will not be relevant, so we fix an arbitrary ordering,
and denote by $e-1$ the edge preceding $e$ in this ordering.
Moreover, note that since $\sE(x) = \sum_{e \in \edgeset} x(e)$, for every $e \in \edgeset$
we have
\begin{equation}
\label{eq:differences}
    f(x) - f(x^{(e)}) = f(x^{[e]}) - f(x^{[e-1]}) = \fr{x(e) - Y(e)}{\sigma_n}.
\end{equation}
The above equality will simplify the formulas for the error terms greatly, and we will henceforth
be able to mostly ignore the ordering on the edges.

First, let us calculate $\Delta_{1,e}(x)$. 
Recalling the definition from \eqref{eq:Delta1def} and applying \eqref{eq:differences}, we find that
\begin{align}
    \Delta_{1,e}(x) &= \fr{1}{2\sigma_n^2} \E[(x(e) - Y(e))^2] \\
    &= \fr{1}{2\sigma_n^2} ((1-2p) x(e) + p),
    \label{eq:squarecalculation}
\end{align}
since $Y \sim \cG(n,p)$ and $x(e)^2 = x(e)$.
Recalling the definition of $b$ from \eqref{eq:bdef}, this yields
\begin{align}
\label{eq:bcalc}
    b &= \fr{1}{2\sigma_n^2} \sum_{e \in \edgeset} \E[(1-2p) X(e) + p] \\
    &= \fr{\binom{n}{2} p (1-p)}{\sigma_n^2} + O\left(n^{-\fr{1}{2}} \sqrt{\log n}\right)
\end{align}
by Proposition \ref{prop:marginal} which approximates $\E[X(e)]$, also using the fact that $\sigma_n = \Theta(n)$.
Since $b$ is of constant order, Theorem \ref{thm:global_clt_intro}/\ref{thm:global_clt}
will follow if we can bound each $\delta_i$ appropriately.

Next let us calculate $\Delta_{2,e}(x)$, recalling its definition from \eqref{eq:Delta2def}.
Recalling also the expanded definition of the tilting function $\sg$ in \eqref{eq:smallgexpanded}
and also using \eqref{eq:differences}, we have
\begin{align}
    \Delta_{2,e}(x) &= \sum_{j=1}^K \fr{\beta_j}{2 n^{\sv_j-2}}
    \E \left[
        \left(
            \sN_{G_j}(x) - \sN_{G_j}(x^{(e)})
            - 2 \se_j p^{\se_j-1} n^{\sv_j-2}
            (x(e) - Y(e)) \right)
            \fr{x(e) - Y(e)}{\sigma_n}
    \right] \\
    &\quad - \fr{1}{\sigma_n} \E \left[(\sR(x) - \sR(x^{(e)}))(x(e) - Y(e))\right].
\end{align}
Now note that
\begin{equation}
    \sN_{G_j}(x) - \sN_{G_j}(x^{(e)}) = (x(e) - Y(e)) \cdot \partial_e \sN_{G_j}(x)
    = (x(e) - Y(e)) \sN_{G_j}(x,e),
\end{equation}
where $\sN_G(x,e)$ is the number of homomorphisms of $G$ in $x^{+e}$ using the edge $e$.
So by the same calculation as in \eqref{eq:squarecalculation}, we find that
\begin{align}
\label{eq:Delta2calcglob}
    \Delta_{2,e}(x) &= \sum_{j=1}^K \fr{\beta_j}{2 n^{\sv_j-2} \sigma_n} ((1-2p) x(e) + p)
    \left( \sN_{G_j}(x,e) - 2 \se_j p^{\se_j-1} n^{\sv_j-2} \right) \\
    &\quad - \fr{1}{\sigma_n} \E \left[(\sR(x) - \sR(x^{(e)}))(x(e) - Y(e))\right]
\end{align}
When we plug in $x = X$ later, the latter term will be $e^{-\Omega(n^2)}$ in probability,
by Lemma \ref{lem:Rfacts}, since we have $|\sR(x) - \sR(x^{(e)}) | \leq \partial_e \sR(x)$.
However, we carry this term through for now since much of the work will be done in deterministic
equalities before we plug in $x = X$.

With the above preparation finished, we may begin bounding the error terms one-by-one.
Most are relatively straightforward, and similar calculations appeared already in \cite{fang2024normal}.
However, we present these details in a more streamlined way which will also serve as a warm-up for the
vertex degree CLT, which is novel and a bit more algebraically technical.

\subsubsection{Bounding $\delta_0$}
\label{sec:stein_global_delta0}

Note first that by \eqref{eq:smallgexpanded}, we have
\begin{align}
    \left| \sg(x) - \sg(x^{(e)}) \right|
    &\leq |\partial_e \sg(x)| \\
    &\leq \sum_{j=1}^K \beta_j \left( \fr{\sN_{G_j}(x,e)}{n^{\sv_j-2}} + 2 \se_j p^{\se_j-1} \right) 
    + |\partial_e \sR(x)| \\
    &\lesssim 1 + |\partial_e \sR(x)|,
\end{align}
since $\sN_{G_j}(x,e) \lesssim n^{\sv_j-2}$ uniformly in $x$.
The above expression is $O(1)$ uniformly in $x$ by Lemma \ref{lem:Rfacts}.
Additionally, since $\sigma_n = \Theta(n)$, we have $|f(x) - f(x')| \lesssim \fr{1}{n}$
whenever $x,x'$ differ by at most one edge.
Therefore, recalling the definition of $\delta_0$ from \eqref{eq:delta0def}, and using
the fact that $e^{|\sg(x) - \sg(x^{(e)})|} = O(1)$ uniformly in $x$ from the above calculation,
we have
\begin{align}
    \delta_0 &= \sum_{e \in \edgeset} \E \bigg[
        e^{\left|\sg(X) - \sg(X^{(e)})\right|}
        \left( \sg(X) - \sg(X^{(e)}) \right)^2 \\
    &\qquad \qquad \qquad \qquad \times \left(\left|\sg(X) - \sg(X^{(e)})\right|
        + \left|f(X) - f(X^{(e)})\right|\right)
        \left|f(X^{[e]}) - f(X^{[e-1]})\right|
    \bigg] \\
    &\lesssim \fr{1}{n} \sum_{e \in \edgeset}
    \E \left[
        \left|\partial_e \sg(X)\right|^2
        \left( \left|\partial_e \sg(X)\right| + \fr{1}{n} \right)
    \right]. \\
    \label{eq:del0globeq1}
    &\lesssim \fr{1}{n} \sum_{e \in \edgeset}
        \E \left[ \left| n^2 \partial_e \sH(X) - 2 \sH'(p) \right|^3 \right]
    +
    \fr{1}{n^2} \sum_{e \in \edgeset}
        \E \left[ \left| n^2 \partial_e \sH(X) - 2 \sH'(p) \right|^2 \right]
    + e^{-\Omega(n^2)},
\end{align}
where we applied formula \eqref{eq:smallgdef} for $\sg$ and Lemma \ref{lem:Rfacts} at the last step.
Now, by Lemma \ref{lem:global_step_concentration} and a union bound, since 
\begin{equation}
    n^2 \partial_e \sH(X) = \sum_{j=0}^K \fr{\beta_j \sN_{G_j}(X,e)}{n^{\sv_j-2}}
    \qquad \text{and} \qquad
    2 \sH'(p) = 2 \sum_{j=0}^K \beta_j \se_j p^{\se_j-1},
\end{equation}
for any small enough $\alpha > 0$ there is some $c > 0$ for which
\begin{equation}
\label{eq:hconc}
    \P\left[\left| n^2 \partial_e \sH(X) - 2 \sH'(p) \right| > n^{-\fr{1}{2}+\alpha} \right]
    \leq \Exp{- c n^{2\alpha}}.
\end{equation}
This implies that for any $\eps > 0$ we have
\begin{equation}
    \E \left[\left|n^2 \partial_e \sH(X) - 2 \sH'(p) \right|^2 \right] \lesssim n^{-1 + \eps},
    \qquad \text{and} \qquad
    \E \left[\left|n^2 \partial_e \sH(X) - 2 \sH'(p) \right|^3 \right] \lesssim n^{-\fr{3}{2} + \eps}.
\end{equation}
Plugging these into \eqref{eq:del0globeq1},
we obtain $\delta_0 \lesssim n^{-\fr{1}{2} + \eps}$ for any $\eps > 0$.

We remark that \eqref{eq:hconc} is essentially the only thing preventing us from achieving a CLT
in the phase coexistence case when the graphs $G_0, G_1, \dotsc, G_K$ are not necessarily forests,
since it ultimately relies on Proposition \ref{prop:Eproduct}.
If \eqref{eq:hconc} could be proved without this assumption, then the proof presented here could be 
adapted to the more general case.
This would follow by using the version of Proposition \ref{prop:global_hajek_intro} with a non-explicit
H\'ajek projection factor, and so in that case the variance proxy $\sigma_n^2$ would also be non-explicit.

\subsubsection{Bounding $\delta_1$}
\label{sec:stein_global_delta1}

Using again the fact that $|f(x) - f(x')| \lesssim \fr{1}{n}$ when $x,x'$ differ at one edge,
and recalling the definition of $\delta_1$ from \eqref{eq:delta1def}, we have
\begin{equation}
    \delta_1 = \sum_{e \in \edgeset} \E \left[
        \left( f(X) - f(X^{(e)}) \right)^2
        \left| f(X^{[e]}) - f(X^{[e-1]}) \right|
    \right]
    \lesssim \sum_{e \in \edgeset} \fr{1}{n^3} \lesssim \fr{1}{n}.
\end{equation}

\subsubsection{Bounding $\delta_1'$}
\label{sec:stein_global_delta1prime}

We may take $D_e^*(x,y) = \fr{1}{\sigma_n}$, since the edge indicators are all
bounded by $1$ and so \eqref{eq:dstarcond} holds with this choice.
For the first term of $\delta_1'$ (recalling the definition from \eqref{eq:delta1primedef}),
similar to the bound for $\delta_0$, we have
\begin{align}
\label{eq:globaldelta1primea}
    &\sum_{e \in \edgeset} \E \left[
        e^{\left|\sg(X) - \sg(X^{(e)})\right|}
        D_e^*(X,Y)
        \left|f(X) - f(X^{(e)})\right|
        \left|\sg(X) - \sg(X^{(e)})\right|
    \right] \\
   &\qquad \qquad \lesssim
   \fr{1}{n^2} \sum_{e \in \edgeset} \E\left[|n^2 \partial_e \sH(X) - 2 \sH'(p)|\right]
   + e^{-\Omega(n^2)},
\end{align}
and this is $\lesssim n^{-\fr{1}{2} + \eps}$ for any $\eps > 0$ by \eqref{eq:hconc}.
For the second term in $\delta_1'$, we have
\begin{equation}
\label{eq:globaldelta1primeb}
    \E\left[
        \left| \sum_{e \in \edgeset} \Econd{D_e^*(X,Y) (f(X) - f(X^{(e)}))}{X} \right|
    \right]
    = \fr{1}{\sigma_n^2} \E \left[ \left| \sum_{e \in \edgeset} (X(e) - p) \right| \right]
    \lesssim n^{-\fr{1}{2}} \sqrt{\log n},
\end{equation}
using Proposition \ref{prop:marginal} (recalling Remark \ref{rmk:equivalence}).
This means that $\delta_1' \lesssim n^{-\fr{1}{2} + \eps}$ for any $\eps > 0$.

\subsubsection{Bounding $\delta_2$}
\label{sec:stein_global_delta2}

Recall the definition of $\delta_2$ from \eqref{eq:delta2def}.
By the calculation \eqref{eq:squarecalculation} of $\Delta_{1,e}(x)$, we have
\begin{equation}
    \delta_2^2 = \Var \left[ \sum_{e \in \edgeset} \Delta_{1,e}(X) \right]
    = \fr{(1-2p)^2}{4 \sigma_n^4} \Var[ \sE(X) ].
\end{equation}
Now note that $\Var[\sE(X)] \lesssim n^2$ by \cite[Proposition 5.14]{winstein2025concentration}
(which is a consequence of the Lipschitz concentration Theorem \ref{thm:lipschitzconcentration}),
or by an application of Lemma \ref{lem:basic_concentration}, our ``basic'' Lipschitz concentration
inequality (recalling Remark \ref{rmk:equivalence}).
So, since $\sigma_n^4 = \Theta(n^4)$, we obtain $\delta_2 \lesssim n^{-1}$.

\subsubsection{Bounding $\delta_3$}
\label{sec:stein_global_delta3}

This final bound is the most technical, and it will be even more so in the proof of
Theorem \ref{thm:local_clt_intro}, the CLT for vertex degrees.
The reader should thus carefully follow these calculations (most of which have already
appeared in \cite{fang2024normal}) as a warm-up for that later proof.

To begin with, let us define
\begin{equation}
    \Delta_3(x) \coloneqq \sum_{e \in \edgeset} \Delta_{2,e}(x) - (1-b) f(x),
\end{equation}
so that $\delta_3 = \sqrt{\Var[\Delta_3(X)]}$, recalling the definition from \eqref{eq:delta3def}.
Then, using \eqref{eq:bcalc} and \eqref{eq:Delta2calcglob} for the expressions of $b$ and $\Delta_{2,e}(x)$
respectively, we have
\begin{align}
    \Delta_3(x) &= \sum_{j=1}^K \fr{\beta_j}{2 n^{\sv_j-2} \sigma_n}
    \sum_{e \in \edgeset} ((1-2p) x(e) + p) \left(\sN_{G_j}(x,e) - 2 \se_j p^{\se_j-1} n^{\sv_j-2} \right) \\
    &\qquad - \left(1 - \fr{\binom{n}{2} p(1-p)}{\sigma_n^2} + O(n^{-\fr{1}{2}} \sqrt{\log n}) \right) \fr{\sE(x) - \E[\sE(X)]}{\sigma_n} \\
    &\qquad - \fr{1}{\sigma_n} \E \left[ \sum_{e \in \edgeset} (\sR(x) - \sR(x^{(e)})) (x(e) - Y(e)) \right].
\end{align}
Rearranging this gives
\begin{align}
    \label{eq:delta30line1}
    \Delta_3(x) &= \sum_{j=1}^K \fr{\beta_j (1-2p)}{2 n^{\sv_j-2} \sigma_n}
    \left( \sum_{e \in \edgeset} x(e) \sN_{G_j}(x,e) - 2 \se_j p^{\se_j-1} n^{\sv_j-2} \sE(x) \right) \\
    \label{eq:delta30line2}
    &\qquad + \sum_{j=1}^K \fr{\beta_j p}{2 n^{\sv_j-2} \sigma_n}
    \left( \sum_{e \in \edgeset} \sN_{G_j}(x,e) - 2 \se_j p^{\se_j-1} n^{\sv_j-2} \binom{n}{2} \right) \\
    &\qquad - \left(1 - \fr{\binom{n}{2} p(1-p)}{\sigma_n^2} + O(n^{-1/2} \sqrt{\log n}) \right) \fr{\sE(x) - \E[\sE(X)]}{\sigma_n} \\
    &\qquad - \fr{1}{\sigma_n} \E \left[ \sum_{e \in \edgeset} (\sR(x) - \sR(x^{(e)})) (x(e) - Y(e)) \right].
\end{align}
Now in line \eqref{eq:delta30line1} we would like to replace
\begin{equation}
    \sum_{e \in \edgeset} x(e) \sN_{G_j}(x,e) \qquad \text{with} \qquad \se_j \sN_{G_j}(x),
\end{equation}
as in \cite{fang2024normal}; however, in that work they restricted to \emph{injective} homomorphisms,
under which the above two expressions would be equal since each homomorphism of $G_j$ in $x$ would use exactly $\se_j$
different edges in $\edgeset$.
However, we make no such restriction and as such we will need to carry the error term through the
rest of the proof and later bound its variance appropriately.
We do something similar in line \eqref{eq:delta30line2} to replace
\begin{equation}
    \sum_{e \in \edgeset} \sN_{G_j}(x,e) \qquad \text{with} \qquad
    \sum_{\phi \in \cE_j} \sN_{G_j \setminus \phi} (x).
\end{equation}
Again the above two expressions are equal if we only count injective homomorphisms, since
every injective homomorphism of $G_j$ in $x^{+e}$ mapping an edge $\phi \in \cE_j$ to $e \in \edgeset$.
restricts to an injective homomorphism of $G_j \setminus \phi$ in $x$.

All this is to say that we have
\begin{align}
    \label{eq:d3line1}
    \Delta_3(x)
    &= \sum_{j=1}^K \fr{\beta_j (1-2p) \se_j}{2 n^{\sv_j-2} \sigma_n}
    \left( \sN_{G_j}(x) - 2 p^{\se_j-1} n^{\sv_j-2} \sE(x) \right) \\
    \label{eq:d3line2}
    &\qquad + \sum_{j=1}^K \fr{\beta_j p}{2 n^{\sv_j-2} \sigma_n}
    \sum_{\phi \in \cE_j} \left( \sN_{G_j \setminus \phi}(x) - 2 p^{\se_j-1} n^{\sv_j-2} \binom{n}{2} \right) \\
    &\qquad - \left(1 - \fr{\binom{n}{2} p(1-p)}{\sigma_n^2} + O(n^{-\fr{1}{2}} \sqrt{\log n}) \right) \fr{\sE(x) - \E[\sE(X)]}{\sigma_n} \\
    \label{eq:d3line-3}
    &\qquad + \sum_{j=1}^K \fr{\beta_j (1-2p) \se_j}{2 n^{\sv_j-2} \sigma_n} \left(\sum_{e \in \edgeset} x(e) \sN_{G_j}(x,e) - \se_j \sN_{G_j}(x) \right) \\
    \label{eq:d3line-2}
    &\qquad + \sum_{j=1}^K \fr{\beta_j p}{2 n^{\sv_j-2} \sigma_n} \left( \sum_{e \in \edgeset} \sN_{G_j}(x,e) - \sum_{\phi \in \cE_j} \sN_{G_j \setminus \phi}(x) \right) \\
    \label{eq:d3line-1}
    &\qquad - \fr{1}{\sigma_n} \E \left[ \sum_{e \in \edgeset} (\sR(x) - \sR(x^{(e)})) (x(e) - Y(e)) \right].
\end{align}
Now we will replace the second terms in the parentheses on lines \eqref{eq:d3line1} and \eqref{eq:d3line2}
with the correct H\'ajek projection of the corresponding first terms, and we obtain
\begin{align}
    \Delta_3(x)
    &= \sum_{j=1}^K \fr{\beta_j (1-2p) \se_j}{2 n^{\sv_j-2} \sigma_n}
    \left( \sN_{G_j}(x) - 2 \se_j p^{\se_j-1} n^{\sv_j-2} \sE(x) \right) \\
    &\qquad + \sum_{j=1}^K \fr{\beta_j p}{2 n^{\sv_j-2} \sigma_n}
    \sum_{\phi \in \cE_j} \left( \sN_{G_j \setminus \phi}(x) - 2 (\se_j-1) p^{\se_j-2} n^{\sv_j-2} \sE(x) \right) \\
    &\qquad + \sum_{j=1}^K \beta_j \left( (1-2p) \se_j^2 p^{\se_j-1} - (1-2p) \se_j p^{\se_j-1} + p \se_j (\se_j - 1) p^{\se_j-2} \right) \fr{\sE(x)}{\sigma_n} \\
    &\qquad - \sum_{j=1}^K \fr{\beta_j \se_j p^{\se_j}}{\sigma_n} \binom{n}{2} \\
    &\qquad - \left(1 - \fr{\binom{n}{2} p(1-p)}{\sigma_n^2} + O(n^{-\fr{1}{2}} \sqrt{\log n}) \right) \fr{\sE(x) - \E[\sE(X)]}{\sigma_n} \\
    &\qquad + \mathsf{remainder}(x),
\end{align}
where $\mathsf{remainder}(x)$ denotes the sum of lines \eqref{eq:d3line-3}, \eqref{eq:d3line-2} and \eqref{eq:d3line-1}
from the previous display.
Now let us combine terms in a more convenient fashion and recall the definition of $\hsN_G$ from \eqref{eq:hatngdef} to get
\begin{align}
    \Delta_3(x)
    &= \sum_{j=1}^K \fr{\beta_j (1-2p) \se_j}{2 n^{\sv_j-2} \sigma_n} \hsN_{G_j}(x) \\
    &\qquad + \sum_{j=1}^K \fr{\beta_j p}{2 n^{\sv_j-2} \sigma_n} \sum_{\phi \in \cE_j} \hsN_{G_j \setminus \phi}(x) \\
    \label{eq:sigmareason}
    &\qquad + \left( 2 p (1-p) \sum_{j=1}^K \beta_j \se_j (\se_j-1) p^{\se_j-2} - 1 + \fr{\binom{n}{2} p (1-p)}{\sigma_n^2} + O(n^{-\fr{1}{2}}\sqrt{\log n}) \right) \fr{\sE(x)}{\sigma_n} \\
    \label{eq:deterministic}
    &\qquad - \sum_{j=1}^K \fr{\beta_j \se_j p^{\se_j}}{\sigma_n} \binom{n}{2}
    + \fr{\E[\sE(X)]}{\sigma_n} \left(1 - \fr{\binom{n}{2} p(1-p)}{\sigma_n^2} + O(n^{-\fr{1}{2}} \sqrt{\log n}) \right) \\
    &\qquad + \mathsf{remainder}(x).
\end{align}
Now, by the definition \eqref{eq:sigmadef_global} of $\sigma_n^2$, line \eqref{eq:sigmareason} is simply
$O\left(n^{-\fr{1}{2}} \sqrt{\log n}\right) \cdot \fr{\sE(x)}{\sigma_n} = O\left(n^{-\fr{3}{2}}\sqrt{\log n}\right) \cdot \sE(x)$;
in fact,  this is the reason for that definition.

Finally we are ready to begin bounding the variance of $\Delta_3(X)$.
First note that the line \eqref{eq:deterministic} is deterministic so we can ignore it.
Further, by Jensen's inequality we have
\begin{equation}
    \left( \sum_{i=1}^m a_i \right)^2 \leq m \sum_{i=1}^m a_i^2
\end{equation}
for arbitrary real numbers $a_1, \dotsc, a_m$, and so we obtain
\begin{align}
    \label{eq:d32line1}
    \fr{\Var[\Delta_3(X)]}{K + \sum_{j=1}^K \se_j + 2}
    &\leq \sum_{j=1}^K \fr{\beta_j^2 (1-2p)^2}{4 n^{2\sv_j-4} \sigma_n^2} \Var\left[\hsN_{G_j}(X)\right] \\
    \label{eq:d32line2}
    &\qquad + \sum_{j=1}^K \fr{\beta_j^2 p^2}{4 n^{2\sv_j-4} \sigma_n^2} \sum_{\phi \in \cE_j} \Var\left[\hsN_{G_j \setminus \phi}(X)\right] \\
    \label{eq:d32line3}
    &\qquad + O(n^{-3} \log n) \cdot \Var[\sE(X)] \\
    &\qquad + \Var[\mathsf{remainder}(X)].
\end{align}
Now, by Proposition \ref{prop:global_hajek_intro} (recalling Remark \ref{rmk:equivalence}),
for any $\eps > 0$ and any graph $G$ with $\sv$ vertices
(which must be a forest in the phase coexistence case by Assumption \ref{assumption}) we have
\begin{equation}
    \Var\left[\hsN_G(X)\right] \lesssim n^{2 \sv - 3 + \eps}.
\end{equation}
Since we assumed that all $G_0, G_1, \dotsc, G_K$ are forests in the phase coexistence case,
all graphs $G_j \setminus \phi$ are also forests in this case, so the above inequality holds for all relevant graphs $G$.
Let us plug this into lines \eqref{eq:d32line1} and \eqref{eq:d32line2}
and use the fact that $\Var[\sE(X)] \lesssim n^2$ (as mentioned in Section \ref{sec:stein_global_delta2})
in line \eqref{eq:d32line3}.
Using also the fact that $\sigma_n^2 = \Theta(n^2)$, we find that
\begin{equation}
\label{eq:d3varbound1}
    \Var[\Delta_3(X)] \lesssim n^{-1+\eps}
    + n^{-1} \log n + \Var[\mathsf{remainder}(X)],
\end{equation}
and it only remains to bound the variance of $\mathsf{remainder}(X)$.

Recall that $\mathsf{remainder}(x)$ is the sum of lines \eqref{eq:d3line-3}, \eqref{eq:d3line-2},
and \eqref{eq:d3line-1}, so we have
\begin{align}
    \label{eq:d33line1}
    \fr{\Var[\mathsf{remainder}(X)]}{2K+1}
    &\leq \sum_{j=1}^K \fr{\beta_j^2 (1-2p)^2 \se_j^2}{4 n^{2\sv_j-4} \sigma_n^2}
    \Var\left[ \sum_{e \in \edgeset} X(e) \sN_{G_j}(X,e) - \se_j \sN_{G_j}(X) \right] \\
    \label{eq:d33line2}
    &\qquad + \sum_{j=1}^K \fr{\beta_j^2 p^2}{4 n^{2\sv_j-4} \sigma_n^2}
    \Var\left[ \sum_{e \in \edgeset} \sN_{G_j}(X,e) - \sum_{\phi \in \cE_j} \sN_{G_j \setminus \phi}(X) \right] \\
    \label{eq:d33line3}
    &\qquad + \fr{1}{\sigma_n^2} \Var \left[ \Econd{\sum_{e \in \edgeset} (\sR(X) - \sR(X^{(e)})) (X(e) - Y(e))}{X} \right].
\end{align}
First, by Lemma \ref{lem:Rfacts}, line \eqref{eq:d33line3} is $\Exp{-\Omega(n^2)}$, so it remains to consider
the other two lines \eqref{eq:d33line1} and \eqref{eq:d33line2}, for which we will apply the ``basic''
Lipschitz concentration inequality, Lemma \ref{lem:basic_concentration} (again recalling Remark \ref{rmk:equivalence}).

Let us first consider \eqref{eq:d33line1}; we would like to construct a Lipschitz vector for the function
\begin{equation}
    h_G(x) \coloneqq \sum_{e \in \edgeset} x(e) \sN_{G}(x,e) - \se \sN_{G}(x),
\end{equation}
for some graph $G$ with $\sv$ vertices and $\se$ edges.
First note as discussed above that if we were only counting \emph{injective} homomorphisms, the above function
would be zero.
So the difference consists of a sum of terms, each of order $1$, for each \emph{non-injective} homomorphism.
Thus for each $e \in \edgeset$ uniformly,
\begin{equation}
    \left| \partial_e h_G(x) \right| \lesssim n^{\sv - 3},
\end{equation}
since there are $\lesssim n^{\sv-3}$ non-injective homomorphisms of $G$ in $x$ which have the edge $e$
in their image.
This means that the Lipschitz vector $\cL_G$ for $h_G$ satisfies $\| \cL_G \|_\infty \lesssim n^{\sv-3}$,
and so the result of Lemma \ref{lem:basic_concentration} is that
\begin{equation}
    \P \left[ |h_G(X) - \E[h_G(X)]| > \lambda \right] \leq \Exp{- c\fr{\lambda^2}{n^{2\sv-6} n^2}}
\end{equation}
for some $c > 0$, whenever $\lambda \geq e^{-\Omega(n)}$ and $\lambda \lesssim n^{\sv-3+\fr{3}{2}}$.
In particular,
\begin{equation}
\label{eq:noninjvar1}
    \Var[h_G(X)] \lesssim n^{2\sv-4},
\end{equation}
using also the fact that $h_G(x)$ is uniformly polynomially bounded.

Next we turn to \eqref{eq:d33line2}, i.e.\ we must construct a Lipschitz vector for the function
\begin{equation}
    h_G'(x) \coloneqq \sum_{e \in \edgeset} \sN_{G_j}(x,e) - \sum_{\phi \in \cE_j} \sN_{G_j \setminus \phi}(x).
\end{equation}
Again, all of the injective homomorphisms cancel in the above difference, so $h_G'(x)$ is a sum of
order-$1$ terms corresponding to each \emph{non-injective} homomorphism.
Thus the same reasoning as in the previous paragraph applies and we obtain
\begin{equation}
\label{eq:noninjvar2}
    \Var[h_G'(X)] \lesssim n^{2\sv-4}
\end{equation}
as well.

Applying \eqref{eq:noninjvar1} and \eqref{eq:noninjvar2} in lines \eqref{eq:d33line1} and \eqref{eq:d33line2}
respectively (in the bound for $\Var[\mathsf{remainder}(X)]$)
and using again the fact that $\sigma_n^2 = \Theta(n^2)$, we obtain
\begin{equation}
\label{eq:varremainder}
    \Var[\mathsf{remainder}(X)] \lesssim n^{-2}.
\end{equation}
Plugging this into \eqref{eq:d3varbound1}, we finally arrive at the bound
\begin{equation}
    \Var[\Delta_3(X)] \lesssim n^{-1+\eps}
\end{equation}
for any $\eps > 0$.
This means that $\delta_3 = \sqrt{\Var[\Delta_3(X)]} \lesssim n^{-\fr{1}{2} + \eps}$ for any $\eps > 0$.

\subsubsection{Conclusion of the proof}
\label{sec:stein_global_conclusion}

From the previous six sections, for any $\eps > 0$ we have
\begin{equation}
    b \gtrsim 1, \qquad
    \delta_0 \lesssim n^{-\fr{1}{2} + \eps}, \qquad
    \delta_1 \lesssim n^{-1}, \qquad
    \delta_1' \lesssim n^{-\fr{1}{2} + \eps}, \qquad
    \delta_2 \lesssim n^{-1}, \qquad
    \delta_3 \lesssim n^{-\fr{1}{2} + \eps}.
\end{equation}
Therefore, Theorem \ref{thm:stein} immediately implies
Theorem \ref{thm:global_clt_intro}/\ref{thm:global_clt}.
\subsection{Quantitative CLT for vertex degree}
\label{sec:stein_local}

In this final section we turn to the proof of Theorem \ref{thm:local_clt_intro}, whose statement has
been reproduced below for the reader's convenience.
This is a quantitative central limit theorem for the degree $\deg_v(X)$ of any deterministic vertex $v \in [n]$
of the $n$-vertex random graph $X$.
As before, we take $X$ to have the distribution \eqref{eq:tilt} with $Y \sim \cG(n,p)$
for $p \in U_\beta$, and with $\sg$ given by \eqref{eq:smallgdef}.
We again make the ``phase-uniqueness or forest'' assumption \ref{assumption}, which is only used when applying
Proposition \ref{prop:Eproduct}, including its implicit use
via Proposition \ref{prop:marginal} and Lemma \ref{lem:local_step_concentration}.

Recall the notation that $G_j = (\cV_j, \cE_j)$ for $j = 0, \dotsc, K$ are the graphs in the ERGM specification
\eqref{eq:hdef}, and moreover that $\sv_j = |\cV_j|$ and $\se_j = |\cE_j|$.
Additionally, for $\rho \in \cV_j$, we denote by $\sd_\rho$ the degree of $\rho$ in $G_j$.
Again, since $G_0$ is always a single edge, the $j=0$ term naturally drops out of many of the sums.

\begin{theorem}[Restatement of Theorem \ref{thm:local_clt_intro}]
\label{thm:local_clt}
Let us make Assumption \ref{assumption}.
As in \eqref{eq:sigmadef_local_intro}, set
\begin{equation}
\label{eq:sigmadef_local}
    \varsigma_n^2 = \left( 1 - p(1-p) \sum_{j=1}^K \beta_j p^{\se_j-2} \sum_{\rho \in \cV_j} \sd_\rho (\sd_\rho-1) \right)^{-1} \times p (1-p) (n-1).
\end{equation}
Then for any deterministic vertex $v \in [n]$ and any $\eps > 0$, if $X$ is distributed as \eqref{eq:tilt}
with $Y \sim \cG(n,p)$ and $\sg$ given in \eqref{eq:smallgdef}, under the ``phase-uniqueness or forest'' assumption \ref{assumption}, we have
\begin{equation}
    \d \left( \fr{\deg_v(X) - \E[\deg_v(X)]}{\varsigma_n}, Z \right) \lesssim n^{-\fr{1}{4}+\eps},
\end{equation}
for both $\d = \dW$ and $\d = \dK$.
Recalling Remark \ref{rmk:equivalence}, this is equivalent to the original statement of
Theorem \ref{thm:local_clt_intro} which was stated for $X \sim \condmeas$.
\end{theorem}

As in the proof of Theorem \ref{thm:global_clt_intro}/\ref{thm:global_clt}, the reason for
the form of $\varsigma_n^2$ will become clear in Section \ref{sec:stein_local_delta3}, before which
we will only use the fact that $\varsigma_n = \Theta(\sqrt{n})$.
Again, however, we should first check that $\varsigma_n$ is well-defined, i.e.\ that
\begin{equation}
\label{eq:varsigmawelldef}
    p(1-p) \sum_{j=1}^K \beta_j p^{\se_j-2} \sum_{\rho \in \cV_j} \sd_\rho (\sd_\rho-1) < 1.
\end{equation}
This follows from the fact that $\sigma_n^2$ from \eqref{eq:sigmadef_global_intro}
and \eqref{eq:sigmadef_global} is well-defined, as shown below the statement of Theorem \ref{thm:global_clt}.
Indeed, we have
\begin{equation}
    \sum_{\rho \in \cV_j} \sd_\rho (\sd_\rho - 1) \leq \sum_{\rho \in \cV_j} \sd_\rho (\se_j - 1)
    = 2 \se_j (\se_j-1)
\end{equation}
by the handshake lemma, which yields \eqref{eq:varsigmawelldef}, recalling 
the fact \eqref{eq:sigmawelldef} that $\sigma_n^2$ is well-defined.

As in the proof of Theorem \ref{thm:global_clt_intro}/\ref{thm:global_clt}, we will apply
Theorem \ref{thm:stein} to the choice
\begin{equation}
    f(X) = \fr{\deg_v(X) - \E[\deg_v(X)]}{\varsigma_n},
\end{equation}
meaning we must bound $b$ given in \eqref{eq:bdef} from below, and the error terms
$\delta_0, \delta_1, \delta_1', \delta_2, \delta_3$ given in \eqref{eq:delta0def}-\eqref{eq:delta3def}
from above.

\subsubsection{Preliminary calculations}
\label{sec:stein_local_prelim}

Again we begin by calculating $b$ given in \eqref{eq:bdef} as well as $\Delta_{1,e}(x)$
and $\Delta_{2,e}(x)$ given by \eqref{eq:Delta1def} and \eqref{eq:Delta2def}.
As before, the ordering on edges $e \in \edgeset$ does not matter so we just fix an ordering
arbitrarily and let $e-1$ denote the edge preceding $e$.
Similarly to the previous case, since $\deg_v(x) = \sum_{e \ni v} x(e)$, for any
$e \ni v$ we have
\begin{equation}
\label{eq:differenceslocal}
    f(x) - f(x^{(e)}) = f(x^{[e]}) - f(x^{[e-1]}) = \fr{x(e) - Y(e)}{\varsigma_n},
\end{equation}
and these differences are $0$ for all $e \not \ni v$;
recall the definitions of $x^{(e)}$ and $x^{[e]}$ from
\eqref{eq:xparendef} and \eqref{eq:xbrackdef} respectively.
The above equality will simplify our calculations as in the previous case.

Let us now calculate $\Delta_{1,e}(x)$, recalling the definition from \eqref{eq:Delta1def}.
First, if $v \notin e$ then $\Delta_{1,e}(x) = 0$; on the other hand, if $v \in e$ then
\begin{equation}
    \Delta_{1,e}(x) = \fr{1}{2\varsigma_n^2} ((1-2p) x(e) + p)
\end{equation}
as in the calculation \eqref{eq:squarecalculation}.
This gives
\begin{align}
\label{eq:bcalcloc}
    b &= \fr{1}{2\varsigma_n^2} \sum_{e \ni v} \E[(1-2p)X(e) + p] \\
    &= \fr{(n-1) p (1-p)}{\varsigma_n^2} + O(n^{-\fr{1}{2}} \sqrt{\log n})
\end{align}
by Proposition \ref{prop:marginal}.
As in the previous proof, this is order-one.

Next, if $v \notin e$ then we also have $\Delta_{2,e}(x) = 0$;
if $v \in e$ then by the same calculations as in
Section \ref{sec:stein_global_prelim},
\begin{align}
\label{eq:Delta2calcloc}
    \Delta_{2,e}(x) &= \sum_{j=1}^K \fr{\beta_j}{2 n^{\sv_j-2} \varsigma_n} ((1-2p) x(e) + p)
    \left( \sN_{G_j}(x,e) - 2 \se_j p^{\se_j-1} n^{\sv_j-2} \right) \\
    &\quad - \fr{1}{\varsigma_n} \E \left[(\sR(x) - \sR(x^{(e)}))(x(e) - Y(e))\right].
\end{align}
Now we turn to bounding the error terms one-by-one

\subsubsection{Bounding $\delta_0$}
\label{sec:stein_local_delta0}

We remind the reader that $|\sg(x) - \sg(x^{(e)})|$ is $O(1)$ uniformly in $x$ by the calculation in
Section \ref{sec:stein_global_delta0}.
Now, since $\varsigma_n = \Theta(\sqrt{n})$, by \eqref{eq:differenceslocal} we have
$|f(x) - f(x^{(e)})|, |f(x^{[e]}) - f(x^{[e-1]})| \lesssim n^{-\fr{1}{2}}$ if $v \in e$,
and these expressions are zero otherwise.
Therefore, by a similar calculation as in Section \ref{sec:stein_global_delta0}, recalling the
definition of $\delta_0$ from \eqref{eq:delta0def}, we find that
\begin{equation}
    \delta_0 \lesssim
    n^{-\fr{1}{2}} \sum_{e \ni v} \E \left[ \left| n^2 \partial_e \sH(X) - 2 \sH'(p) \right|^3 \right]
    + n^{-1} \sum_{e \ni v} \E \left[ \left| n^2 \partial_e \sH(X) - 2 \sH'(p) \right|^2 \right]
    + e^{-\Omega(n^2)}.
\end{equation}
Now applying the concentration inequality \eqref{eq:hconc} again, we find that
$\delta_0 \lesssim n^{-1 + \eps}$ for any $\eps > 0$.

Again, if the concentration inequality \eqref{eq:hconc} could be extended to the non-forest phase coexistence
case, then the strategy used here would be able to prove a CLT in that case, although with a non-explicit
variance proxy $\varsigma_n^2$.

\subsubsection{Bounding $\delta_1$}
\label{sec:stein_local_delta1}

Recalling the definition of $\delta_1$ from \eqref{eq:delta1def} and using again the fact
that $|f(x) - f(x^{(e)})|, |f(x^{[e]}) - f(x^{[e-1]})| \lesssim n^{-\fr{1}{2}}$ if $v \in e$
and these differences are zero otherwise, we obtain

\begin{equation}
    \delta_1 = \sum_{e \ni v} \E \left[
        \left( f(X) - f(X^{(e)}) \right)^2
        \left| f(X^{[e]}) - f(X^{[e-1]}) \right|
    \right]
    \lesssim n \cdot n^{-\fr{3}{2}} \lesssim n^{-\fr{1}{2}}.
\end{equation}

\subsubsection{Bounding $\delta_1'$}
\label{sec:stein_local_delta1prime}

We take $D_e^*(x,y) = \fr{1}{\varsigma_n}$ like before, but this time $\varsigma_n = \Theta(\sqrt{n})$.
Recalling the definition of $\delta_1'$ from \eqref{eq:delta1primedef} and using again the fact
that $e^{|\sg(x) - \sg(x^{(e)})|} = O(1)$, the first term in $\delta_1'$ is bounded, up to constants, by
\begin{equation}
    \fr{1}{n} \sum_{e \ni v} \E \left[ \left| n^2 \partial_e \sH(X) - 2 \sH'(p) \right| \right]
    + e^{-\Omega(n^2)},
\end{equation}
using Lemma \ref{lem:Rfacts}.
This is $\lesssim n^{-\fr{1}{2} + \eps}$ for any $\eps > 0$ by \eqref{eq:hconc}.
The second term in $\delta_1'$ is $\lesssim n^{-\fr{1}{2}} \sqrt{\log n}$ by the same reasoning as
in Section \ref{sec:stein_global_delta1prime}.
So $\delta_1' \lesssim n^{-\fr{1}{2} + \eps}$ for any $\eps > 0$.

\subsubsection{Bounding $\delta_2$}
\label{sec:stein_local_delta2}

Recalling the definition of $\delta_2$ from \eqref{eq:delta2def} and the calculation \eqref{eq:Delta2calcloc},
we have
\begin{equation}
    \delta_2^2 = \Var \left[ \sum_{e \in \edgeset} \Delta_{1,e}(X) \right] = \fr{(1-2p)^2}{4 \varsigma_n^4} \Var[\deg_v(X)].
\end{equation}
Now, since $\varsigma_n^4 = \Theta(n^2)$ and $\Var[\deg_v(X)] \lesssim n$ by the concentration for \emph{local}
Lipschitz functions afforded by Theorem \ref{thm:lipschitzconcentration} (recalling also Remark \ref{rmk:equivalence}),
we obtain $\delta_2 \lesssim n^{-\fr{1}{2}}$.

\subsubsection{Bounding $\delta_3$}
\label{sec:stein_local_delta3}

As in the proof of Theorem \ref{thm:global_clt_intro}/\ref{thm:global_clt}, bounding $\delta_3$ is the most
technical portion of the proof.
The strategy is broadly similar here as in the previous case, but there are more algebraic manipulations
required.
As before, let us define
\begin{equation}
    \Delta_3(x) \coloneqq \sum_{e \in \edgeset} \Delta_{2,e}(x) - (1-b) f(x),
\end{equation}
so that $\delta_3 = \sqrt{\Var[\Delta_3(X)]}$, recalling the definition from \eqref{eq:delta3def}.
Then, using \eqref{eq:bcalcloc} and \eqref{eq:Delta2calcloc}, we have
\begin{align}
    \Delta_3(x) &= \sum_{j=1}^K \fr{\beta_j}{2 n^{\sv_j-2} \varsigma_n}
    \sum_{e \ni v} ((1-2p)x(e) + p) (\sN_{G_j}(x,e) - 2 \se_j p^{\se_j-1} n^{\sv_j-2}) \\
    &\qquad - \left(1 - \fr{(n-1)p(1-p)}{\varsigma_n^2} + O(n^{-\fr{1}{2}} \sqrt{\log n}) \right) \fr{\deg_v(x) - \E[\deg_v(X)]}{\varsigma_n} \\
    &\qquad - \fr{1}{\varsigma_n} \E \left[ \sum_{e \ni v} (\sR(x) - \sR(x^{(e)})) (x(e) - Y(e)) \right]
\end{align}
Rearranging terms yields
\begin{align}
    \label{eq:d3locline1}
    \Delta_3(x) &= \sum_{j=1}^K \fr{\beta_j (1-2p)}{2 n^{\sv_j-2} \varsigma_n} \left( \sum_{e \ni v} x(e) \sN_{G_j}(x,e) - 2 \se_j p^{\se_j-1} n^{\sv_j-2} \deg_v(x) \right) \\
    \label{eq:d3locline2}
    &\qquad + \sum_{j=1}^K \fr{\beta_j p}{2 n^{\sv_j-2} \varsigma_n} \left( \sum_{e \ni v} \sN_{G_j}(x,e) - 2 \se_j p^{\se_j-1} n^{\sv_j-2} (n-1) \right) \\
    &\qquad - \left(1 - \fr{(n-1)p(1-p)}{\varsigma_n^2} + O(n^{-\fr{1}{2}} \sqrt{\log n}) \right) \fr{\deg_v(x) - \E[\deg_v(X)]}{\varsigma_n} \\
    &\qquad - \fr{1}{\varsigma_n} \E \left[ \sum_{e \ni v} (\sR(x) - \sR(x^{(e)})) (x(e) - Y(e)) \right].
\end{align}
Now, since every homomorphism $G_j \to x$ using some $e \ni v$ assigns some vertex $\rho \in \cV_j$ to $v$, we have
\begin{align}
    \sum_{e \ni v} x(e) \sN_{G_j}(x,e)
    &= \sum_{\rho \in \cV_j} \sum_{e \ni v}  x(e) \sN_{G_j}^{\rho \to v}(x,e) + O(n^{\sv_j-2}) \\
    \label{eq:line1rep}
    &= \sum_{\rho \in \cV_j} \sd_\rho \sN_{G_j}^{\rho \to v} (x) + O(n^{\sv_j-2}).
\end{align}
The error term here arises from the homomorphisms which are not injective;
note that each \emph{injective} homomorphism counted by $\sN_{G_j}^{\rho \to v}(x)$ leads to $\sd_\rho$
injective homomorphisms counted by $\sum_{e \ni v} x(e) \sN_{G_j}^{\rho \to v}(x,e)$, since there are
$\sd_\rho$ edges $e \ni v$ which each contribute once for this homomorphism.
Further, there are $\lesssim n^{\sv_j-2}$ non-injective homomorphisms overall which map some vertex to $v$.
Similarly, we have
\begin{equation}
    \label{eq:line2rep}
    \sum_{e \ni v} \sN_{G_j}(x,e) = \sum_{\rho \in \cV_j} \sum_{\phi \ni \rho} \sN_{G_j \setminus \phi}^{\rho \to v} (x)
    + O(n^{\sv_j-2}),
\end{equation}
since each injective homomorphism of $G_j$ in $x$ mapping $\phi$ to $e \ni v$ in $x^{+e}$ restricts to an
injective homomorphism of $G_j \setminus \phi$ in $x$ mapping some vertex of $\phi$ to $v$.
All this is to say that we may replace the first terms inside the parentheses on lines \eqref{eq:d3locline1}
and \eqref{eq:d3locline2} with the corresponding terms given in \eqref{eq:line1rep} and \eqref{eq:line2rep}
respectively, incurring an overall error which is bounded (uniformly over $x$) by
$O(\varsigma_n^{-1}) = O(n^{-\fr{1}{2}})$.
Note that in this case, unlike in Section \ref{sec:stein_global_delta3}, we do not need to apply a
concentration inequality to these errors later since they are small enough uniformly in $x$.

Furthermore, by the handshake lemma we can expand $2 \se_j$ as $\sum_{\rho \in \cV_j} \sd_\rho$, and so we obtain
\begin{align}
    \label{eq:d3loc2line1}
    \Delta_3(x) &= \sum_{j=1}^K \fr{\beta_j(1-2p)}{2 n^{\sv_j-2} \varsigma_n} \sum_{\rho \in \cV_j} \sd_\rho
    \left( \sN_{G_j}^{\rho \to v} (x) - p^{\se_j-1} n^{\sv_j-2} \deg_v(x) \right) \\
    \label{eq:d3loc2line2}
    &\qquad + \sum _{j=1}^K \fr{\beta_j p}{2 n^{\sv_j-2} \varsigma_n} \sum_{\rho \in \cV_j} \sum_{\phi \ni \rho}
    \left( \sN_{G_j \setminus \phi}^{\rho \to v} (x) - p^{\se_j-1} n^{\sv_j-2} (n-1) \right) \\
    &\qquad - \left(1 - \fr{(n-1)p(1-p)}{\varsigma_n^2} + O(n^{-\fr{1}{2}} \sqrt{\log n}) \right) \fr{\deg_v(x) - \E[\deg_v(X)]}{\varsigma_n} \\
    &\qquad - \fr{1}{\varsigma_n} \E \left[ \sum_{e \ni v} (\sR(x) - \sR(x^{(e)})) (x(e) - Y(e)) \right] + O(n^{-\fr{1}{2}}).
\end{align}
Now we will again insert the correct H\'ajek projections, replacing the second terms in the parentheses
in lines \eqref{eq:d3loc2line1} and \eqref{eq:d3loc2line2} to obtain
\begin{align}
    \Delta_3(x) &= \sum_{j=1}^K \fr{\beta_j (1-2p)}{2 n^{\sv_j-2} \varsigma_n} \sum_{\rho \in \cV_j} \sd_\rho
    \left( \sN_{G_j}^{\rho \to v} (x) - \sd_\rho p^{\se_j-1} n^{\sv_j-2} \deg_v(x) \right) \\
    &\qquad + \sum_{j=1}^K \fr{\beta_j p}{2 n^{\sv_j-2} \varsigma_n} \sum_{\rho \in \cV_j} \sum_{\phi \ni \rho}
    \left( \sN_{G_j \setminus \phi}^{\rho \to v} (x) - (\sd_\rho-1) p^{\se_j-2} n^{\sv_j-2} \deg_v(x) \right) \\
    &\qquad + \sum_{j=1}^K \fr{\beta_j}{2} \left( (1-2p) p^{\se_j-1} \sum_{\rho \in \cV_j}
    \sd_\rho (\sd_\rho-1) + p \cdot p^{\se_j-2} \sum_{\rho \in \cV_j} \sd_\rho (\sd_\rho-1) \right) \fr{\deg_v(x)}{\varsigma_n} \\
    &\qquad - \sum_{j=1}^K \fr{\beta_j \se_j p^{\se_j}}{\varsigma_n} (n-1) \\
    &\qquad - \left(1 - \fr{(n-1)p(1-p)}{\varsigma_n^2} + O(n^{-\fr{1}{2}} \sqrt{\log n}) \right) \fr{\deg_v(x) - \E[\deg_v(X)]}{\varsigma_n} \\
    &\qquad - \fr{1}{\varsigma_n} \E \left[ \sum_{e \ni v} (\sR(x) - \sR(x^{(e)})) (x(e) - Y(e)) \right] + O(n^{-\fr{1}{2}}).
\end{align}
Now again let us combine terms in a convenient fashion to obtain
\begin{align}
    \hspace{-1cm}
    \Delta_3(x) &= \sum_{j=1}^K \fr{\beta_j (1-2p)}{2 n^{\sv_j-2} \varsigma_n} \sum_{\rho \in \cV_j} \sd_\rho \hsN_{G_j}^{\rho \to v}(x) \\
    &\qquad + \sum_{j=1}^K \fr{\beta_j p}{2 n^{\sv_j-2} \varsigma_n} \sum_{\rho \in \cV_j} \sum_{\phi \ni \rho} \hsN_{G_j \setminus \phi}^{\rho \to v}(x) \\
    \label{eq:varsigmareason}
    &\qquad + \left( p (1-p) \sum_{j=1}^K \beta_j p^{\se_j -2} \sum_{\rho \in \cV_j} \sd_\rho (\sd_\rho-1) - 1 + \fr{(n-1)p(1-p)}{\varsigma_n^2} + O(n^{-\fr{1}{2}} \sqrt{\log n})\right) \fr{\deg_v(x)}{\varsigma_n} \\
    \label{eq:localdeterministic}
    &\qquad - \sum_{j=1}^K \fr{\beta_j \se_j p^{\se_j}}{\varsigma_n} (n-1) + \fr{\E[\deg_v(X)]}{\varsigma_n}
    \left(1 - \fr{(n-1)p(1-p)}{\varsigma_n^2} + O(n^{-\fr{1}{2}} \sqrt{\log n}) \right) \\
    &\qquad - \fr{1}{\varsigma_n} \E \left[ \sum_{e \ni v} (\sR(x) - \sR(x^{(e)})) (x(e) - Y(e)) \right] + O(n^{-\fr{1}{2}}).
\end{align}
The line \eqref{eq:varsigmareason} is
$O(n^{-\fr{1}{2}} \sqrt{\log n}) \cdot \fr{\deg_v(x)}{\varsigma_n} = O(n^{-1} \sqrt{\log n}) \cdot \deg_v(x)$
by definition of $\varsigma_n^2$ (as before, this is in fact the reason for this precise definition of $\varsigma_n$).
Since the line \eqref{eq:localdeterministic} is deterministic, we obtain
\begin{align}
    \label{eq:d3loc3line1}
    \fr{\Var[\Delta_3(X)]}{\sum_{j=1}^K \sv_j + \sum_{j=1}^K \sum_{\rho \in \cV_j} \sd_\rho + 3} &\leq
    \sum_{j=1}^K \fr{\beta_j^2 (1-2p)^2}{4 n^{2\sv_j-4} \varsigma_n^2} \sum_{\rho \in \cV_j} \sd_\rho^2 \Var \left[ \hsN_{G_j}^{\rho \to v} (X) \right] \\
    \label{eq:d3loc3line2}
    &\qquad + \sum_{j=1}^K \fr{\beta_j^2 p^2}{4 n^{2\sv_j-4} \varsigma_n^2} \sum_{\rho \in \cV_j} \sum_{\phi \ni \rho} \Var \left[ \hsN_{G_j \setminus \phi}^{\rho \to v} (X) \right] \\
    &\qquad + O(n^{-2} \log n) \cdot \Var[\deg_v(X)] \\
    &\qquad + \fr{1}{\varsigma_n^2} \Var \left[ \Econd{\sum_{e \ni v} (\sR(X) - \sR(X^{(e)})) (X(e) - Y(e))}{X} \right] + O(n^{-1}).
\end{align}
Finally, by Proposition \ref{prop:local_hajek_intro} (recalling Remark \ref{rmk:equivalence}),
for any $\eps > 0$ and any finite graph $G$ with $\sv$ vertices
(which must be a forest in the phase coexistence case) we have
\begin{equation}
    \Var \left[ \hsN_G^{\rho \to v} (X) \right] \lesssim n^{2\sv-\fr{7}{2}+\eps}.
\end{equation}
Plugging this into lines \eqref{eq:d3loc3line1} and \eqref{eq:d3loc3line2}
and using the fact that $\varsigma_n^2 = \Theta(n)$ and $\Var[\deg_v(X)] \lesssim n$, and that the
term above involving $\sR$ is exponentially small in $n^2$ by Lemma \ref{lem:Rfacts}, we obtain
\begin{equation}
    \Var[\Delta_3(X)] \lesssim n^{-\fr{1}{2}+\eps} + n^{-1} \log n + \Exp{-\Omega(n^2)} + n^{-1},
\end{equation}
meaning that $\delta_3 = \sqrt{\Var[\Delta_3(X)]}\lesssim n^{-\fr{1}{4} + \eps}$ for any $\eps > 0$.

\subsubsection{Conclusion of the proof}
\label{sec:stein_local_conclusion}

From the previous six sections, for any $\eps > 0$ we have
\begin{equation}
    b \gtrsim 1, \qquad
    \delta_0 \lesssim n^{-1 + \eps}, \qquad
    \delta_1 \lesssim n^{-\fr{1}{2}}, \qquad
    \delta_1' \lesssim n^{-\fr{1}{2} + \eps}, \qquad
    \delta_2 \lesssim n^{-\fr{1}{2}}, \qquad
    \delta_3 \lesssim n^{-\fr{1}{4} + \eps}.
\end{equation}
Therefore, Theorem \ref{thm:stein} implies
Theorem \ref{thm:local_clt_intro}/\ref{thm:local_clt}.

\bibliography{references}
\bibliographystyle{alpha}

\end{document}